\documentclass[11pt]{article}
\input xy
\xyoption{all}
\usepackage{amsmath, amsthm, amssymb}
\usepackage{showlabels}
\usepackage{xcolor}
\usepackage{titlesec}

\usepackage{hyperref}
\usepackage[titles]{tocloft}
\usepackage{enumitem}

\usepackage{etoolbox}
\makeatletter
\patchcmd{\l@section}  {#1}  {\S#1}  {}{}\makeatother

\usepackage{enumitem}
\newlist{abbrv}{itemize}{1}
\setlist[abbrv,1]{label=,labelwidth=1in,align=parleft,itemsep=0.1\baselineskip,leftmargin=!}

 \numberwithin{equation}{section}

\usepackage[nottoc,notlot,notlof]{tocbibind}

\titleformat{\section}
  {\normalfont\Large\bfseries}{\S\thesection}{1em}{}
  
  \titleformat{\subsection}
{\normalfont\large\bfseries}{\S \thesubsection}{1em}{}

\usepackage{dutchcal}

\def\sqr#1#2{{\vcenter{\vbox{\hrule height .#2pt
                             \hbox{\vrule width .#2pt height#1pt\kern#1pt
                                   \vrule width .#2pt}
                             \hrule height .#2pt}}}}

\newtheorem{theorem}[equation] {Theorem}

\theoremstyle{definition}
\newtheorem{definition}[equation]{Definition}

\theoremstyle{definition}
\newtheorem{remark} [equation] {Remark}

\theoremstyle{recall}
\newtheorem{recall} [equation] {Recall}

\theoremstyle{theorem}
\newtheorem{example}[equation]{Example}

\theoremstyle {proposition}
\newtheorem{proposition}[equation]{Proposition}

\theoremstyle {proposition}
\newtheorem{corollary}[equation]{Corollary}

\theoremstyle {defintion}
\newtheorem{lemma}[equation]{Lemma}

\theoremstyle{definition}
\newtheorem{notation}[equation]{Notation}

\def\tic{\star}
\def\tac{\star\star}

\def\uno{1\!\!\!\!1}
\def\goths{\frak sl_2}
\def\gothd{\hat{{\frak sl_3}}^{\!\!\chi}}
\def\u{{\cal U}}
\def\uz{\u_{\Z}}
\def\hu{\hat{\cal U}}
\def\huz{\hat{\cal U}_{\Z}}
\def\tu{\tilde{\cal U}}
\def\tuz{\tilde{\cal U}_{\Z}}
\def\duz{^*{\cal U}_{\Z}}
\def\um{\tilde{\cal U}_{\Z,M}}

\def\gothg{\frak g}

\def\n{\noindent}
\def\k{\mathcal{h}}

\def\compo{\,{\scriptstyle\circ}\,}

\def\N{{\mathbb{N}}
}
\def\Z{{\mathbb{Z}}
}
\def\Q{{\mathbb{Q}}}

\def\C{{\mathbb{C}}
}

\begin{document} 
\title{On the integral form
of rank 1 Kac-Moody algebras.
}

\date{}
\author{Ilaria Damiani, Margherita Paolini}
\maketitle

\begin{abstract}

\n In this paper we shall prove that the $\Z$-subalgebra generated by the divided powers of the Drinfeld generators $x_r^{\pm}$ ($r\in\Z$) of the Kac-Moody algebra of type $A_2^{(2)}$ is an integral form (strictly smaller than Mitzman's, see \cite{DM}) of the enveloping algebra, we shall exhibit a basis generalizing the one provided in \cite{HG} for the untwisted affine Kac-Moody algebras 
and we shall determine explicitly the commutation relations. Moreover we prove that both in the untwisted and in the twisted case the positive (respectively negative) imaginary part of the integral form is an algebra of polynomials over $\Z$.

\end{abstract}

\small \small \tableofcontents 

\section{Introduction} \label{intr}
\vskip .5truecm

\n Recall that the twisted affine Kac-Moody algebra of type $A_2^{(2)}$ is $\gothd$, the $\chi$-invariant subalgebra of $\hat{{\frak sl_3}}$ where $\chi$ is the non trivial Dynkin diagram automorphism of $A_2$ (see \cite{VK}) and denote by $\tu$ its enveloping algebra $\u(\gothd)$.

\n The aim of this paper is to give a basis over $\Z$ of the $\Z$-subalgebra of $\tu$ 
 generated by the divided powers of the Drinfeld generators $x_r^{\pm}$'s ($r\in\Z$) (see definitions \ref{a22} and \ref{thuz}),
thus proving that this $\Z$-subalgebra is an integral form of $\tu$.
\vskip .3 truecm

\n The integral forms for finite dimensional semisimple Lie algebras were first introduced by Chevalley in \cite{Ch} for the study of the Chevalley groups and of their representation theory. 

\n The construction of the ``divided power''-$\Z$-form 
 for 
 the simple finite dimensional Lie algebras  
 is due to  Kostant (see \cite{Ko})
 ; it has been generalized to the untwisted affine Kac-Moody algebras by Garland in \cite{HG} as we shall quickly recall.

\n Given a simple Lie algebra ${{\frak g}}_0$ and the corresponding untwisted affine Kac-Moody algebra ${{\frak g}}={{\frak g}}_0\otimes\C[t,t^{-1}]\oplus\C c$ provided with an (ordered) Chevalley basis, the $\Z$-subalgebra $\uz$ of $\u=\u({{\frak g}})$ generated by the divided powers of the  
real root vectors is an integral form of $\u$; a $\Z$-basis of this integral form (hence its $\Z$-module structure) can be described by decomposing $\uz$ as tensor product of its $\Z$-subalgebras relative respectively to the real root vectors ($\uz^{re,+}$ and $\uz^{re,-}$), to the imaginary root vectors ($\uz^{im,+}$ and $\uz^{im,-}$) and to the Cartan subalgebra ($\uz^{{\frak h}}$):

$\uz^{re,+}$ has a basis $B^{re,+}$ consisting of the (finite) ordered products of divided powers of the distinct positive real root vectors 
and $(\uz^{re,-},B^{re,-})$ can be described in the same way: $$B^{re,\pm}=\{x_{\pm\beta_1}^{(k_{\beta_1})}\cdot ...\cdot x_{\pm\beta_N}^{(k_{\beta_N})}|N\geq 0,\ \beta_1>...>\beta_N>0\ {\rm{real\ roots}},\ k_{\beta_j}> 0\ \forall j\}.$$
Here a real root $\beta$ of ${{\frak g}}$ is said to be positive if there exists a positive root $\alpha$ of ${{\frak g}}_0$ such that $\beta=\alpha$ or $\beta-\alpha$ is imaginary; $x_{\beta}$ is the Chevalley generator corresponding to the real root $\beta$.

A basis $B^{\frak h}$ of $\uz^{{\frak h}}$, which is commutative, consists of the products  
of the ``binomials'' 
of the (Chevalley) generators $h_{i
}$ ($i\in I$) of the Cartan subalgebra of ${{\frak g}}$: 
$$B^{\frak h}=\left\{
\prod_i{h_{i
}\choose k_i}\Big|k_i\geq 0\ \forall i\right\};$$
it is worth remarking that $\uz^{{\frak h}}$ is not an algebra of polynomials.

$\uz^{im,+}$ (and its symmetric $\uz^{im,-}$) is commutative, too; as $\Z$-module it is isomorphic to the tensor product of the $\u_{i,\Z}^{im,+}$'s (each factor corresponding to the $i^{th}$ copy of $\u(\hat\goths)$ inside $\u$), 
so that it is enough to describe it in the rank 1 case: the basis $B^{im,+}$ of $\uz^{im,+}(\hat\goths)$ provided by Garland can be described as a set of finite products of the elements 
$\Lambda_k{(\xi(m))}$ ($r\in\N$, $m>0$), where the $\Lambda_k{(\xi(m))}$'s are the elements of $\u^{im,+}=\C[h_{
r}(=h
\otimes t^r)|r>0]$ defined recursively (for all $m\neq 0$) by 
$$\Lambda_{-1}{(\xi(m))}=1,\ \ k\Lambda_{k-1}{(\xi(m))} =\sum_{r\geq 0,s>0\atop r+s=k}\Lambda_{r-1}{(\xi(m))} {h_{ms}}:$$
$$B
^{im,+}=\left\{\prod_{m>0}\Lambda_{k_m-1}{(\xi(m))}|k_m\geq 0\ \forall m,\ \#\{m>0|k_m\neq 0\}<\infty\right\};$$

\n It is not clear from this description that $\uz^{im,+}$ and $\uz^{im,-}$ are algebras of polynomials.

\n Thanks to the isomorphism of $\Z$-modules
$$\uz\cong\uz^{re,-}\otimes_{\Z}\uz^{im,-}\otimes_{\Z}\uz^{{\frak h}}\otimes_{\Z}\uz^{im,+}\otimes_{\Z}\uz^{re,+}$$
a $\Z$-basis $B$ of $\uz$ is produced as multiplication of $\Z$-bases of these subalgebras:
$$B=B^{re,-}B^{im,-}B^{{\frak h}}B^{im,+}B^{re,+}.$$

\n The same result has been proved for all the twisted affine Kac-Moody algebras by Mitzman in \cite{DM}, where the author provides a deeper comprehension and a compact description of the commutation formulas by means of a drastic  simplification of both the relations and their proofs. This goal is achieved remarking that the generating series of the elements involved in the basis  
can be expressed as suitable exponentials, observation that allows to apply very general tools of calculus, such as the well known properties
$$x\exp(y)=\exp(y)\exp([\cdot,y])(x)$$
if $\exp(y)$ and $\exp([\cdot,y])(x)$ are well defined,
and 
$$D(\exp(f))=D(f)\exp(f)$$ if $D$ is a derivation such that $[D(f),f]=0$.

\n Here, too, it is not yet clear that $\uz^{im,\pm}$ are algebras of polynomials.

\n However this property, namely  $\uz^{im,+}=\Z[\Lambda_{k-1}=\Lambda_{k-1}{(\xi(1))}=p_{k,1}|k>0]$, is stated in Fisher-Vasta's PhD thesis (\cite{F}), where the author describes the results of Garland for the untwisted case and of Mitzman for $A_2^{(2)}$ aiming at a better understanding of the commutation formulas.
Yet the proof is missing: the theorem describing the integral form is based on observations which seem to forget some necessary commutations, those between $(x_r^+)^{(k)}$ and $(x_s^-)^{(l)}$ when $|r+s|>1$; in \cite{F} only the cases $r+s=0$ and $r+s=\pm 1$ are considered, the former producing the binomials appearing in $B^{{\frak h}}$, the latter producing the elements $p_{n,1}$
(and their corresponding negative elements in $\uz^{im,-}$).

\vskip .3 truecm

\n Comparing the Kac-Moody presentation of the affine Kac-Moody algebras with its ``Drinfeld'' presentation as current algebra, one can notice a difference between the untwisted and twisted case, which is at the origin of our work. As in the simple finite dimensional case, also in the affine cases the generators of $\uz$ described above are redundant: the $\Z$-subalgebra of $\u$ generated by $\{e_i^{(k)},f_i^{(k)}|i\in I,\ k\in\N\}$, obviously contained in $\uz$, is actually equal to $\uz$. 

\n On the other hand, the situation changes when we move to the Drinfeld presentation and study the $\Z$-subalgebra $\duz$ of $\u$ generated by the divided powers of the Drinfeld generators $
(x_{i,r}^{\pm})^{(k)}$: indeed, while in the untwisted case it is still true that $\uz=$ $\duz$ and (also in the twisted case) it is always true that $\duz\subseteq\uz$, in general we get two different $\Z$-subalgebras of $\u$; more precisely $\duz\subsetneq\uz$ in case $A_{2n}^{(2)}$, that is when there exists a vertex $i$ whose corresponding rank 1 subalgebra 
is not a copy of $\u(\hat\goths)$ but is a copy of $\u(\gothd)$.
Thus in order to complete the description of $\duz$ we need to study the case of $A_2^{(2)}$.
\vskip .3 truecm
\n In the present paper we prove that the $\Z$-subalgebra generated by $$\{(x_r^+)^{(k)},(x_r^-)^{(k)}|r\in\Z,k\in\N\}$$ is an integral form of the enveloping algebra also in the case of $A_2^{(2)}$, we exhibit a basis generalizing the one provided in \cite{HG} and in \cite{DM} and determine the commutation relations in a compact yet explicit formulation (see theorem \ref{trmA22} and appendix \ref{appendA}). We use the same approach as Mitzman's, with a further simplification consisting in the remark that an element of the form $G(u,v)=\exp(xu)\exp(yv)$ is characterized by two properties: $G(0,v)=\exp(yv)$ and ${dG\over du}=xG$.

\n Moreover, studying the rank 1 cases we prove that, both in the untwisted and in the twisted case, $\uz^{im,+}$ and $\duz^{im,+}$ are algebras of polynomials: as stated in \cite{F}, the generators of $\uz^{im,+}$ are the elements $\Lambda_k$ introduced in \cite{HG} and \cite{DM} (see proposition \ref{tmom} and remark \ref{tmfv}); the generators of $\duz^{im,+}$ in the case $A_2^{(2)}$ are elements defined formally as the $\Lambda_k$'s after a deformation of the $h_r$'s (see definition \ref{thuz} and remark \ref{hdiversi}): describing $\duz^{im,+}(\gothd)$ (denoted by $\tuz^{0,+}$) has been the hard part of this work.
\vskip .3 truecm
\n We work over $\Q$ and dedicate a preliminary particular care to the description of some integral forms of $\Q[x_i|i\in I]$ and of their properties and relations when they appear in some non commutative situations, properties that will be repeatedly used for the computations in ${\frak g}$:  fixing the notations helps to understand the construction in the correct setting. With analogous care we discuss the symmetries arising both in $\hat\goths$ and in $\gothd$. We chose to recall also the case of $\goths$ and to give in a few lines the proof of the theorem describing its divided power integral form in order to present in this easy context the tools that will be used in the more complicated affine cases.

\vskip .3 truecm

\n The paper is organized as follows.

Section \ref{intgpl} is devoted to review the description of some integral forms of the algebra of polynomials (polynomials over $\Z$, divided powers,``binomials'' and symmetric functions, see \cite{IM}): they are introduced together with their generating series as exponentials of suitable series with null constant term, and their properties are rigorously stated, thus preparing to their use in the Lie algebra setting.

\n We have inserted here, in proposition \ref{tmom}, a result about the stability of the symmetric functions with integral coefficients under the homomorphism $\lambda_m$ mapping $x_i$ to $x_i^m$ ($m>0$ fixed), which is almost trivial in the symmetric function context; it is a straightforward consequence of this observation that $\uz^{im,+}$ is an algebra of polynomials and so is $^*\uz^{im,+}$ in the twisted case. We also provide a direct, elementary proof of this proposition (see proposition \ref{tdmom}).

In section \ref{ncn} we collect some computations in non commutative situations that we shall systematically refer to in the following sections.

Section \ref{sld} deals with the case of $\goths$: the one-page formulation and proof that we present (see theorem \ref{trdc}) inspire the way we study $\hat\goths$ and $\gothd$, and offer an easy introduction to the strategy followed also in the harder affine cases: decomposing our $\Z$-algebra as a tensor product of commutative subalgebras; describing these commutative structures thanks to the examples introduced in section \ref{intgpl}; and glueing the pieces together applying the results of section \ref{ncn}.

\n Even if the results of this section imply the commutation rules between $(x_r^+)^{(k)}$ and $(x_{-r}^+)^{(l)}$ ($r\in\Z$, $k,l\in\N$) in the enveloping algebra of $\hat\goths$ (see remark \ref{hrs}), it is worth remarking that section \ref{slh} does not depend on section \ref{sld}, and can be read independently (see remark \ref{hev}).

In section \ref{slh} we discuss the case of $\hat\goths$. 

\n The first part of the section is devoted to the choice of the notations in $\hu=\u(\hat\goths)$; to the definition of its (commutative) subalgebras $\hu^{\pm}$ (corresponding to the real component of $\hu$), $\hu^{0,\pm}$ (corresponding to the imaginary component), $\hu^{0,0}$ (corresponding to the Cartan), of their integral forms $\huz^{\pm}$, $\huz^{0,\pm}$, $\huz^{0,0}$, and of the $\Z$-subalgebra $\huz$ of $\hu$; and to a detailed reminder about the useful symmetries (automorphisms, antiautomorphisms, homomorphisms and triangular decomposition) thanks to which we can get rid of redundant computations.

\n In the second part of the section  the apparently
tough computations involved in the commutation relations are reduced to four formulas  whose proofs are contained in a few lines:
proposition \ref{zzk}, proposition \ref{pum}, lemma \ref{limt}, and proposition \ref{exefh}, (together with proposition \ref{tmom}) are all what is needed to show that $\huz$ is an integral form of $\hu$, to recognize that the imaginary (positive and negative) components $\huz^{0,\pm}$ of $\huz$ are the algebras of polynomials $\Z[\Lambda_k(\xi(\pm 1))|k\geq 0]=\Z[\hat h_{\pm k}|k>0]$, and to exhibit a $\Z$-basis of $\huz$ (see theorem \ref{trm}).

In section \ref{ifa22} we finally present the case of $A_2^{(2)}$. 

\n As for $\hat\goths$ we first evidentiate some general structures of $\u(\gothd)$ (that we denote here $\tu$ in order to distinguish it from $\hu=\u(\hat\goths)$): notations, subalgebras and symmetries. Here we introduce the elements $\tilde h_k$ through the announced deformation of the formulas defining the elements $\hat h_k$'s (see definition \ref{thuz} and remark \ref{hdiversi}). We also describe a $\Q[w]$-module structure on a Lie subalgebra $L$ of $\gothd$ (see definitions \ref{sottoalgebraL} and \ref{qwmodulo}), thanks to which we can further simplify the notations.

\n In addition, in remark \ref{emgg} we recall the embeddings of $\hu$ inside $\tu$ thanks to which a big part of the work can be translated from section \ref{slh}.

\n The heart of the problem is thus reduced to the commutation of $\exp(x_0^+u)$ with $\exp(x_1^-v)$ (which is technically more complicated than for $A_1^{(1)}$ since it is a product involving a higher number of factors) and to deducing from this formula the description of the imaginary part of the integral form as the algebra of the polynomials in the $\tilde h_k$'s. To the solution of this problem, which represents the central contribution of this work, we dedicate subsection \ref{sottosezione}, where we concentrate, perform and explain the necessary computations.
\vskip .15 truecm
\n At the end of the paper some appendices are added for the sake of completeness.

In appendix \ref{appendA} we collect all the straightening formulas: since not all of them are necessary to our proofs and in the previous sections we only computed those which were essential for our argument, we give here a complete explicit picture of the commutation relations.

Appendix \ref{appendB} is devoted to the description of a $\Z$-basis of $\Z^{(sym)}[h_r|r>0]$ alternative to that introduced in the example \ref{rvsf}.

\n $\Z^{(sym)}[h_r|r>0]$ is the algebra of polynomials $\Z[\hat h_k|k>0]$, and as such it has a $\Z$-basis consisting of the monomials in the $\hat h_k$'s, which is the one considered in our paper. But, as mentioned above, this algebra, that we are naturally interested in because it is isomorphic to the imaginary positive part of the integral form of the rank 1 Kac-Moody algebras, was not recognized by Garland and Mitzman as an algebra of polynomials: in this appendix the $\Z$-basis they introduce is studied from the point of view of the symmetric functions and thanks to this interpretation it is easily proved to generate freely the same $\Z$-submodule of $\Q[h_r|r>0]$ as the monomials in the $\hat h_k$'s.

In appendix \ref{appendC} we compare the Mitzman integral form of the enveloping algebra of type $A_2^{(2)}$ with the one studied here, proving the inclusion stated above. We also show that our commutation relations imply Mitzman theorem, too.

Finally, in order to help the reader to orientate in the notations and to find easily their definitions, we conclude the paper with an index of symbols, collected in appendix \ref{appendD}.

\vskip .3 truecm
\n The study of the integral form of the affine Kac-Moody algebras from the point of view of the Drinfeld presentation, which differs from the one defined through the Kac-Moody presentation (\cite{HG} and \cite{DM}) in the case $A_{2}^{(2)}$ as outilined above, is motivated by the interest in the representation theory over $\Z$, since for the affine Kac-Moody algebras the notion of highest weight vector with respect to the $e_i$'s has been usefully replaced with that defined through the action of the $x_{i,r}^+$'s (see the works of Chari and Pressley \cite{C} and \cite{CP2}): in order to study what happens over the integers it is useful to work with an integral form defined in terms of the same $x_{i,r}^+$'s.

\n This work is also intended to be the preliminary classical step in the project of constructing and describing the {\it quantum} integral form for the twisted affine quantum algebras (with respect to the Drinfeld presentation). It is a joint project with Vyjayanthi Chari (see also \cite{CP}), 
who proposed it during a period of three months that she passed as a visiting professor at the Department of Mathematics of the University of Rome ``Tor Vergata''.

\n The commutation relations involved are extremely complicated and appear to be unworkable by hands without a deeper insight; we hope that a simplified approach can open a viable way to work in the quantum setting.

\vskip .5 truecm

\section{Integral form and commutative examples} \label{intgpl}

\vskip .5 truecm
\n In this section we give the definition of integral form and summarize, fixing the notations useful to our purpose, some well known commutative examples (deeply studied and systematically exposed in \cite{IM}), which will play a central role in the non commutative enveloping algebra of finite and affine Kac-Moody algebras.

\vskip .3 truecm

\begin{definition} \label{intu} 

\n Let $U$ be a $\Q$-algebra. An integral form of $U$ is a $\Z$-algebra $U_{\Z}$ such that:

i) $U_{\Z}$ is a free $\Z$-module;

ii) $U=\Q\otimes_{\Z}U_{\Z}$.

\n In particular an integral form of $U$ is (can be identified to) a $\Z$-subalgebra of $U$, and a $\Z$-basis of an integral form of $U$ is a $\Q$-basis of $U$.

\end{definition}

\vskip .3 truecm
\begin{example} \label{clpol}

\n Of course $\Z[x_i|i\in I]$ is an integral form of $\Q[x_i|i\in I]$ with basis the set of monomials in the $x_i$'s, namely 
$\{{\bf{x}}^{{\bf{k}}}=\prod_{i\in I}x_i^{k_i}\}$ where ${\bf{k}}:I\to\N$ is finitely supported (that is $\#\{i\in I|k_i\neq 0\}<\infty$). 

\n If $\{y_i\}_{i\in I}$ and $\{x_i\}_{i\in I}$ are $\Z$-bases of the same $\Z$-module, then $\Z[x_i|i\in I]=\Z[y_i|i\in I]$.

\end{example}

This can be said also as follows:

\n {Let $M$ be a free $\Z$-module and $V=\Q\otimes_{\Z}M$} and consider the functor $S=$ ``symmetric algebra'' from the category of $\Z$-modules (respectively $\Q$-vector spaces) to the category of commutative unitary $\Z$-algebras (respectively commutative unitary $\Q$-algebras).

\n Then $SM$ is an integral form of $SV$ and $SM\cap V=M$.

\n By definition, every integral form of $SV$ containing $M$ contains $SM$, that is $SM$ is the least integral form of $SV$ containing $M$.
\vskip .3 truecm
\n We are interested in other remarkable integral forms of $SV$ containing $M$.

\begin{remark} \label{srinv}
\n Let $U$ be a unitary $\Z$-algebra and $f(u)\in U[[u]]$. Then: 

\n 1) If $f(u)\in 1+uU[[u]]$ then:

i) $f(u)$ is invertible in $U[[u]]$;

ii) the coefficients of $f(u)$, those of $f(-u)$ and those of $f(u)^{-1}$ generate the same $\Z$-subalgebra of $U$;

\n 2) If $f(u)\in uU[[u]]$ then ${\rm{exp}}(f(u))$ is a well defined element of $1+uU[[u]]$;

\n 3) If $f(u)\in 1+uU[[u]]$ then ${\rm{ln}}(f(u))$ is a well defined element of $uU[[u]]$;

\n 4) ${\rm{exp}}\compo{\rm{ln}}\big|_{1+uU[[u]]}=id$ and ${\rm{ln}}\compo {\rm{exp}}\big|_{uU[[u]]}=id$.
\end{remark}

\vskip .3 truecm
\begin{notation} \label{ntdvd}
\n Let $a$ be an element of a unitary $\Q$-algebra $U$. The divided powers of $a$ are the elements $$a^{(k)}={a^k\over k!}\, \, (k\in\N).$$
Remark that the generating series of the $a^{(k)}$'s is $\exp(a u)$, that is
\begin{equation}\label{gensexp}
\sum_{k\geq 0}a^{(k)}u^k=\exp(au)
\end{equation}.
\end{notation}

\vskip .3 truecm
\begin{example} \label{dvdpw}

\n Let $\{x_i\}_{i\in I}$ be a $\Z$-basis of $M$. Then it is well known and trivial that:

i) The $\Z$-subalgebra $S^{(div)}M\subseteq SV$ generated by $\{x^{(k)}\}_{x\in M,k\in\N}$ contains $M$; 

ii) $S^{(div)}M\cap V=M$;

iii) $\{x_i^{(k)}\}_{i\in I,k\in\N}$ is a set of algebra-generators (over $\Z$) of $S^{(div)}M$;

iv) the set $\{{\bf{x}}^{({\bf{k}})}=\prod_{i\in I}x_i^{(k_i)}|{\bf{k}}:I\to\N$ is finitely supported$\}$ is a $\Z$-basis of $S^{(div)}M$.

v) $S^{(div)}M$ is an integral form of $SV$ (called the algebra of the divided powers of $M$).

\n $S^{(div)}M$ is also denoted $\Z^{(div)}[x_i|i\in I]$.

\n Remark that if $m(u)=\sum_{r\in\N}m_ru^r\in uM[[u]]$ then
\begin{equation} \label{divpoweq}
m(u)^{(k)}\in S^{(div)}M[[u]]\,\,\forall k\in\N
\end{equation}
or equivalently
$${\rm{exp}}(m(u))\in S^{(div)}M[[u]].$$
The viceversa is obviously also true:
\begin{equation}\label{sdivmv}m(u)\in uV[[u]],\ {\rm{exp}}(m(u))\in S^{(div)}M[[u]]\Leftrightarrow m(u)\in uM[[u]].\end{equation}

\end{example}

\vskip .3 truecm

\begin{notation} \label{ntbin}

\n Let $a$ be an element of a unitary $\Q$-algebra $U$. The ``binomials'' of $a$ are the elements $${a\choose k}={a(a-1)\cdot...\cdot(a-k+1)\over k!}\, \,\,\,\, \,\,\,(k\in\N).$$
Notice that ${a\choose k}$ is the image of ${x\choose k}$ through the evaluation $ev_a:\Q[x]\to U$ mapping $x$ to $a$.
\end{notation}

\n Now consider the series $\exp(a\ln(1+u))$: this is a well defined element of $U[[u]]$ (because $a\ln(1+u)\in uU[[u]]$) whose coefficients are polynomials in $a$; this means that with the notations above 
$$\exp(a\ln(1+u))=ev_a(\exp(x{\rm{ln}}(1+u))).$$
In particular if we want to prove that for all $U$ and for all $a\in U$ the generating series of the ${a\choose k}$'s is $\exp(a\ln(1+u))$
it is enough to prove the claim in the case $a=x\in\Q[x]$, and to this aim it is enough to compare the evaluations on an infinite subset of $\Q$ (for istance on $\N$), thus reducing the proof to the trivial observation that $\forall n\in\N$
$$\sum_{k\in\N}{n\choose k}u^k=(1+u)^n=\exp(n\ln(1+u)).$$
Thus in general the generating series $\exp(a\ln(1+u))$ of the ${a\choose k}$'s can and will be denoted as $(1+u)^a$; more explicitly
\begin{equation}\label{gensbin}\sum_{k\in\N}{a\choose k}u^k=(1+u)^a=\exp\Big(\sum_{r>0}(-1)^{r-1}{a\over r}u^r\Big).\end{equation}

\n It is clear from the definition of $(1+u)^a$ that if $a$ and $b$ are commuting elements of $U$ then $$(1+u)^{a+b}=(1+u)^a(1+u)^b.$$
It is also clear that the $\Z$-submodule of $U$ generated by the coefficients of $(1+u)^{a+m}$ ($a\in U$, $m\in\Z$) depends only on $a$ and not on  $m$; it is actually a $\Z$-subalgebra of $U$: indeed for all $k,l\in\N$
$${a\choose k}{a-k\choose l}={k+l\choose k}{a\choose k+l}.$$

\n More precisely for each $m\in\Z$ and $n\in\N$ the $\Z$-submodule of $U$ generated by the ${a+m\choose k}$'s for $k=0,...,n$ ($a\in U$) depends only on $a$ and $n$ and not on  $m$.

\n Finally notice that in $U[[u]]$ we have ${{\rm{d}}\over{\rm{d}}u}(1+u)^a=a(1+u)^{a-1}$. \vskip .3 truecm
\begin{example} \label{binex}

\n Let $\{x_i\}_{i\in I}$ be a $\Z$-basis of $M$. Then it is well known and trivial that:

i) The $\Z$-subalgebra $S^{(bin)}M\subseteq SV$ generated by $\{{x\choose k}
\}_{x\in M,k\in\N}$ contains $M$; 

ii) $\{{x_i\choose k}\}_{i\in I,k\in\N}$ is a set of algebra-generators (over $\Z$) of $S^{(bin)}M$;

iii) the set $\{{{\bf{x}}\choose{\bf{k}}}=\prod_{i\in I}{x_i\choose k_i}\}|{\bf{k}}:I\to\N$ finitely supported$\}$ is a $\Z$-basis of $S^{(bin)}M$.

iv) $S^{(bin)}M\cap V=M$;

v) $S^{(bin)}M$ is an integral form of $SV$ (called the algebra of binomials of $M$).

\n $S^{(bin)}M$ is also denoted $\Z^{(bin)}[x_i|i\in I]$.
\end{example}

\vskip .3 truecm
\begin{example} {\,(Review of the symmetric functions, see \cite{IM})} \label{rvsf}

\n Let $n\in\N$. It is well known that $\Z[x_1,...,x_n]^{{\cal{S}}_n}$ is an integral  form of $\Q[x_1,...,x_n]^{{\cal{S}}_n}$ and that
$\Z[x_1,...,x_n]^{{\cal{S}}_n}=\Z[e_1^{[n]},...,e_n^{[n]}]$, where the (algebraically independent for $k=1,...,n$) elementary symmetric polynomials $e_k^{[n]}$'s are defined by
\begin{equation}\label{mcd}
\prod_{i=1}^n(T-x_i)=\sum_{k\in\N}(-1)^ke_k^{[n]}T^{n-k}
\end{equation} 
and are homogeneous of degree $k$, that is $e_k^{[n]}\in\Z[x_1,...,x_n]_k^{{\cal{S}}_n}\subseteq\Q[x_1,...,x_n]_k^{{\cal{S}}_n}$.

\n It is also well known that for $n_1\geq n_2$ the natural projection $$\pi_{n_1,n_2}:\Q[x_1,...,x_{n_1}]\to\Q[x_1,...,x_{n_2}]$$ defined
by $$\pi_{n_1,n_2}(x_i)=\begin{cases}x_i&{\rm{if\, }}i\leq n_2\\ 0 &{\rm{otherwise}}\end{cases}$$ is such that
$\pi_{n_1,n_2}(e_k^{[n_1]})=
e_k^{[n_2]}
$ for all $k\in\N$.
Then $$\bigoplus_{d\geq 0}\varprojlim\Z[x_1,...,x_n]_d^{{\cal{S}}_n}=\Z[e_1,...,e_k,...]\, \, (e_k\, {\rm{inverse\, limit\, of\, the\, }}e_k^{[n]})$$ is an integral form of $\oplus_{d\geq 0}\varprojlim\Q[x_1,...,x_n]_d^{{\cal{S}}_n}$, which is called the algebra of the symmetric functions.

\n Moreover the elements 
$$p_r^{[n]}=\sum_{i=1}^nx_i^r\in\Z[x_1,...,x_n]^{{\cal S}_n}\, \, (r>0,\, n\in\N)$$
and their inverse limits $p_r\in\Z[e_1,...,e_k,...]$ ($\pi_{n_1,n_2}(p_r^{[n_1]})=p_r^{[n_2]}$ for all $r>0$ and all $n_1\geq n_2$) give another set of generators of the $\Q$-algebra of the symmetric functions: the $p_r$'s are algebraically independent and $$\bigoplus_{d\geq 0}\varprojlim\Q[x_1,...,x_n]_d^{{\cal S}_n}=\Q[p_1,...,p_r,...].$$
Finally
$\Z[e_1,...,e_k,...]$ is an integral form of $\Q[p_1,...,p_r,...]$ containing $p_r$ for all $r>0$ (more precisely a linear combination of the $p_r$'s lies in $\Z[e_1,...,e_k,...]$ if and only if it has integral coefficients), the relation between the $e_k$'s and the $p_r$'s being given by:
$$\sum_{k\in\N}(-1)^ke_ku^k={\rm{exp}}\Big(-\sum_{r>0}{p_r\over r}u^r\Big).$$

\n In this context we use the notation $$\Z[e_k|k>0]=\Z^{(sym)}[p_r|r>0
]\subseteq\Q[p_r|r>0]
;$$
to stress the dependence of the $e_k$'s on the $p_r$'s we set $e=\hat p$, that is
\begin{equation} \label{dfhp}
\hat p(u)=\sum_{k\in\N}\hat p_ku^k={\rm{exp}}\Big(\sum_{r>0}(-1)^{r-1}{p_r\over r}u^r\Big)
\end{equation}
and \begin{equation}\label{zdfhp}\Z[\hat p_k|k>0
]=\Z^{(sym)}[p_r|r>0].\end{equation}
\end{example}

\vskip .3 truecm
\begin{remark} \label{spsym}
\noindent
\n With the notations above, let $\varphi:\Q[p_1,...,p_r,...]\to U$ be an algebra-homomorphism and $a=\varphi(p_1)$:

i) if $\varphi(p_r)=0$ for $r>1$ then $\varphi(\hat p_k)=a^{(k)}$ for all $k\in\N$;

ii) if $\varphi(p_r)=a$ for all $r>0$ then $\varphi(\hat p_k)={a\choose k}$ for all $k\in\N$.

\n Hence 
$\Z^{(sym)}$ is a generalization of both $\Z^{(div)}$ and $\Z^{(bin)}$.

\end{remark}

\vskip .3 truecm

\begin{remark}\label{funtorialita}
Let $M$ be the $\Z$-module with basis $\{p_r|r>0\}$ and, as above, $V=\Q\otimes_{\Z}M$. Then:

\n i) as for the functors $S$, $S^{(div)}$ and $S^{(bin)}$, we have $\Z^{(sym)}[p_r|r>0]\cap V=M$;

\n ii) unlike the functors $S$, $S^{(div)}$ and $S^{(bin)}$, $\Z^{(sym)}[p_r|r>0]$ depends on $\{p_r|r>0\}$ and not only on $M$: for instance 
$$\Z^{(sym)}[-p_1,p_r|r>1]\neq \Z^{(sym)}[p_r|r>0]$$
(it is easy to check that these integral forms are different for example in degree 3);

\n iii) not all the sign changes 
of the $p_r$'s produce different $\Z^{(sym)}$-forms of $\Q[p_r|r>0]$: 
$$\Z^{(sym)}[(-1)^rp_r|r>0]=\Z^{(sym)}[-p_r|r>0]=\Z^{(sym)}[p_r|r>0]$$
since 
$$\exp\left(\sum_{r>0}(-1)^{r-1}{(-1)^rp_r\over r}u^r\right)=\exp\left(\sum_{r>0}(-1)^{r-1}{p_r\over r}(-u)^r\right)$$ and
$$\exp\left(\sum_{r>0}(-1)^{r-1}{-p_r\over r}u^r\right)=\exp\left(\sum_{r>0}(-1)^{r-1}{p_r\over r}u^r\right)^{-1}$$
(see remark \ref{srinv},1),ii)).
\end{remark}

\vskip .3 truecm

\n In general it is not trivial to understand whether an element of $\Q[p_r|r>0]$ belongs or not to $\Z^{(sym)}[p_r|r>0]$; proposition \ref{tmom} gives an answer to this question, which is generalized in proposition \ref{convoluzioneintera} (the examples in remark \ref{funtorialita}, ii) and iii) can be obtained also as applications of proposition \ref{convoluzioneintera}).

\vskip .3 truecm
\begin{proposition} \label {tmom}
\n Let us fix $m>0$ and let $\lambda_m:\Q[p_r|r>0]\to\Q[p_r|r>0]$ be the algebra homomorphism defined by
$\lambda_m(p_r)=p_{mr}$ for all $r>0$. 

\n Then $\Z^{(sym)}[p_r|r>0]$ $(=\Z[\hat p_k|k>0])$ is $\lambda_m$-stable.
\begin{proof}
\n For $n\in\N$ let $\lambda_m^{[n]}:\Q[x_1,...,x_n]\to\Q[x_1,...,x_n]$ be the algebra homomorphism defined by $\lambda_m^{[n]}(x_i)=x_i^m$ for all $i=1,...,n$. 

\n We obviously have that
$$\Z[x_1,...,x_n]\ \ {\rm{is\ }}\lambda_m^{[n]}{\rm{-stable}},$$
$$\Q[x_1,...,x_n]_d\ \ {\rm{is\ mapped\ to\ }}\Q[x_1,...,x_n]_{md}\,\,\forall d\geq 0,$$
$$\lambda_m^{[n]}\circ\sigma=\sigma\circ\lambda_m^{[n]}\,\,\forall n\in\N,\, \sigma\in{\cal{S}}_n,$$
$$\pi_{n_1,n_2}\circ\lambda_m^{[n_1]}=\lambda_m^{[n_2]}\circ\pi_{n_1,n_2}\,\,\forall n_1\geq n_2,$$
$$\lambda_m^{[n]}(p_r^{[n]})=p_{mr}^{[n]}\,\,\forall n\in\N,\,r>0,$$
hence there exist the limits of the $\lambda_m^{[n]}\big|_{\Q[x_1,...,x_n]_d^{{\cal S}_n}}$'s: their direct sum over $d\geq 0$  stabilizes $\oplus_{d\geq 0}\lim\Z[x_1,...,x_n]_d^{{\cal S}_n}=\Z[\hat p_k|k>0]$ and is $\lambda_m$.

\n In particular $\lambda_m(\hat p_k)\in\Z[\hat p_l|l>0]$ $\forall k\in\N$.

\end{proof}
\end{proposition}

\vskip .3 truecm
\n We also propose a second, direct, proof of proposition \ref{tmom}, which provides in addition an explicit expression of the $\lambda_m(\hat p_k)$'s in terms of the $\hat p_l$'s.

\begin{proposition}\label{tdmom}
Let $m$ and $\lambda_m$ be as in proposition \ref{tmom} and $\omega \in\C$ a primitive $m^{th}$ root of 1. Then
$$\lambda_m(\hat p(-u^m))=\prod_{j=0}^{m-1}\hat p(-\omega^j u)\in\Z[\hat p_k|k>0][[u]].$$

\begin{proof}
The equality in the statement is an immediate consequence of
$$\sum_{j=0}^{m-1}\omega^{jr}=\begin{cases}m&{\rm{if}}\ m|r\\0&{\rm{otherwise}},\end{cases}$$
so that 
$$-\sum_{j=0}^{m-1}\sum_{r>0}{p_r\over r}\omega^{jr}u^r=-\sum_{r>0}{p_{mr}\over r}u^{mr}=\lambda_m\left(-\sum_{r>0}{p_r\over r}(u^m)^r\right),$$
whose exponential is the claim.

\n Then for all $k>0$
$$\lambda_m(\hat p_k)\in\Q[\hat p_l|l>0]\cap\Z[\omega][\hat p_l|l>0]=\Z[\hat p_l|l>0]$$
since $\Q\cap\Z[\omega]=\Z$.
\end{proof}

\end{proposition}

\n In order to characterize the functions $a:\Z_+\to\Q$ such that 
$$\Z^{(sym)}[a_rp_r|r>0]\subseteq\Z^{(sym)}[p_r|r>0]$$
we introduce the notation \ref{hcappucciof}, where we rename the $p_r$'s into $h_r$ since in the affine Kac-Moody case the $\Z^{(sym)}$-construction describes the imaginary component of the integral form. Moreover from now on $p_i$ will denote a positive prime number.
\begin{notation}\label{hcappucciof}
Given $a:\Z_+\to\Q$ set $$\sum_{k\geq 0}\hat h^{\{ a \} }_ku^k=\hat h^{\{a\}}(u)=\exp\left(\sum_{r>0}(-1)^{r-1}{a_rh_r\over r}u^r\right);$$
$\uno$ denotes the function defined by $\uno_r=1$ for all $r\in\Z_+$; 

\n for all $m>0$
$\uno^{(m)}$ denotes the function defined by $\uno^{(m)}_r=\begin{cases}m&{\rm{if}}\ m|r\\0&{\rm{otherwise}}.\end{cases}$

\n Thus 
$\hat h^{\{\uno\}}(u)=\hat h(u)$ and $\hat h^{\{\uno^{(m)}\}}(-u)=\lambda_m(\hat h(-u^m))$.

\end{notation}
\vskip .3 truecm
\begin{recall}

\n The convolution product $*$ in the ring of the $\Q$-valued arithmetic functions 
$${\cal A}r=\{f:\Z_+\to\Q\}$$
is defined by
$$(f*g)(n)=\sum_{r,s:\atop rs=n}f(r)g(s).$$
The M\"obius function $\mu:\Z_+\to\Q$ defined by 
$$\mu\left(\prod_{i=1}^n p_i^{r_i}\right)=\begin{cases}(-1)^n&{\rm {if}}\ r_i=1\ \forall i\\0&{\rm{otherwise}}\end{cases}$$
{\centerline{(where the $p_i$'s are distinct positive prime integers and $r_i\geq 1$ for all $i$)}}
is the inverse of $\uno$ in the ring of the arithmetic functions.
\end{recall}

\vskip .3 truecm

\begin{proposition}\label{convoluzioneintera}
\n Let $a:\Z_+\to\Q$ be any function; then, with the notations fixed in \ref{hcappucciof},
$$\hat h^{\{a\}}_k\in\Z[\hat h_l|l>0]\ \ \forall k>0\Leftrightarrow n|(\mu*a)(n)\in\Z\ \ \forall n>0.$$

\begin{proof}
\n Remark that $a=\uno*\mu*a$, that is
$$\forall n>0\ a_n=\sum_{m|n}(\mu*a)(m)=\sum_{m|n}{(\mu*a)(m)\over m}m=\sum_{m>0}{(\mu*a)(m)\over m}\uno^{(m)}_n,$$
which means
$$a=\sum_{m>0}{(\mu*a)(m)\over m}\uno^{(m)}.$$
Let $k_m={(\mu*a)(m)\over m}$ for all $m>0$, choose $m_0>0$ such that $k_m\in\Z$ $\forall m<m_0$ and set $a^{(0)}=\sum_{m<m_0}k_m\uno^{(m)}$, $a'=a-a^{(0)}$, so that 
$$\hat h^{\{a\}}(u)=\hat h^{\{a'\}}(u)\hat h^{\{a^{(0)}\}}(u),$$ 
and, by proposition \ref{tmom} (see also notation \ref{hcappucciof}),
$$\ \ \ \hat h^{\{a^{(0)}\}}(u)\in\Z[\hat h_k|k>0][[u]].$$
It follows that

\n i) $\hat h^{\{a\}}(u)\in\Z[\hat h_k|k>0][[u]]\Leftrightarrow\hat h^{\{a'\}}(u)\in\Z[\hat h_k|k>0][[u]]$.

\n ii) $\forall n<m_0$ $\hat h^{\{a'\}}_n=0$, so that $\hat h^{\{a\}}_n=\hat h^{\{a^{(0)}\}}_n\in\Z[\hat h_k|k>0]$;

in particular $\hat h^{\{a\}}(u)\in\Z[\hat h_k|k>0][[u]]$ if $k_m\in\Z$ $\forall m>0$.

\n iii) $a'_{m_0}=(\mu*a)(m_0)=n_0k_{m_0}$ so that $\hat h^{\{a'\}}_{m_0}=k_{m_0}h_{m_0}$, which belongs to $\Z[\hat h_k|k>0]$ if and only if
$k_{m_0}\in\Z$ (see remark \ref{funtorialita},i));

in particular $\hat h^{\{a\}}(u)\not\in\Z[\hat h_k|k>0][[u]]$ if $\exists m_0\in\Z_+$ such that $k_{m_0}\not\in\Z$.

\end{proof}
\end{proposition}
\begin{proposition}\label{emmepiallaerre}
\n Let $a:\Z_+\to\Z$ be a function satisfying the condition
$$p^r|a_{mp^r}-a_{mp^{r-1}}\ \ \forall p,m\in\Z_+\ \ {\rm{with}}\ \  p\ \ {\rm{prime\ and}}\ (m,p)=1.$$
Then $n|(\mu*a)(n)$ $\forall n\in\Z_+$.

\begin{proof}

\n 
The condition $1|(\mu*a)(1)$ is equivalent to the condition $a_1\in\Z$.

\n For $n>1$ remark that $$n|(\mu*a)(n)\Leftrightarrow p^r|(\mu*a)(n)\ \ \forall p\ {\rm{prime}},\ r>0\ {\rm{such\ that}}\ p^r||n.$$
Recall that if $P
$ is the set of the prime factors of $n$ and $p\in P$ then
$$(\mu*a)(n)=\sum_{S\subseteq P}(-1)^{\#S}a_{{n\over\prod_{q\in S}q}}=$$
\begin{equation}\label{mpr}=\sum_{S'\subseteq P\setminus\{p\}}(-1)^{\#S'}(a_{{n\over\prod_{q\in S'}q}}-a_{{n\over p\prod_{q\in S'}q}}).\end{equation}
The claim follows from the remark that 
$p^r||n$ if and only if $p^r||{n\over\prod_{q\in S'}q}$. 

\end{proof}
\end{proposition}
\begin{remark} \label{vicelambda}
The viceversa of proposition \ref{emmepiallaerre} is trivially true, too, and is immediately proved applying 
(\ref{mpr}) 
to the minimal $n>0$ such that there exists $p|n$ and $r>0$ ($p^r|n$, $n=mp^r$) not satisfying the hypothesis of the statement.
\end{remark}
\n Proposition \ref{tmom} will play an important role in the study of the commutation relations in the enveloping algebra of $\hat\goths$ (see remarks \ref{stuz},vi)  and \ref{exev}) and of $\gothd$ (see remark \ref{ometiomecap} and proposition \ref{sttuz},iv)).
 
\n Proposition
 \ref{convoluzioneintera} is based on and generalizes proposition \ref{tmom}; it is a key tool in the study of the integral form in the case of $A_2^{(2)}$, see corollary \ref{hcappucciod}.

\n A more precise connection between the integral form $\Z^{(sym)}[h_r|r>0]$ of $\Q[h_r|r>0]$ and the homomorphisms $\lambda_m$'s, namely another $\Z$-basis of $\Z^{(sym)}[h_r|r>0]$ (basis defined in terms of the elements $\lambda_m(\hat h_k)$'s and arising from Garland's and Mitzman's description of the integral form of the affine Kac-Moody algebras) is discussed in appendix \ref{appendB}.

\vskip .5 truecm
\section{Some non commutative cases} \label{ncn}

\vskip .5truecm
\n We start this section with a basic remark.

\begin{remark} \label{dbsg}
\n i) Let $U_1$, $U_2$ be two $\Q$-algebras, with integral forms respectively $\tilde U_1$ and $\tilde U_2$. Then $\tilde U_1\otimes_{\Z}\tilde U_2$ is an integral form of the $\Q$-algebra $U_1\otimes_{\Q}U_2$
.

\n ii) Let $U$ be an associative unitary $\Q$-algebra (not necessarily commutative) 
and $U_1,U_2\subseteq U$ be two $\Q$-subalgebras such that $U\cong U_1\otimes_{\Q} U_2$ as $\Q$-vector spaces. If $\tilde U_1,\tilde U_2$ are integral forms of $U_1,U_2$, then $\tilde U_1\otimes_{\Z} \tilde U_2$ is an integral form of $U$ if and only if $\tilde U_2\tilde U_1\subseteq\tilde U_1\tilde U_2$.
\end{remark}

\n Remark \ref{dbsg},ii) suggests that if we have a (linear) decomposition of an algebra $U$ as an ordered tensor product of polynomial algebras $U_i$ ($i=1,...,N$), that is we have a linear isomorphism
$$U\cong U_1\otimes_{\Q}...\otimes_{\Q} U_N,$$ then one can tackle the problem of finding an integral form of $U$ by studying the commutation relations among the elements of some suitable integral forms of the $U_i$'s.

\n Glueing together in a non commutative way the different integral forms of the algebras of polynomials discussed in section \ref{intgpl} is the aim of this section, which collects the preliminary work of the paper: the main results of the following sections are applications of the formulas found here.

\vskip .3 truecm
\begin{notation}\label{lard}
\n Let $U$ be an associative $\Q$-algebra and $a\in U$. 

\n We denote by $L_a$ and $R_a$ respectively the left and right multiplication by $a$; of course $L_a-R_a=[a,\cdot]=-[\cdot,a]$.
\end{notation}

\vskip .3 truecm
\begin{lemma} \label{cle}
\n Let $U$ be an associative unitary $\Q$-algebra.

\n Consider elements $a,b,c\in U$, 
$f,g\in End(U)$ and $\alpha(u)\in U[[u]]$. 
Then:

i) if ${\exp}(f)$ and ${\exp}(g)$ converge and $[f,g]=0$ we have $${\rm{exp}}(f\pm g)={\rm{exp}}(f){\rm{exp}}(g)^{\pm 1};$$

ii) $[L_a,R_a]=0$;

iii) if $f$ is an algebra-homomorphism and $f(a)=a$ we have

$$[f,L_a]=[f,R_a]=0;$$

iv) if ${\exp}(a)$ converges so do ${\rm{exp}}(L_a)$ and ${\rm{exp}}(R_a)$, and we have
$${\rm{exp}}(L_a)=L_{{\rm{exp}}(a)},\,\,{\rm{exp}}(R_a)=R_{{\rm{exp}}(a)},\,\,{\rm{exp}}(R_a)=L_{{\rm{exp}}(a)}{\rm{exp}}([\cdot,a]);$$

v) if $\exp(a)$ and ${\rm{exp}}(c)$ converge we have $$ab=bc\Leftrightarrow {\rm{exp}}(a)b=b{\rm{exp}}(c);$$

vi) if $\exp(b)$ converges and $[b,c]=0$ we have $$[a,b]=c\Leftrightarrow a{\rm{exp}}(b)={\rm{exp}}(b)(a+c)
;$$

vii) if ${\rm{exp}}(a)$, ${\rm{exp}}(b)$ and ${\rm{exp}}(c)$ converge and $[a,c]=[b,c]=0$ then
$$[a,b]=c\Leftrightarrow {\rm{exp}}(a){\rm{exp}}(b)={\rm{exp}}(b){\rm{exp}}(a){\rm{exp}}(c)$$

viii) if ${\exp}(a)$, ${\rm{exp}}(b)$ and ${\rm{exp}}(c)$ converge and $[a,c]=[b,c]=0$ then
$$[a,b]=c\Rightarrow {\rm{exp}}(a+b)={\rm{exp}}(a){\rm{exp}}(b){\rm{exp}}(-c/2);$$

ix) if ${\rm{d}}:U\to U$ is a derivation and $[a,{\rm{d}}(a)]=0$ we have 
$${\rm{d}}({\rm{exp}}(a))={\rm{d}}(a){\rm{exp}}(a)={\rm{exp}}(a){\rm{d}}(a).$$

x) if $\alpha(u)=\sum_{r\in\N}\alpha_ru^r
$ ($\alpha_r\in U$ $\forall r\in\N$) we have
$${{\rm{d}}\over{\rm{d}}u}\alpha(u)=\alpha(u)b\Leftrightarrow \alpha(u)=\alpha_0{\rm{exp}}(bu)$$
and
$${{\rm{d}}\over{\rm{d}}u}\alpha(u)=b\alpha(u)\Leftrightarrow \alpha(u)={\rm{exp}}(bu)\alpha_0.$$

\begin{proof}
Statements v) and vi) are immediate consequence respectively of the fact that  for all $n\in\N$:

v) $a^nb=bc^n$;

vi) $ab^{(n)}=b^{(n)}a+b^{(n-1)}c$.

\n vii) follows from v) ad vi).

\n viii) follows from vii):
$$(a+b)^{(n)}=\sum_{r,s,t:\atop r+s+2t=n}{(-1)^t\over 2^t}a^{(r)}b^{(s)}c^{(t)}.$$

\n The other points are obvious.
\end{proof}
\end{lemma}

\vskip .3 truecm
\begin{proposition} \label{bdm}

\n Let us fix $m\in\Z$ and consider the $\Q$-algebra structure on $U=\Q[x]\otimes_{\Q}\Q[h]$ given by $xh=(h-m)x$. 

\n Then $\Z^{(div)}[x]\otimes_{\Z}
\Z^{(bin)}[h]$ and $\Z^{(bin)}[h]\otimes_{\Z}
\Z^{(div)}[x]$ are integral forms of $U$: their images in $U$ are closed under multiplication, and coincide. Indeed
\begin{equation}\label{fru}x^{(k)}{h\choose l}={h-mk\choose l}x^{(k)}\,\,\forall k,l\in\N\end{equation}
or equivalently, with a notation that will be useful in the following,
\begin{equation}\label{fu}
{\rm{exp}}(xu)(1+v)^h=(1+v)^h {\rm{exp}}\left({xu\over (1+v)^m}\right).
\end{equation}

\begin{proof}
The relation between $x$ and $h$ 
can be written as
$$xP(h)=P(h-m)x$$ and $$x^{(k)}P(h)=P(h-mk)x^{(k)}$$ for all $P\in\Q[h]$ and for all $k >0$. In particular it holds for $P(h)={h\choose l}$, that is
\begin{equation} \label{xvh}
x(1+v)^h=(1+v)^{h-m}x=(1+v)^h{x\over(1+v)^m}
\end{equation}
and \begin{equation}\label{xvh2}x^{(k)}(1+v)^h=(1+v)^h\left({x\over(1+v)^m}\right)^{(k)}.\end{equation}
\n The conclusion follows multiplying by $u^k$ and summing over $k$. 
\end{proof}
\end{proposition}

\vskip .3 truecm
\begin{proposition} \label{jhg}
\n Let us fix $m\in\Z$ and consider the $\Q$-algebra structure on $$U=\Q[x]\otimes_{\Q}\Q[z]\otimes_{\Q}\Q[y]$$ 
defined by $[x,z]=[y,z]=0$, $[x,y]=mz$.

\n Then $\Z^{(div)}[x]\otimes_{\Z}\Z^{(div)}[z]\otimes_{\Z}\Z^{(div)}[y]$ is an integral form of $U$.
\begin{proof}
Since $z$ commutes with $x$ and $y$ we just have to straighten $y^{(r)}x^{(s)}$. Thus the claim is a straightforward consequence of lemma \ref{cle},vii):
\begin{equation}\label{strxx}\exp(yu)\exp(xv)=\exp(xv)\exp(zuv)^{-m}\exp(yu).\end{equation}
\end{proof}
\end{proposition}

\vskip .3 truecm

\begin{proposition} \label{heise}

\n Let us fix $m,l\in\Z$ and consider the $\Q$-algebra structure on $U=\Q[h_r|r<0]\otimes_{\Q}\Q[h_0,c]\otimes_{\Q}\Q[h_r|r>0]$ given by 
$$[c,h_r]=0,\,\,[h_r,h_s]=\delta_{r+s,0}r(m+(-1)^rl)c\,\,\forall r,s\in\Z.$$

\n Then, recalling the notation $\Z[\hat h_{\pm k}|k>0]=\Z^{(sym)}[h_{\pm r}|r>0]$ and defining $U_{\Z}$ to be the $\Z$-subalgebra of $U$ generated by 
$U_{\Z}^{\pm}=\Z^{(sym)}[h_{\pm r}|r>0]$ and 
$U_{\Z}^0=\Z^{(bin)}[h_0,c]$, we have that
\begin{equation} \label{hhh}
\hat h_+(u)\hat h_-(v)=\hat h_-(v)(1-uv)^{-mc}(1+uv)^{-lc}\hat h_+(u)
\end{equation}
and $U_{\Z}=U_{\Z}^-U_{\Z}^0U_{\Z}^+$, so that
$$U_{\Z}\cong\Z^{(sym)}[h_{-r}|r>0]\otimes_{\Z}\Z^{(bin)}[h_0,c]\otimes_{\Z}\Z^{(sym)}[h_{r}|r>0]$$
is an integral form of $U$.
\begin{proof}
\ref{hhh} follows from lemma \ref{cle}, vii) remarking that 
$$\Big[\sum_{r>0}(-1)^{r-1}{h_r\over r}u^r,\sum_{s>0}(-1)^{s-1}{h_{-s}\over s}v^s\Big]=c\sum_{r>0}{m+(-1)^rl\over r}u^rv^r=$$
$$=-mc{\rm {ln}}(1-uv)-lc{\rm {ln}}(1+uv).$$

\n Of course $U_{\Z}^0U_{\Z}^-=U_{\Z}^-U_{\Z}^0$ is a $\Z$-subalgebra of $U$, $U_{\Z}^-U_{\Z}^0U_{\Z}^+\subseteq U_{\Z}$, $U_{\Z}$ is generated by
$U_{\Z}^-U_{\Z}^0U_{\Z}^+$ as $\Z$-algebra and $U_{\Z}^-U_{\Z}^0U_{\Z}^+\cong U_{\Z}^-\otimes_{\Z}U_{\Z}^0\otimes_{\Z}U_{\Z}^+$ as $\Z$-modules.

\n  Hence we need to prove that $U_{\Z}^-U_{\Z}^0U_{\Z}^+$ is a $\Z$-subalgebra of $U$, or equivalently that it is closed under left multiplication by $U_{\Z}^+$ (because it is obviously
closed under left multiplication by $U_{\Z}^-U_{\Z}^0$), which is a straightforward consequence of \ref{hhh}.

\end{proof}

\end{proposition}

\vskip .3 truecm

\begin{lemma}\label{lhlh}
Let $U$ be a $\Q$-algebra, $T:U\to U$ an automorphism, $$f\in\sum_{r>0}\Z T^ru^r\subseteq End(U[[u]]),$$
$h\in uU[[u]]$ and $x\in U$ such that
$T(h)=h$ and $[x,h]=f(x)$.
Then
$$x{\rm{exp}}(h)=\exp(h)\cdot \exp(f)(x).$$
\begin{proof}
By proposition \ref{cle},iv) 
$$x\exp(h)=\exp(h)\exp([\cdot,h])(x),$$
so we have to prove that 
$\exp([\cdot,h])(x)=\exp(f)(x),$ or equivalently that
$[\cdot,h]^n(x)=f^n(x)$ for all $n\in\N$.

\n If $n=0,1$ the claim is obvious; if $n>1$, $f^{n-1}(x)=\sum_{r>0}a_rT^ru^r(x)$ with $a_r\in\Z$ for all $r>0$, $f$ commutes with $T$, and by the inductive hypothesis
$$[\cdot,h]^n(x)=[f^{n-1}(x),h]=\left[\sum_{r>0}a_rT^{r}u^r(x),h\right]=$$
$$=\sum_{r>0}a_ru^rT^r([x,h])=\sum a_ru^rT^r f(x)=f\sum a_r u^rT^r(x)=f(f^{n-1}(x))=f^n(x).$$

\end{proof}
\end{lemma}

\begin{proposition} \label{hh}

\n Let us fix integers $m_d$'s ($d>0$) and consider elements $\{h_r,\ x_s|r>0,s\in\Z\}$ in a $\Q$-algebra $U$ 
such that
$$[h_r,x_s]=\sum_{d|r}dm_dx_{r+s}\ \ \forall r>0, s\in\Z.$$

\n Let $T$ be an algebra automorphism of $U$ such that $$T(h_r)=h_r\,\,{\rm{and}}\,\, T(x_s)=x_{s-1}\,\,\forall r>0, s\in\Z.$$

\n Then, recalling the notation $\Z[\hat h_k|k>0]=\Z^{(sym)}[h_r|r>0]$, we have that
\begin{equation}\label{cxh}
x_r\hat h_+(u)=\hat h_+(u)\cdot\left(\prod_{d>0}(1-(-T^{-1}u)^d)^{-m_d}\right)(x_r).
\end{equation}
If moreover the subalgebras of $U$ generated by $\{h_r|r>0\}$ and $\{x_r|r\in\Z\}$ are isomorphic respectively to 
$\Q[h_r|r>0]$ and $\Q[x_r|r\in\Z]$ and 
there is a $\Q$-linear isomorphism
$U\cong\Q[h_r|r>0]\otimes_{\Q}\Q[x_r|r\in\Z]$
then
$$\Z^{(sym)}[h_r|r>0]\otimes_{\Z}\Z^{(div)}[x_r|r\in\Z]$$ is an integral form of $U$.
\begin{proof}
This is an application of lemma \ref{lhlh}: let $h=\sum_{r>0}(-1)^{r-1}{h_r\over r}u^r$; then
$$[x_0,h]=\sum_{r>0}{(-1)^r\over r}u^r\sum_{d|r}dm_dT^{-r}(x_0)=$$
$$=\sum_{d>0}\sum_{s>0}{(-1)^{ds}\over s}m_dT^{-ds}u^{ds}(x_0)=f(x_0)$$
where $$f=
-\sum_{d>0}m_d\ln(1-(-1)^{d}T^{-d}u^d).$$
Then
$$x_0\hat h_+(u)=\hat h_+(u)\cdot\exp(f)(x_0)=\hat h(u)\cdot\left(\prod_{d>0}(1-(-T^{-1}u)^d)^{-m_d}\right)(x_0),$$
and the analogous statement for $x_r$ follows applying $T^{-r}$.

\n Remark that $\prod_{d>0}(1-(-T^{-1}u)^d)^{-m_d}=\sum_{r\geq 0}a_rT^{-r}u^r$ with $a_r\in\Z$ $\forall r\in\N$; the hypothesis on the commutativity of the subalgebra generated by the $x_r$'s implies that 
$(\sum_{r\geq 0}a_rx_ru^r)^{(k)}$ lies in the subalgebra of $U$ generated by the divided powers $\{x_r^{(k)}|r\in\Z,k\geq 0\}$, which allows to conclude the proof thanks to the last hypotheses on the structure of $U$.

\end{proof}
\end{proposition}

\vskip .3 truecm

\begin{remark} \label{praff}

\n Proposition \ref{hh}, implies proposition \ref{bdm}: indeed when $m_1=m$, $m_d=0$ $\forall d>1$ we have a projection
$h_r\mapsto h, x_r\mapsto x$, which  maps 
$\exp(x_0u)$ to $\exp(xu)$, $\hat h(u)$ to $(1+u)^h$ and $T$ to the identity.
\end{remark}

\vskip .5 truecm
\section{The integral form of $\goths$ ($A_1$)} \label{sld}
\vskip .5truecm

\n The results about $\goths$ and the $\Z$-basis of the integral form $\u_{\Z}(\goths)$ of its enveloping algebra $\u(\goths)$ are well known (see \cite{Ko} and \cite{S}). Here we recall the description of $\uz(\goths)$ in terms of the non-commutative generalizations described in section \ref{ncn}, with the notations of the commutative examples given in section \ref{intgpl}.

\n The proof expressed in this language has the advantage to be easily generalized to the affine case.
\vskip .3 truecm
\begin{definition} \label{sl2}

\n $\goths$ (respectively $\u(\goths)$) is the Lie algebra (respectively the associative algebra) over $\Q$ generated by $\{e,f,h\}$ with relations $$[h,e]=2e,\, [h,f]=-2f,\, [e,f]=h.$$

\n $\uz(\goths)$ is the $\Z$-subalgebra of $\u(\goths)$ generated by $\{e^{(k)},f^{(k)}| \; k\in\N\}$. 
\end{definition}

\vskip .3 truecm
\begin{theorem}\label{trdc}
\n Let $\u^+$, $\u^-$, $\u^0$ denote the $\Q$-subalgebras of $\u(\goths)
$ generated respectively by $e$, by $f$, by $h$.

\n Then $\u^+\cong\Q[e]$, $\u^-\cong\Q[f]$, $\u^0\cong\Q[h]$ and $\u(\goths)\cong\u^-\otimes\u^0\otimes\u^+$; moreover 
\begin{equation} \label{usldi}
\uz(\goths)\cong\Z^{(div)}[f]\otimes_{\Z}\Z^{(bin)}[h]\otimes_{\Z}\Z^{(div)}[e]
\end{equation}
is an integral form of $\u(\goths)$.
\begin{proof}
Thanks to proposition \ref{bdm}, we just have to study the commutation between 
$e^{(k)}$ and $f^{(l)}$ for $k,l\in\N$.

\n Let us recall the commutation relation
\begin{equation} \label{efu}
e\exp(fu)=\exp(fu)(e+hu-fu^2)
\end{equation}
which is a direct application of lemma \ref{cle},iv) and of the relations 
$[e,f]=h$, $[h,f]=-2f$ and $[f,f]=0$.

\n We want to prove that in $\u(\goths)[[u,v]]$
\begin{equation} \label{cef}
{\rm{exp}}(eu){\rm{exp}}(fv)={\rm{exp}}\Big({fv\over 1+uv}\Big)(1+uv)^h{\rm{exp}}\Big({eu\over 1+uv}\Big).
\end{equation}
Let $F(u)={\rm{exp}}\Big({fv\over 1+uv}\Big)(1+uv)^h{\rm{exp}}\Big({eu\over 1+uv}\Big).$

\n It is obvious  that $F(0)={\rm{exp}}(fv)$; hence our claim is equivalent to $${{\rm{d}}\over{\rm{d}}u}F(u)=eF(u).$$
To obtain this result we derive remarking lemma \ref{cle},ix) and then apply formulas \ref{xvh} and \ref{efu}:
$${{\rm{d}}\over{\rm{d}}u}F(u)=$$
$$={\rm{exp}}\Big({fv\over 1+uv}\Big)(1+uv)^h{e\over(1+uv)^2}{\rm{exp}}\Big({eu\over 1+uv}\Big)+$$
$$+{\rm{exp}}\Big({fv\over 1+uv}\Big)\Big({hv\over 1+uv}-{fv^2\over(1+uv)^2}\Big)(1+uv)^h{\rm{exp}}\Big({eu\over 1+uv}\Big)=$$
$$={\rm{exp}}\Big({fv\over 1+uv}\Big)\Big(e+{hv\over 1+uv}-{fv^2\over(1+uv)^2}\Big)(1+uv)^h{\rm{exp}}\Big({eu\over 1+uv}\Big)=$$
$$=eF(u).$$
Remarking that $${xu\over 1+v}\in\Z[x][[u,v]],\,\,{\rm{hence}}\,\,\Big({xu\over 1+v}\Big)^{(k)}\in\Z^{(div)}[x][[u,v]]\,\,\forall k\in\N,$$
it follows that the right hand side of \ref{usldi} is an integer form of $\u(\goths)$ (containing $\uz(\goths)$).

\n Finally remark that inverting the exponentials on the right hand side, the formula (\ref{cef}) gives an expression of $(1+uv)^h$ in terms of the divided powers of $e$ and $f$, so that 
$\Z^{(bin)}[h] \subseteq \uz(\goths)$, which completes the proof.

\end{proof}
\end{theorem}

\vskip .5 truecm
\section{The integral form of $\hat{\goths}$ ($A_1^{(1)}$)} \label{slh}
\vskip .5truecm

\n The results about $\hat{\goths}$ and the integral form $\hu_{\Z}$ of its enveloping algebra $\hu$ are due to Garland (see \cite{HG}). Here we simplify the description of the imaginary positive component of $\hu_{\Z}$ 
proving that it is an algebra of polynomials over $\Z$ and give a compact and complete proof of the assertion that the set given in theorem \ref{trm} is actually a $\Z$-basis of $\huz$. This proof has the advantage, following \cite{DM}, to reduce the long and complicated commutation formulas to compact, simply readable and easily proved ones. It is evident from this approach that the results for $\hat{\goths}$ are generalizations of those for $\goths$, so that the commutation formulas 
arise naturally recalling the homomorphism 
\begin{equation} \label{evaluation}
ev:\hat{\goths}=\goths\otimes\Q[t^{\pm 1}]\oplus\Q c\to\goths\otimes\Q[t^{\pm 1}]\to\goths
\end{equation}
induced by the evaluation of $t$ at 1
.

\n On the other hand these results and the strategy for their proof will be shown to be in turn generalizable to $\gothd$.

\n As announced in the introduction, the proof of theorem \ref{trm} is based on a few results:
proposition \ref{zzk}, proposition \ref{pum}, lemma \ref{limt}, and proposition \ref{exefh}.

\vskip .3 truecm

\begin{definition} \label{hs2}
\n $\hat{\goths}$ (respectively $\hu
$) is the Lie algebra (respectively the associative algebra) over $\Q$ generated by $\{x_r^+,x_r^-,h_r,c|r\in\Z\}$ with relations $$c\,\,\,{\rm{is\,\,central}},$$
$$[h_r,h_s]=2r\delta_{r+s,0}c,\,\,\,[h_r,x_s^{\pm}]=\pm 2x_{r+s}^{\pm}$$
$$[x_r^+,x_s^+]=0=[x_r^-,x_s^-],$$
$$[x_r^+,x_s^-]=h_{r+s}
+r\delta_{r+s,0}c.$$
\n Notice that $\{x_r^+,x_r^-|r\in\Z\}$ generates $\hu$.

\n $\hu^+$, $\hu^-$, $\hu^0$ are the subalgebras of $\hu$ generated respectively by $\{x_r^+|r\in\Z\}$, $\{x_r^-|r\in\Z\}$, $\{c,h_r|r\in\Z\}$.

\n $\hu^{0,+}$, $\hu^{0,-}$, $\hu^{0,0}$, are the subalgebras of $\hu$ (of $\hu^0$) generated respectively by $\{h_r|r>0\}$, $\{h_r|r<0\}$, $\{c,h_0\}$.
\end{definition}
\vskip .3 truecm

\begin{remark} \label{hefp}
\n $\hu^+$, $\hu^-$ are (commutative) algebras of polynomials:
$$\hu^+\cong\Q[x_r^+|r\in\Z],\,\,\,\hu^-\cong\Q[x_r^-|r\in\Z];$$
$\hu^0$ is not commutative: $[h_r,h_{-r}]=2rc$;

\n $\hu^{0,+}$, $\hu^{0,-}$, $\hu^{0,0}$, are (commutative) algebras of polynomials:
$$\hu^{0,+}\cong\Q[h_r|r>0],\,\,\,\hu^{0,-}\cong\Q[h_r|r<0],\,\,\,\hu^{0,0}\cong\Q[c,h_0];$$
Moreover we have the following ``triangular'' decompositions:
$$\hu\cong\hu^-\otimes\hu^0\otimes\hu^+,$$
$$\hu^0\cong\hu^{0,-}\otimes\hu^{0,0}\otimes\hu^{0,+}.$$
Remark that the images in $\hu$ of $\hu^-\otimes\hu^0$ and $\hu^0\otimes\hu^+$ are  subalgebras of $\hu$ and the images of 
$\hu^{0,-}\otimes\hu^{0,0}$ and $\hu^{0,0}\otimes\hu^{0,+}$ are commutative subalgebras of $\hu^0$.

\end{remark}

\vskip .3 truecm

\begin{definition} \label{hto}
\n $\hu$ is endowed with the following anti/auto/homo/morphisms:

\n $\sigma$ is the antiautomorphism defined on the generators by:
$$x_r^+\mapsto x_r^+,\,\,\,x_r^-\mapsto x_r^-,\,\,\,(\Rightarrow h_r\mapsto-h_r,\,\,\,c\mapsto -c);$$
$\Omega$ is the antiautomorphism defined on the generators by:
$$x_r^+\mapsto x_{-r}^-,\,\,\,x_r^-\mapsto x_{-r}^+,\,\,\,(\Rightarrow h_r\mapsto h_{-r},\,\,\,c\mapsto c);$$
\n $T$ is the automorphism defined on the generators by:
$$x_r^+\mapsto x_{r-1}^+,\,\,\,x_r^-\mapsto x_{r+1}^-,\,\,\,(\Rightarrow h_r\mapsto h_r-\delta_{r,0}c,\,\,\,c\mapsto c);$$
\n for all $m\in\Z$, $\lambda_m$ is the homomorphism defined on the generators by:
$$x_r^+\mapsto x_{mr}^+,\,\,\,x_r^-\mapsto x_{mr}^-,\,\,\,(\Rightarrow h_r\mapsto h_{mr},\,\,\,c\mapsto mc).$$
\end{definition}
\vskip .3 truecm

\begin{remark} \label{hti}
\n $\sigma^2={\rm{id}}_{\hu}$, $\Omega^2={\rm{id}}_{\hu}$, $T$ is invertible of infinite order;

\n $\lambda_{-1}^2=\lambda_1={\rm{id}}_{\hu}$; $\lambda_m$ is not invertible if $m\neq\pm 1$; $\lambda_0=ev$ (through the identification
$<x_0^+,x_0^-,h_0>\cong<e,f,h>$).
\end{remark}

\begin{remark} \label{htc}
\vskip .3 truecm
\n $\sigma\Omega=\Omega\sigma$, $\sigma T=T\sigma$, $\sigma\lambda_m=\lambda_m\sigma$ for all $m\in\Z$;

\n $\Omega T=T\Omega$, $\Omega\lambda_m=\lambda_m\Omega$ for all $m\in\Z$;

\n $\lambda_m T^{\pm 1}=T^{\pm m}\lambda_m$ for all $m\in\Z$;

\n $\lambda_m\lambda_n=\lambda_{mn}$, for all $m,n\in\Z$. 
\end{remark}

\vskip .3 truecm

\begin{remark} \label{hbs}

\n $\sigma\big|_{\hu^{\pm}}={\rm{id}}_{\hu^{\pm}},\,\,\,\sigma(\hu^{0,\pm})=\hu^{0,\pm},\,\,\,\sigma(\hu^{0,0})=\hu^{0,0}$.

\n $\Omega(\hu^{\pm})=\hu^{\mp},\,\,\,\Omega(\hu^{0,\pm})=\hu^{0,\mp},\,\,\,\Omega\big|_{\hu^{0,0}}={\rm{id}}_{\hu^{0,0}}$.

\n $T(\hu^{\pm})=\hu^{\pm},\,\,\,T\big|_{\hu^{0,\pm}}={\rm{id}}_{\hu^{0,\pm}},\,\,\,
T(\hu^{0,0})=\hu^{0,0}$.

\n For all $m\in\Z$ $\lambda_m(\hu^{\pm})\subseteq\hu^{\pm},\,\,\,\lambda_m(\hu^0)=\hu^0,\,\,\,\lambda_m(\hu^{0,0})\subseteq\hu^{0,0}$, $$\lambda_m(\hu^{0,\pm})\subseteq
\begin{cases}\hu^{0,\pm}&{\rm{if}}\,m>0\cr \hu^{0,\mp}&{\rm{if}}\, m<0\cr \hu^{0,0}&{\rm{if}}\,m=0.
\end{cases}$$

\end{remark}

\vskip .3 truecm

\begin{definition}\label{hhuz}
\n Here we define some $\Z$-subalgebras of $\hu$:

\n $\huz$ is the $\Z$-subalgebra of $\hu
$ generated by $\{(x_r^{+})^{(k)},(x_r^{-})^{(k)}|r\in\Z,k\in\N\}$;

\n $\huz^{\pm}=\Z^{(div)}[x_r^{\pm}|r\in\Z]$; 

\n $\huz^{0,0}=\Z^{(bin)}[h_0,c]$;

\n $\huz^{0,\pm}=\Z^{(sym)}[h_{\pm r}|r>0]$;

\n $\huz^0$ 
is the $\Z$-subalgebra of $\hu$ generated by $\huz^{0,-}$, $\huz^{0,0}$ and $\huz^{0,+}$.

The notations are those of section \ref{intgpl}. 
\end{definition}

\vskip .3 truecm
\n We want to prove that $\huz^0 =\huz^{0,-}\huz^{0,0}\huz^{0,+}$, so that it is an integral form of $\hu^0$, and that $\huz=\huz^-\huz^0\huz^+$, so that $\huz$ is an integral form of $\hu$.

\n As in the case of $\goths$,
working in $\hu[[u]]$ (see the notation below) simplifies enormously the proofs and gives a deeper insight to the question.

\vskip .3 truecm

\begin{notation}\label{hgens}
\n We shall consider the following elements in $\hu[[u]]$:
$$x^+(u)=\sum_{r\geq 0}x_r^+u^r=\sum_{r\ge0} T^{-r}u^r(x_0^+),$$
$$x^-(u)=\sum_{r\geq 0}x_{r+1}^-u^r=\sum_{r\ge0} T^{r}u^r(x_1^-),$$
$$h_{\pm}(u)=\sum_{r\geq 1}(-1)^{r-1}{h_{\pm r}\over r}u^r,$$
$$\hat h_{\pm}(u)={\exp}(h_{\pm}(u))=\sum_{r\geq 0}\hat h_{\pm r} u^r.$$
\end{notation}

\vskip .3 truecm{}
\begin{remark} \label{nev}
\n Notice that $ev \circ T=ev$ and 
$$ev(x^+(-u))=ev\left({{1}\over{1+T^{-1}u}}x_0^+\right)={e\over 1+u},$$
$$ev(x^-(-u))=ev\left({{T}\over{1+Tu}}x_0^-\right)={f\over 1+u},$$
$$ev(h_{\pm}(u))=h{\rm{ln}}(1+u),$$
$$ev(\hat h_{\pm}(u))=(1+u)^h.$$
\end{remark}
\vskip .3 truecm

\begin{remark} \label{stuz}
Here we list some obvious remarks.

\n i) $\huz^{\pm}\subseteq\huz\cap\hu^{\pm}$ and $\huz$ is the $\Z$-subalgebra of $\hu$ generated by $\huz^+\cup\huz^-$; 

\n ii) $\huz^{\pm}$, $\huz^{0,0}$, $\huz^{0,\pm}$ and $\huz^{0,\pm}\huz^{0,0}=\huz^{0,0}\huz^{0,\pm}$ are integral forms respectively of 
$\hu^{\pm}$, $\hu^{0,0}$, $\hu^{0,\pm}$ and $\hu^{0,\pm}\hu^{0,0}=\hu^{0,0}\hu^{0,\pm}$;

\n iii) $\huz$ and $\huz^{0,0}$ are  
stable under $\sigma$, $\Omega$, $T^{\pm 1}$, $\lambda_m $ for all $m\in\Z$;

\n iv) $\huz^{\pm}$ is stable under $\sigma$, $T^{\pm 1}$, $\lambda_m $ for all $m\in\Z$ and $\Omega(\huz^{\pm})=\huz^{\mp}$;

\n v) $\huz^{0,\pm}$ is stable under $\sigma$, $T^{\pm 1}$ and $\Omega(\huz^{0,\pm})=\lambda_{-1}(\huz^{0,\pm})=\huz^{0,\mp}$: more precisely
$$\sigma(\hat h_{\pm}(u))\!=\!\hat h_{\pm}(u)^{-1},\,\Omega(\hat h_{\pm}(u))\!=\!\lambda_{-1}(\hat h_{\pm}(u))\!=\!\hat h_{\mp}(u),\,T^{\pm 1}(\hat h_{\pm}(u))\!=\!\hat h_{\pm}(u);$$
vi)
for $m\in\Z$
$$\lambda_m(\huz^{0,\pm})\subseteq
\begin{cases}\huz^{0,\pm}&{\rm{if}}\,m>0\cr \huz^{0,\mp}&{\rm{if}}\, m<0\cr \huz^{0,0}&{\rm{if}}\,m=0,
\end{cases}$$
thanks to v), to proposition \ref{tmom} and to remarks \ref{htc} and \ref{nev}.
\end{remark}

\vskip.3 truecm
\begin{remark} \label{tmfv}
\n 
The elements $\hat h_k$'s with $k>0$ generate the same $\Z$-subalgebra of $\hu$ as the elements $\Lambda_k$'s ($k\geq 0$) defined in \cite{HG}. 

\n Indeed let 
$$\sum_{n\geq 0}p_nu^n=P(u)=\hat h(-u)^{-1};$$ then 
remarks \ref{srinv},1,ii) and \ref{funtorialita},iii) imply that $\Z[\hat h_k|k>0]=\Z[p_{n}|n>0]$; but 
$${{\rm{d}}\over{\rm{d}}u}P(u)=P(u)\sum_{r>0}h_ru^{r-1},$$
that is
$$p_0=1,\ \ p_n={1\over n}\sum_{r=1}^nh_rp_{n-r}\ \forall n>0,$$
hence $p_n=\Lambda_{n-1}$ $\forall n\geq 0$.

\n On the other hand applying $\lambda_m$ we get 
$$\lambda_m(p_0)=1,\ \ \lambda_m(p_n)={1\over n}\sum_{r=1}^nh_{rm}\lambda_m(p_{n-r}),$$  so that
$\lambda_m(p_n)=\lambda_m(\Lambda_{n-1})=\Lambda_{n-1}{(\xi(m))}$ (see \cite{HG}).

\end{remark}

\vskip .3 truecm

\begin{remark} \label{hrs}
\n Remark that for all $r\in\Z$ the subalgebra of $\hat{\goths}$ generated by $$\{x_r^+,x_{-r}^-,h_0+rc\}$$ maps isomorphically onto $\goths$ through the evaluation homomorphism $ev$ (see formula \ref{evaluation}). On the other hand for each $r\in\Z$ there is an injection $\u(\goths)\to\hu$:
$$e\mapsto x_r^+,\,\,f\mapsto x_{-r}^-,\,\,h\mapsto h_0+rc.$$
In particular theorem \ref{trdc}, implies that the elements ${h_0+rc\choose k}$ belong to $\huz$ for all $r\in\Z, k\in\N$ (thus, remarking that the elements ${c\choose k}$'s are central and the example \ref{binex}, we get that $\huz^{0,0}\subseteq\huz$)
and proposition \ref{bdm} implies that 
$\huz^{0,0}\huz^+$ and $\huz^-\huz^{0,0}$ are integral forms respectively of $\hu^{0,0}\hu^+$ and $\hu^-\hu^{0,0}$.
\end{remark}
\vskip .3 truecm

\vskip .3 truecm
\begin{proposition} \label{zzk}
\n The following identity holds in $\hu$:

$$\hat h_+(u)\hat h_-(v)=\hat h_-(v)(1-uv)^{-2c}\hat h_+(u).$$

\n $\huz^0=\huz^{0,-}\huz^{0,0}\huz^{0,+}$:
it 
is an integral form of $\hu^0$.

\begin{proof}

\n Since $[h_r,h_s]=2r\delta_{r+s,0}c$,
the claim is proposition \ref{heise} with $m\!=\!2$, $l\!=\!0$.
\end{proof}
\end{proposition}

\vskip .3 truecm
\begin{proposition}\label{pum}
\n The following identity holds in $\hu$:
\begin{equation} \label{xup}
x_0^+\hat h_+(u)=
\hat h_+(u)(1+T^{-1}u)^{-2}(x_0^+). 
\end{equation}

\n Hence for all $k\in\N$
\begin{equation} \label{xup2}
(x_0^+)^{(k)}\hat h_+(u)= \hat h_+(u)((1+T^{-1}u)^{-2}(x_0^+))^{(k)}.
\end{equation}
\begin{proof}
The claim follows from proposition \ref{hh} with $m_1=2$, $m_d=0 \; \forall d>1$, and from \ref{divpoweq}.

\end{proof}
\end{proposition}
\vskip .3 truecm

\begin{remark}
The identity (\ref{xup}) can be written as $$x_0^+\hat h_+(u)=\hat h_+(u){{\rm{d}}\over{\rm{d}}u}(ux^+(-u)).$$ Indeed 
$$(1+T^{-1}u)^{-2}(x_0^+)=\sum_{r\in\N}(-1)^r(r+1)x_r^+u^r={{\rm{d}}\over{\rm{d}}u}(ux^+(-u)).$$
\end{remark}

\begin{remark} Remark that the identity (\ref{xup2}) is the affine version of 
\begin{equation} \label{xup3}
e^{(k)}(1+u)^h=(1+u)^{h}\left({e\over (1+u)^2}\right)^{(k)}
\end{equation} (see equation (\ref{xvh2})); indeed $ev$ maps (\ref{xup2}) to (\ref{xup3}).
\end{remark}

\begin{corollary}\label{cum}
\n $\huz^+\huz^{0,\pm}\subseteq\huz^{0,\pm}\huz^+$ and $\huz^{\pm}\huz^0=\huz^0\huz^{\pm}$.

\n Then $\huz^0\huz^+$ and $\huz^-\huz^0$ are integral forms respectively of $\hu^0\hu^+$ and $\hu^-\hu^0$.
\begin{proof}
Applying $T^{-r}$ to (\ref{xup2}),  we find that $(x_r^+)^{(k)}\hat h_+(u)\subseteq\hat h_+(u)\huz^+[[u]]$ $\forall r\in\Z,k\in\N$, hence $\huz^+\hat h_+(u)\subseteq\hat h_+(u)\huz^+[[u]]$ and $\huz^+\huz^{0,+}\subseteq\huz^{0,+}\huz^+$. From this, applying $\lambda_{-1}$ we get $\huz^+\huz^{0,-}\subseteq\huz^{0,-}\huz^+$, hence $\huz^+\huz^0\subseteq\huz^0\huz^+$ thanks to remark \ref{hrs}. Finally applying $\Omega$ we obtain that $\huz^0\huz^-\subseteq\huz^-\huz^0$ and applying $\sigma$ we get the reverse inclusions.
\end{proof}
\end{corollary}

\vskip .3 truecm

We are now left to prove that $\huz^+\huz^-\subseteq\huz^-\huz^0\huz^+$ and that $\huz^0\subseteq\huz$.

\n To this aim we study the commutation relations between $(x_r^+)^{(k)}$ and $(x_s^-)^{(l)}$ or equivalently between ${\rm{exp}}(x_r^+u)$ and ${\rm{exp}}(x_s^-v)$. 
\vskip .3 truecm
\begin{remark}\label{exev}
\n Remark \ref{hrs}, implies that ${\rm{exp}}(x_r^+u){\rm{exp}}(x_{-r}^-v)\in\huz^-\huz^0\huz^+$ for all $r\in\Z$.

\n In order to prove a similar result for
${\rm{exp}}(x_r^+u){\rm{exp}}(x_s^-v)$ when $r+s\neq 0$ remark that in general $${\rm{exp}}(x_r^+u){\rm{exp}}(x_s^-v)=T^{-r}\lambda_{r+s}({\rm{exp}}(x_0^+u){\rm{exp}}(x_1^-v)),$$ so that remark \ref{stuz},iv),v),vi) allows us to 
reduce to the case $r=0$, $s=1$.

\n This case will turn out to be enough also to prove that $\huz^0\subseteq\huz$.
\end{remark}

\vskip .3 truecm
\begin{remark} \label{hev}
\n In the study of the commutation relations in $\huz$ remark that 
$$ev(\exp(x_0^+u)\exp(x_1^-v))=\exp(eu)\exp(fv)$$ and that straightening $\exp(x_0^+u)\exp(x_1^-v)$ through the triangular decomposition $\hu\cong\hu^-\otimes\hu^0\otimes\hu^+$ we get an element of 
$\hu[[u,v]]$ 
whose coefficients involve $x_{r+1}^-,h_{r+1},\, x_r^+$ with $r\geq 0$ and whose image through $ev$ is 
$$\exp\Big({fv\over 1+uv}\Big)(1+uv)^h\exp\Big({eu\over 1+uv}\Big)$$
(see remark \ref{nev}).

\n Viceversa once we have such an expression for 
$\exp(x_0^+u)\exp(x_1^-v)$ applying $T^{-r}\lambda_{r+s}$ we can deduce from it the identity (\ref{cef}) and the expression for 
$\exp(x_r^+u)\exp(x_s^-v)$ for all $r,s\in\Z$ (also in the case $r+s=0$).

\n Remark that 
$${\rm{exp}}(vx^-(-uv))\hat h_+(uv){\rm{exp}}(ux^+(-uv))$$
is an element of $\hu[[u,v]]$ which has the required properties (see remark \ref{nev}) and belongs to $\huz^-\huz^0\huz^+$. 

\vskip .3 truecm
\end{remark}

Our aim is to prove that 
$$\exp(x_0^+u)\exp(x_1^-v)=\exp(vx^-(-uv))\hat h_+(uv)\exp(ux^+(-uv)).$$

\vskip .3 truecm
\begin{lemma} \label{limt}
\n In $\hu[[u,v]]$ we have
$$x_0^+\exp(vx^-(-uv))=\exp(vx^-(-uv))\Big(x_0^++{{\rm{d}}h_+(uv)\over{\rm{d}}u}+{{\rm{d}}vx^-(-uv)\over{\rm{d}}u}\Big).$$
\begin{proof}
The claim follows from lemma \ref{cle},iv) remarking that 
$$[x_0^+,vx^-(-uv)]=v\sum_{r\in\N}h_{r+1}(-uv)^r={{\rm{d}}\over{\rm{d}}u}\sum_{r\in\N}{h_{r+1}\over r+1}(-1)^{r}(uv)^{r+1}={{\rm{d}}h_+(uv)\over{\rm{d}}u},$$
$$\Big[{{\rm{d}}h_+(uv)\over{\rm{d}}u},vx^-(-uv)\Big]=-2v^2\sum_{r,s\in\N}x_{r+s+2}^+(-uv)^{r+s}=$$
$$=-2v^2\sum_{r\in\N}(r+1)x_{r+2}^-(-uv)^r=2{{\rm{d}}vx^-(-uv)\over{\rm{d}}u}$$
and 
$$\Big[{{\rm{d}}vx^-(-uv)\over{\rm{d}}u},vx^-(-uv)\Big]=0.$$
\end{proof}
\end{lemma}

\vskip .3 truecm
\begin{proposition}\label{exefh}
\n In $\hu[[u,v]]$ we have
$${\rm{exp}}(x_0^+u){\rm{exp}}(x_1^-v)={\rm{exp}}(vx^-(-uv))\hat h_+(uv){\rm{exp}}(ux^+(-uv)).$$
\begin{proof}
Let $F(u)={\rm{exp}}(vx^-(-uv))\hat h_+(uv){\rm{exp}}(ux^+(-uv))$.
It is clear that $F(0)={\rm{exp}}(x_1^-v)$, so that it is enough to prove that 
$${{\rm{d}}\over{\rm{d}}u}F(u)=x_0^+F(u).$$
Remark that, thanks to the derivation rules (lemma \ref{cle},ix)), to proposition \ref{pum}, and to lemma \ref{limt}, we have:
$${{\rm{d}}\over{\rm{d}}u}F(u)={\rm{exp}}(vx^-(-uv))\hat h_+(uv){{{\rm d}}\over{\rm{d}}u}(ux^+(-uv)){\rm{exp}}(ux^+(-uv))+$$
$$+{\rm{exp}}(vx^-(-uv))\Big({{{\rm d}}\over{\rm{d}}u}h_+(uv)+{{{\rm d}}\over{\rm{d}}u}(vx^-(-uv))\Big)\hat h_+(uv){\rm{exp}}(ux^+(-uv))=$$
$$={\rm{exp}}(vx^-(-uv))\Big(x_0^++{{{\rm d}}(h_+(uv)+vx^-(-uv))\over{\rm{d}}u}\Big)\hat h_+(uv){\rm{exp}}(ux^+(-uv))=$$
$$=x_0^+{\rm{exp}}(vx^-(-uv))\hat h_+(uv){\rm{exp}}(ux^+(-uv))=x_0^+F(u).$$
\vskip .3 truecm
\end{proof}
\end{proposition}

\begin{corollary} \label{tfin}
\n $\huz^0\subseteq\huz$. 
\begin{proof}
That $\huz^{0,+}\subseteq\huz$ is a consequence of proposition \ref{exefh} inverting the exponentials (see the proof theorem \ref{trdc}), which implies also (applying $\Omega$) that 
$\huz^{0,-}\subseteq\huz$; the claim then follows thanks to remark \ref{hrs}.
\end{proof}
\end{corollary}

\begin{proposition} \label{strutmodulo}
$\huz^-\huz^0\huz^+$ is a $\Z$-subalgebra of $\hu$ (hence $\huz=\huz^-\huz^0\huz^+$).

\begin{proof}
We want to prove that $\huz^-\huz^0\huz^+$ (which is obviously a $\huz^-$-module and, by corollary \ref{cum}, a $\huz^0$-module) is also a 
$\huz^+$-module, or equivalently that 
$
\huz^+\huz^-\subseteq\huz^-\huz^0\huz^+$. 

\n By proposition \ref{exefh} together with remark \ref{exev}, formula \ref{cef} and remark \ref{hrs} 
we have that $
y_+y_-\in\huz^-\huz^0\huz^+$ in the particular case when $y_+=(x_r^+)^{(k)}$ and $y_-=(x_s^-)^{(l)}$, thus we just need to perform the correct induction to deal with the general $y_{\pm}\in\huz^{\pm}$.

\n Remark that setting $$deg(x_r^{\pm})=\pm 1,\ \ deg(h_r)=deg(c)=0$$
induces a $\Z$-gradation on $\hu$ (since the relations defining $\hu$ are homogeneous) and on $\huz$ (since its generators are homogeneous), which is preserved by $\sigma$, $T^{\pm 1}$ and $\lambda_m$ $\forall m\in\Z$; in particular it induces $\N$-gradations
$$\hu^{\pm}=\bigoplus_{k\in\N}\hu_{\pm k}^{\pm}, \qquad \qquad \hu_{\Z}^{\pm}=\bigoplus_{k\in\N}\hu_{\Z,\pm k}^{\pm}$$
with the properties that 
$$\Omega(\hu_{\Z,\pm k}^{\pm})= \hu_{\Z,\mp k}^{\mp},$$
$$\hu_{\Z,k}^{+}=\sum_{n\in\N\atop k_1+...+k_n=k}\Z(x_{r_1}^{+})^{(k_1)}\cdot ... \cdot(x_{r_n}^{+})^{(k_n)}=\sum_{r\in\Z}\Z(x_r^{+})^{(k)}+\sum_{k_1,k_2>0\atop k_1+k_2=k}\hu_{\Z, k_1}^{+}\hu_{\Z, k_2}^{+},$$
$$\hu_{\Z,k}^{+}\huz^0=\huz^0\hu_{\Z,k}^{+}\ \ {\rm{(because}}\ \hu_{k}\hu^0=\hu^0\hu_{k}\ {\rm{and}}\ \huz^{+}\huz^0=\huz^0\huz^{+}{\rm{)}}$$
and
$$[\hu_k^+,\hu_{-l}^-]\subseteq\sum_{m>0}\hu_{-l+m}^-\hu^0\hu_{k-m}^+\ \ \forall k,l\in\N.$$
We want to prove that 
\begin{equation}\label{uzruzs}\hu_{\Z,k}^+\hu_{\Z,-l}^-\subseteq\sum_{m\geq 0}\hu_{\Z,-l+m}^-\huz^0\hu_{\Z,k-m}^+\ \ \forall k,l\in\N,\end{equation}
the claim being obvious for $k=0$ or $l=0$.

\n Suppose $k\neq 0$, $l\neq 0$ and the claim true for all $(\tilde k,\tilde l)\neq (k,l)$ with $\tilde k\leq k$ and $\tilde l\leq l$ .
Then:

a) proposition \ref{exefh} together with remark \ref{exev}, formula (\ref{cef}) and remark \ref{hrs} imply that $$(x_r^+)^{(k)}(x_s^-)^{(l)}\in\sum_{m\geq 0}\hu_{\Z,-l+m}^-\huz^0\hu_{\Z,k-m}^+\ \ \forall r,s\in\Z;$$

b) if $k_1,k_2>0$ are such that $k_1+k_2=k$ or $l_1,l_2>0$ are such that $l_1+l_2=l$, then
$$\hu_{\Z,k_1}^+\hu_{\Z,k_2}^+\hu_{\Z,-l}^-\subseteq\sum_{m_2\geq 0}\hu_{\Z,k_1}^+\hu_{\Z,-l+m_2}^-\huz^0\hu_{\Z,k_2-m_2}^+\subseteq$$
$$\subseteq\sum_{m_1,m_2\geq 0}\hu_{\Z,-l+m_2+m_1}^-\huz^0\hu_{\Z,k_1-m_1}^+\huz^0\hu_{\Z,k_2-m_2}^+=$$
$$=\sum_{m_1,m_2\geq 0}\hu_{\Z,-l+m_2+m_1}^-\huz^0\hu_{\Z,k_1-m_1}^+\hu_{\Z,k_2-m_2}^+\subseteq\sum_{m\geq 0}\hu_{\Z,-l+m}^-\huz^0\hu_{\Z,k-m}^+$$
and symmetrically applying $\Omega$
$$\hu_{\Z,k}^+\hu_{\Z,-l_1}^-\hu_{\Z,-l_2}^-=
\Omega(\hu_{\Z,l_2}^+\hu_{\Z,l_1}^+\hu_{\Z,-k}^-)\subseteq$$
$$\subseteq\Omega (\sum_{m\ge 0}\hu_{\Z,-k+m}^-\hu_{\Z}^0\hu_{\Z,l-m}^+)=\sum_{m\geq 0}\hu_{\Z,-l+m}^-\huz^0\hu_{\Z,k-m}^+.$$
(\ref{uzruzs}) follows from a) and b).

\end{proof}
\end{proposition}

\vskip .3 truecm
\begin{theorem}\label{trm}
\n The $\Z$-subalgebra $\huz$ of $\hu$ generated by $$\{(x_r^{\pm})^{(k)}|r\in\Z,k\in\N\}$$ is an integral form of $\hu$. 

\n More precisely
$$\huz\cong\huz^-\otimes\huz^0\otimes\huz^+\cong\huz^-\otimes\huz^{0,-}\otimes\huz^{0,0}\otimes\huz^{0,+}\otimes\huz^+$$
and a $\Z$-basis of $\huz$ is given by the product
$$B^-B^{0,-}B^{0,0}B^{0,+}B^+$$
where $B^{\pm}$, $B^{0,\pm}$ and $B^{0,0}$ are the $\Z$-bases respectively of $\huz^{\pm}$, $\huz^{0,\pm}$ and $\huz^{0,0}$ given as follows:
$$B^{\pm}=\Big\{{(\bf{x}}^{\pm})^{({\bf{k}})}=\prod_{r\in\Z}(x_r^{\pm})^{(k_r)} \;|\; {\bf{k}}:\Z\to\N\,\, {\rm{is\, finitely\, supported}}\Big\}$$
$$B^{0,\pm}=\Big\{{\hat{{\bf{h}}}_{\pm}^{\bf{k}}}=\prod_{l\in\Z_+}{\hat h_{\pm l}^{k_l}}\;|\; {\bf{k}}:\Z_+\to\N\,\,{\rm{is\, finitely\, supported}}\Big\}$$
$$B^{0,0}=\Big\{{h_0\choose k}{c\choose\tilde k} \;|\;  k,\tilde k\in\N\Big\}.$$
\end{theorem}

\vskip .5 truecm

\section{The integral form of $\gothd$ ($A_2^{(2)}$) }\label{ifa22}
\vskip .5truecm

\n In this section we describe the integral form  $\tuz$ of the enveloping algebra $\tu$ of the Kac-Moody algebra of type $A_2^{(2)}$ generated by the divided powers of the Drinfeld generators $x_r^{\pm}$; unlike the untwisted case, this integral form is strictly smaller than the one (studied in \cite{DM}) generated by the divided powers of the Chevalley generators $e_0$, $e_1$, $f_0$, $f_1$ (see appendix \ref{appendC}).

\n However, the construction of a $\Z$-basis of $\tuz$ follows the idea of the analogous construction in the case $A_1^{(1)}$, seen in the previous section; this method allows us to overcome the technical difficulties arising in case $A_2^{(2)}$ - difficulties which seem otherwise overwhelming.

\n The commutation relations needed to our aim can be partially deduced from the case $A_1^{(1)}$: indeed, underlining some embeddings of $\hat{\goths}$ into $\gothd$ (see remark \ref{emgg}), the commutation relations in $\hu$ can be directly translated into a class of commutation relations in $\tu$ (see corollary \ref{czzp}, proposition \ref{czzq} and the appendix \ref{appendA} for more details).

\n Yet, there are some differences between $A_1^{(1)}$ and $A_2^{(2)}$. 

\n First of all, the real (positive and negative) components of $\tu$ are no more commutative (this is well known: it happens in all the affine cases different from $A_1^{(1)}$, as well as in all the finite cases different from $A_1$), hence the study of their integral form requires some - easy - additional observations (see lemma \ref{zkp}).

\n The non commutativity of the real components of $\tu$ makes the general commutation formula between the exponentials of  positive and negative Drinfeld generators technically more complicated to compute and express than in the case of $\hat\goths$; nevertheless, general and explicit compact formulas can be given in this case, too, always thanks to  the exponential notation. As already seen, the simplification provided by the exponential approach lies essentially on lemma \ref{cle},iv), which allows to perform the computations in $\tu$ reducing to much simpler computations in $\gothd$, and even, thanks to the symmetries highlighted in definition \ref{tto}, in the Lie subalgebra $L=\gothd  \cap({\frak sl_3}\otimes\Q[t])\subseteq\gothd$ (see definition \ref{sottoalgebraL}). 
Recognizing a $\Q[w]$-module structure on each direct summand of $L=L^-\oplus L^0\oplus L^+$ and unifying them in a 
$\Q[w]$-module structure on $L$ (see definition \ref{qwmodulo}) provides a further simplification in the notations: one could have done the same construction for $\hat\goths$, but we have the feeling that in the case of $\hat\goths$ it would be unnecessary and that on the other hand it is useful to present both formulations.

\n The most remarkable difference with respect to $A_1^{(1)}$ on one hand and to Mitzman's integral form on the other hand lies in the description of the generators of the imaginary (positive and negative) components; it can be surprising that they are not what one could expect: $\tuz^{0,+}\neq\Z^{(sym)}[h_r|r>0]$. More precisely (see remark \ref{hdiversi} and theorem \ref{trmA22}) 
$$\tuz^{0,+}\not\subseteq\Z^{(sym)}[h_r|r>0]\ \ {\rm{ and}}\ \  \Z^{(sym)}[h_r|r>0]\not\subseteq\tuz^{0,+};$$
as we shall show, we need to somehow ``deform'' the $h_r$'s (by changhing some of their signs) to get a basis of $\tuz^{0,+}$ by the $(sym)$-construction (see definition \ref{thuz}, example \ref{rvsf} and remark \ref{funtorialita}).

\n Notice that in order to prove that $\tuz$ is an integral form of $\tu$ and that $B$ is a $\Z$-basis of $\tuz$ (theorem \ref{trmA22}) it is not necessary to find explicitly all the commutation formulas between the basis elements. 
In any case, for completeness, we shall collect them in the appendix \ref{appendA}.

\vskip .3 truecm

\begin{definition} \label{a22}

\n ${\gothd}$ (respectively $\tu$) is the Lie algebra (respectively the associative algebra) over $\Q$ generated by $\{c,h_r,x_r^{\pm},X_{2r+1}^{\pm}|r\in\Z\}$ with relations $$c\,\,\,{\rm{is\,\,central}}$$
$$[h_r,h_s]=\delta_{r+s,0}2r(2+(-1)^{r-1})c$$
$$[h_r,x_s^{\pm}]=\pm 2(2+(-1)^{r-1})x_{r+s}^{\pm}$$
$$[h_r,X_s^{\pm}]=\begin{cases}\pm 4X_{r+s}^{\pm}&{\rm{if}}\ 2|r
\\ 0&{\rm{if}}\ 2\not|r
\end{cases}\leqno{(s\ {\rm{odd}})}$$
$$[x_r^{\pm},x_s^{\pm}]=\begin{cases}0&{\rm{if}}\ 2|r+s
\\ \pm(-1)^sX_{r+s}^{\pm}&{\rm{if}}\ 2\not|r+s
\end{cases}$$
$$[x_r^{\pm},X_s^{\pm}]=[X_r^{\pm},X_s^{\pm}]=0$$
$$[x_r^+,x_s^-]=h_{r+s}
+\delta_{r+s,0}rc$$
$$[x_r^{\pm},X_s^{\mp}]=\pm(-1)^r4x_{r+s}^{\mp}\leqno{(s\ {\rm{odd}})}$$
$$[X_r^+,X_s^-]=8h_{r+s}
+4\delta_{r+s,0}rc\leqno{(r,s\ {\rm{odd}})}$$
\end{definition}

\n Notice that $\{x_r^+,x_r^-|r\in\Z\}$ generates $\tu$.

\n Moreover $\{c,h_r,x_r^{\pm},X_{2r+1}^{\pm}|r\in\Z\}$ is a basis of $\gothd$; hence the ordered monomials in these elements (with respect to any total ordering of the basis) is a PBW-basis of $\tu$.

\n $\tu^+$, $\tu^-$, $\tu^0$ are the subalgebras of $\tu$ generated respectively by $$\{x_r^+|r\in\Z\},\ \{x_r^-|r\in\Z\},\ \{c,h_r|r\in\Z\}.$$

\n $\tu^{\pm,0}$, $\tu^{\pm,1}$ and $\tu^{\pm,c}$ are the subalgebras of $\tu^{\pm}$ generated respectively by
$$\{x_r^{\pm}|r\equiv 0\ (mod\ 2)\},\ \{x_r^{\pm}|r\equiv 1\ (mod\ 2)\}\ {\rm{and}}\ \{X_{2r+1}^{\pm}|r\in\Z\}.$$

\n $\tu^{0,+}$, $\tu^{0,-}$, $\tu^{0,0}$, are the subalgebras of $\tu$ (of $\tu^0$) generated respectively by $$\{h_r|r>0\},\ \{h_r|r<0\},\ \{c,h_0\}.$$

\vskip .3 truecm
\n The following remark is a consequence of trivial applications of the PBW-theorem to different subalgebras of $\gothd$.

\begin{remark} \label{tefp}

\n $\tu^+$ and $\tu^-$ are not commutative: $[x_0^+,x_1^+]=-X_1^+$ and $[x_0^-,x_1^-]=X_1^-$.

\n $\tu^{\pm,0}$, $\tu^{\pm,1}$ and $\tu^{\pm,c}$ are (commutative) algebras of polynomials:
$$\tu^{+,0}\cong\Q[x_{2r}^+\; | \; r \in \Z],\ \ 
\tu^{+,1}\cong\Q[x_{2r+1}^+\; | \; r \in \Z],\ \ 
\tu^{+,c}\cong\Q[X_{2r+1}^+\; | \; r \in \Z],$$
$$\tu^{-,0}\cong\Q[x_{2r}^-\; | \; r \in \Z],\ \ 
\tu^{-,1}\cong\Q[x_{2r+1}^-\; | \; r \in \Z],\ \ 
\tu^{-,c}\cong\Q[X_{2r+1} ^-\; | \; r \in \Z].$$
We have the following ``triangular'' decompositions of $\tu^{\pm}$:
$$\tu^{\pm}\cong\tu^{\pm,0}\otimes\tu^{\pm,c}\otimes\tu^{\pm,1}\cong\tu^{\pm,1}\otimes\tu^{\pm,c}\otimes\tu^{\pm,0}$$
Remark that $\tu^{\pm,c}$ is central in $\tu^{\pm}$, so that the images in $\tu^{\pm}$ of $\tu^{\pm,0}\otimes\tu^{\pm,c}$ and $\tu^{\pm,1}\otimes\tu^{\pm,c}$ are commutative subalgebras of 
$\tu$.
\vskip .15 truecm
\n $\tu^0$ is not commutative: $[h_r,h_{-r}]\neq 0$ if $r\neq 0$;

\n $\tu^{0,+}$, $\tu^{0,-}$, $\tu^{0,0}$, are (commutative) algebras of polynomials:
$$\tu^{0,+}\cong\Q[h_r|r>0],\,\,\,\tu^{0,-}\cong\Q[h_r|r<0],\,\,\,\tu^{0,0}\cong\Q[c,h_0];$$
Moreover we have the following triangular decomposition of $\tu^0$:
$$\tu^0\cong\tu^{0,-}\otimes\tu^{0,0}\otimes\tu^{0,+}\cong\tu^{0,+}\otimes\tu^{0,0}\otimes\tu^{0,-}.$$
Remark that $\tu^{0,0}$ is central in $\tu^0$, so that the images in $\tu^0$ of 
$\tu^{0,-}\otimes\tu^{0,0}$ and $\tu^{0,0}\otimes\tu^{0,+}$ are commutative subalgebras of $\tu$.
\vskip .15 truecm
\n Finally remark the triangular decomposition of $\tu$:
$$\tu\cong\tu^-\otimes\tu^0\otimes\tu^+\cong\tu^+\otimes\tu^0\otimes\tu^-,$$
and observe that the images of $\tu^-\otimes\tu^0$ and $\tu^0\otimes\tu^+$ are  subalgebras of $\tu$. 
\end{remark}

\vskip .3 truecm

\begin{definition}\label{tto}
\n $\gothd$ and $\tu$ are endowed with the following anti/auto/ho\-mo/morphisms:

\n $\sigma$ is the antiautomorphism defined on the generators by:
$$x_r^+\mapsto x_r^+,\,\,\,x_r^-\mapsto x_r^-,\,\,\,(\Rightarrow X_r^{\pm}\mapsto-X_r^{\pm},\,\,\,h_r\mapsto-h_r,\,\,\,c\mapsto -c);$$
$\Omega$ is the antiautomorphism defined on the generators by:
$$x_r^+\mapsto x_{-r}^-,\,\,\,x_r^-\mapsto x_{-r}^+,\,\,\,(\Rightarrow X_r^{\pm}\mapsto X_{-r}^{\mp},\,\,\,h_r\mapsto h_{-r},\,\,\,c\mapsto c);$$
\n $T$ is the automorphism defined on the generators by:
$$x_r^+\mapsto x_{r-1}^+,\,\,\,x_r^-\mapsto x_{r+1}^-,\,\,\,(\Rightarrow X_r^{\pm}\mapsto -X_{r\mp 2}^{\pm},\,\,\,h_r\mapsto h_r-\delta_{r,0}c,\,\,\,c\mapsto c);$$
\n for all odd integer $m\in\Z$, $\lambda_m$ is the homomorphism defined on the generators by:
$$x_r^+\mapsto x_{mr}^+,\,\,\,x_r^-\mapsto x_{mr}^-,\,\,\,(\Rightarrow X_r^{\pm}\mapsto X_{mr}^{\pm},\,\,\,h_r\mapsto h_{mr},\,\,\,c\mapsto mc).$$
Remark that if $m$ is even $\lambda_m$ is not defined on $\tu$, but it is still defined on $\tu^{0,+}=\Q[h_r|r>0]$.
\end{definition}

\vskip .3 truecm
\begin{remark}\label{tti}
\n $\sigma^2={\rm{id}}_{\tu}$, $\Omega^2={\rm{id}}_{\tu}$, $T$ is invertible of infinite order;

\n $\lambda_{-1}^2=\lambda_1={\rm{id}}_{\tu}$; $\lambda_m$ is not invertible if $m\neq\pm 1$.

\end{remark}

\vskip .3 truecm
\begin{remark}\label{tct}
$\sigma\Omega=\Omega\sigma$, $\sigma T=T\sigma$, $\Omega T=T\Omega$. 
Moreover for all $m,n$ odd we have $\sigma\lambda_m=\lambda_m\sigma$, $\Omega\lambda_m=\lambda_m\Omega$,
$\lambda_m T^{\pm 1}=T^{\pm m}\lambda_m$, $\lambda_m\lambda_n=\lambda_{mn}$.

\end{remark}

\vskip .3 truecm

\begin{remark}\label{tsb}
\n $\sigma\big|_{\tu^{\pm,0}}={\rm{id}}_{\tu^{\pm,0}},\,\sigma\big|_{\tu^{\pm,1}}={\rm{id}}_{\tu^{\pm,1}},\,\sigma(\tu^{\pm,c})=\tu^{\pm,c},\,\sigma(\tu^{0,\pm})=\tu^{0,\pm}$, $\sigma(\tu^{0,0})=\tu^{0,0}$.

\n $\Omega(\tu^{\pm,0})\!=\tu^{\mp,0}$, $\Omega(\tu^{\pm,1})\!=\tu^{\mp,1}$, $\Omega(\tu^{\pm,c})\!=\tu^{\mp,c}$, $\Omega(\tu^{0,\pm})\!=\tu^{0,\mp}$, $\Omega\big|_{\tu^{0,0}}\!\!=\!{\rm{id}}_{\tu^{0,0}}$.

\n $T
(\tu^{\pm,0})=\tu^{\pm,1}$,
$T
(\tu^{\pm,1})=\tu^{\pm,0}$, $T
(\tu^{\pm,c})=\tu^{\pm,c}$, $T
\big|_{\tu^{0,\pm}}={\rm{id}}_{\tu^{0,\pm}},\,\,\,
T(\tu^{0,0})=\tu^{0,0}$.

\n For all odd $m\in\Z$:

\n $\lambda_m(\tu^{\pm,0})\subseteq\tu^{\pm,0},\,\,\,\lambda_m(\tu^{\pm,1})\subseteq\tu^{\pm,1},\,\,\,\lambda_m(\tu^{\pm,c})\subseteq\tu^{\pm,c},\,\,\,\lambda_m(\tu^{0,0})\subseteq\tu^{0,0}$, $$\lambda_m(\tu^{0,\pm})\subseteq
\begin{cases}\tu^{0,\pm}&{\rm{if}}\,m>0\cr \tu^{0,\mp}&{\rm{if}}\, m<0.
\end{cases}$$

\end{remark}
\vskip .3 truecm

\begin{definition}\label{sottoalgebraL}
$L$, $L^{\pm}$, $L^0$, $L^{\pm,0}$, $L^{\pm,1}$, $L^{\pm,c}$ are the Lie-subalgebras of $\gothd$ generated by:
$$L: \{x_r^+,x_r^-|r\geq 0\},$$
$$L^+: \{x_r^+|r\geq 0\},\ \ L^-: \{x_r^-|r\geq 0\},\ \ L^0: \{h_r|r\geq 0\},$$
$$L^{+,0}: \{x_{2r}^+|r\geq 0\},\ \ L^{+,1}: \{x_{2r+1}^+|r\geq 0\},\ \ L^{+,c}: \{X_{2r+1}^+|r\geq 0\}.$$
$$L^{-,0}: \{x_{2r}^-|r\geq 0\},\ \ L^{-,1}: \{x_{2r+1}^-|r\geq 0\},\ \ L^{-,c}: \{X_{2r+1}^-|r\geq 0\}.$$

\end{definition}
\vskip .3 truecm
\begin{remark}\label{Lbase}
$L^0$, $L^{\pm,0}$, $L^{\pm,1}$ and $L^{\pm,c}$ are commutative Lie-algebras; for these subalgebras of $L$ the Lie-generators given in definition \ref{sottoalgebraL} are bases over $\Q$. 

\n Moreover we have $\Q$-vector space decompositions
$$L=L^-\oplus L^0\oplus L^+,\ \ L^+=L^{+,0}\oplus L^{+,1}\oplus L^{+,c},\ \ L^-=L^{-,0}\oplus L^{-,1}\oplus L^{-,c}.$$
Finally remark that $L^+$ is $T^{-1}$-stable and that $L^-$ is $T$-stable; more in detail $T^{\mp 1}(L^{\pm,0})=L^{\pm,1}$, $T^{\mp 1}(L^{\pm,1})\subseteq L^{\pm,0}$ (so that $L^{\pm,0}$ and $L^{\pm,1}$ are $T^{\mp 2}$-stable); $L^{\pm,c}$ is $T^{\mp1}$-stable.

\end{remark}
\vskip .3 truecm
\begin{definition}\label{qwmodulo}
$L$ is endowed with the $\Q[w]$-module structure defined by $w\big|_{L^-}=T\big|_{L^-}$, $w\big|_{L^+}=T^{-1}\big|_{L^+}$, $w.h_r=h_{r+1}$ $\forall r\in\N$.
\end{definition}

\vskip .3 truecm

\begin{lemma}\label{qwconti}
Let $\xi_1(w),\xi_2(w)\in\Q[w][[u,v]]$. Then:
$$[\xi_1(w^2).x_0^{\pm},\xi_2(w^2).x_1^{\pm}]=\mp(\xi_1\xi_2)(-w).X_1^{\pm};\leqno{i)}$$
$$
[\xi_1(w).x_0^+,\xi_2(w).x_0^-]=(\xi_1\xi_2)(w).h_0;\leqno{ii)}$$
$$[\xi_1(w).x_0^+,\xi_2(w).X_1^-]=4\xi_1(-w)\xi_2(-w^2).x_1^-;\leqno{iii)}$$
$$[\xi_1(w).h_0,\xi_2(w).x_0^{\pm}]=\pm(4\xi_1(w)-2\xi_1(-w))\xi_2(w).x_0^{\pm}.\leqno{iv)}$$
\begin{proof}
The assertions are just a translation of the defining relations of $\tu$:
$$[x_{2r}^{\pm},x_{2s+1}^{\pm}],\ \ [x_r^+,x_s^-],\ \ [x_r^+,X_{2s+1}^-],\ \ [h_r,x_s^{\pm}].$$
For iv), remark that
$$2(2+(-1)^{r-1})w^r=4w^r-2(-w)^r.$$
\end{proof}
\end{lemma}

\vskip .3 truecm
\begin{definition}\label{thuz}
\n Here we define some $\Z$-subalgebras of $\tu$:

\n $\tuz$ is the $\Z$-subalgebra of $\tu
$ generated by $\{(x_r^{+})^{(k)},(x_r^{-})^{(k)}|r\in\Z,k\in\N\}$;

\n $\tuz^+$ and $\tuz^-$ are the $\Z$-subalgebras of $\tu
$ (and of $\tuz$) generated respectively by $\{(x_r^{+})^{(k)}|r\in\Z,k\in\N\}$, and $\{(x_r^{-})^{(k)}|r\in\Z,k\in\N\}$;

\n $\tuz^{\pm,0}=\Z^{(div)}[x_{2r}^{\pm}|r\in \Z]$; 

\n $\tuz^{\pm,1}=\Z^{(div)}[x_{2r+1}^{\pm}|r\in \Z]$; 

\n $\tuz^{\pm,c}=\Z^{(div)}[X_{2r+1}^{\pm}|r\in \Z]$; 

\n $\tuz^{0,0}=\Z^{(bin)}[h_0,c]$;

\n $\tuz^{0,\pm}=\Z^{(sym)}[\varepsilon_rh_{\pm r}|r>0]$ with $\varepsilon_r=\begin{cases}1&{\rm{if}}\ 4\not|r\\-1&{\rm{if}}\ 4|r\end{cases}$;

\n $\tuz^0$ is the $\Z$-subalgebra of $\tu$ generated by $\tuz^{0,-}$, $\tuz^{0,0}$ and $\tuz^{0,+}$.

The notations are those of section \ref{intgpl}.

\n In particular remark the definition of $\tuz^{0,\pm}$ (where the $\varepsilon_r$'s represent the necessary ``deformation'' announced in the introduction of this section, and discussed in details in proposition \ref{convoluzioneintera}) and introduce the notation
$$\Z[\tilde h_k|\pm k>0]=\Z^{(sym)}[\varepsilon_rh_{\pm r}|r>0]$$ where $$\tilde h_{\pm}(u)=\sum_{k\in\N}\tilde h_{\pm k}u^k=\exp\Big(\sum_{r>0}(-1)^{r-1}{\varepsilon_rh_{\pm r}\over r}u^r\Big).$$
\end{definition}
\begin{remark} \label{hdiversi}
It is worth underlining 
that
$\tilde h_{+}(u)\neq\hat h_{+}(u)$, where
$$\Z[\hat h_k|k>0]=\Z^{(sym)}[h_{r}|r>0],$$
that is
$$\hat h_{+}(u)=\sum_{k\in\N}\hat h_{k}u^k=\exp\Big(\sum_{r>0}(-1)^{r-1}{h_{r}\over r}u^r\Big).$$
More precisely the $\Z$-subalgebras generated respectively by $\{\hat h_k|k>0\}$ and $\{\tilde h_k|k>0\}$ are different and not included in each other: indeed $\tilde h_1=\hat h_1$, $\tilde h_2=\hat h_2$, $\tilde h_3=\hat h_3$ but $\hat h_4\not\in\Z[\tilde h_k|k>0]$ and $\tilde h_4\not\in\Z[\hat h_k|k>0]$ (see propositions \ref{convoluzioneintera} and  \ref{emmepiallaerre} and remark \ref{vicelambda}).

\end{remark}
\vskip .3 truecm
\begin{remark}\label{whtilde}
Let $\xi(w)\in\Q[w][[u]]$; the elements $$\exp(\xi(w^2).x_0^{\pm}),\ \ \exp(\xi(w^2).x_1^{\pm})\ \ {\rm{and}}\ \ \exp(\xi(w).X_1^{\pm})$$ lie respectively in $\tuz^{\pm,0}[[u]]$, $\tuz^{\pm,1}[[u]]$ and $\tuz^{\pm,c}[[u]]$ if and only if $\xi(w)$ has integral coefficients, that is if and only if $\xi(w)\in\Z[w][[u]]$ (see example \ref{dvdpw}).

\n Remark also that $$\hat h_+(u)=\exp(\ln(1+wu).h_0),$$ while 
$$\tilde h_+(u)=\exp\left(\big(\ln(1+uw)+{1\over 2}\ln(1-u^4w^4)\big).h_0\right).$$
\end{remark}

\n Before entering the study of the integral forms just introduced, we still dwell on the comparison between $\tilde h_+(u)$ and $\hat h_+(u)$, proving lemma \ref{ometiomecap}, that will be useful later.

\begin{lemma} \label{emme}
For all $m\in\Z\setminus\{0\}$ we have
$$(1+m^2u)^{{1\over m}}\in 1+mu\Z[[u]].$$
\begin{proof}
$(1+\sum_{r> 0}a_ru^r)^m=1+m^2u$ implies 
$$1+m^2u=1+m\sum_{r> 0}a_ru^r+\sum_{k>1}{m\choose k}\big(\sum_{r> 0}a_ru^r\big)^k.$$
Let us prove by induction on $s$ that $a_s\in m\Z$:

if $s=1$ we have that $ma_1=m^2$;

if $s>1$ the coefficient $c_s$ of $u^s$ in $\sum_{k>1}{m\choose k}\big(\sum_{r> 0}a_ru^r\big)^k$ is a combination with integral coefficients of products of the $a_t$'s with $t<s$, which are all multiple of $m$. Then, since $k\geq 2$, $m^2|c_s$.
But $ma_s+c_s=0$, thus $m|a_s$. 

\end{proof}
\end{lemma}

\begin{lemma} \label{ometiomecap}
Let us consider the integral forms
$\Z[\hat h_k|k>0]$ and $\Z[\tilde h_k|k>0]$ of $\Q[h_r|r>0]$  (see example  \ref{rvsf}, formula \ref{dfhp}, definiton \ref{thuz} and remark  \ref{hdiversi}); for all $m>0$ recall the $\Q$-algebra homomorphism $\lambda_m$ of $\Q[h_r|r>0]$ (see proposition \ref{tmom}) and define the analogous homomorphism $\tilde \lambda_m$ mapping each $\varepsilon_rh_r$ to $\varepsilon_{mr}h_{mr}$ (of course $\Z[\tilde h_k|k>0]$ is $\tilde\lambda_m$-stable $\forall m>0$).

\n We have that:

\n i) if $m$ is odd then $\tilde\lambda_m=\lambda_m$; in particular $\Z[\tilde h_k|k>0]$ is $\lambda_m$-stable;

\n ii) $\lambda_2(\hat h_k)\in\Z[\tilde h_l|l>0]$ for all $k>0$;

\n iii) $\hat h_+(4u)^{{1\over 2}}\in\Z[\tilde h_k|k>0][[u]]$;
\begin{proof}

\n i) If $m$ is odd then $4|mr\Leftrightarrow 4|r$, hence $\varepsilon_{mr}=\varepsilon_r$ $\forall r>0$ and the claim follows from proposition \ref{tmom}

\n ii) By proposition \ref{tmom} we know that $\Z[\tilde h_k|k>0]$ is $\tilde\lambda_2$-stable; but 
$$\tilde\lambda_2(\tilde h_+(u^2))=\exp\sum_{r>0}(-1)^{r-1}{\varepsilon_{2r}h_{2r}\over r}u^{2r}=
\exp\sum_{r>0}{h_{2r}\over r}u^{2r}=\lambda_2(\hat h_+(-u^2))^{-1};$$ equivalently
$$\lambda_2(\hat h_+(u^2))=\tilde\lambda_2(\tilde h_+(-u^2))^{-1},$$
which implies the claim.

\n iii) Remark that $$\hat h_+(u)\tilde h(u)_+^{-1}=
\exp\left(-\sum_{r>0}{2h_{4r}\over 4r}u^{4r}\right)=
\tilde\lambda_4(\tilde h_+(-u^4))^{-{1\over 2}};$$
then
$$\hat h_+(4u)^{{1\over 2}}=\tilde h_+(4u)^{{1\over 2}}\tilde\lambda_4(\tilde h_+(-4^4u^4))^{-{1\over 4}}.$$
Since $\tilde h_+(4u)\in 1+4u\Z[\tilde h_k | k>0][[u]]$ and $\tilde\lambda_4(\tilde h_+(4^4u^4))\in 1+4^4u\Z[\tilde h_k| k>0][[u]]$
we deduce from lemma \ref{emme} that 
$$\tilde h_+(4u)^{{1\over 2}},\ \ \tilde\lambda_4(\tilde h_+(4^4u^4))^{{1\over 4}}\in\Z[\tilde h_k | k>0],$$ which implies the claim.

\end{proof}
\end{lemma}

\vskip .3 truecm
\begin{remark}\label{tzp}

\n It is obvious that
$\tuz^{\pm,0}$, $\tuz^{\pm,1}$, $\tuz^{\pm,c}$, $\tuz^{0,\pm}$ and $\tuz^{0,0}$ are integral forms respectively of
$\tu^{\pm,0}$, $\tu^{\pm,1}$, $\tu^{\pm,c}$, $\tu^{0,\pm}$ and $\tu^{0,0}$.

\n Hence by the commutativity properties we also have that $\tuz^{\pm,0}\tuz^{\pm,c}$ and $\tuz^{\pm,1}\tuz^{\pm,c}$ are integral forms respectively of $\tu^{\pm,0}\tu^{\pm,c}$ and $\tu^{\pm,c}\tu^{\pm,1}$.

\n Analogously $\tuz^{0,0}\tuz^{0,+}$ and $\tuz^{0,-}\tuz^{0,0}$ are integral forms respectively of $\tu^{0,0}\tu^{0,+}$ and $\tu^{0,-}\tu^{0,0}$.

\n We want to prove that:

\n 1) $\tuz^0=\tuz^{0,-}\tuz^{0,0}\tuz^{0,+}$, so that $\tuz^0$ is an integral form of $\tu^0$;

\n 2) $\tuz^{\pm}=\tuz^{\pm,1}\tuz^{\pm,c}\tuz^{\pm,0}$, so that $\tuz^+$ and $\tuz^-$ are integral forms respectively of $\tu^+$ and $\tu^-$;

\n 3) $\tuz=\tuz^-\tuz^0\tuz^+$, so that $\tuz$ is an integral form of $\tu$.
\vskip .3 truecm

\n It is useful to evidentiate the behaviour of the $\Z$-subalgebras introduced above under the symmetries of $\tu$.

\end{remark}

\vskip .3 truecm
\begin{proposition}\label{sttuz}
\n The following stability properties under the action of $\sigma$, $\Omega$, $T^{\pm 1}$ and $\lambda_m$ ($m \in \Z$ odd) hold:

\n i) $\tuz$, $\tuz^+$ and $\tuz^-$ are $\sigma$-stable, $T^{\pm 1}$-stable, $\lambda_m$-stable.

\n \,\,\,\,\,
$\tuz$ is also $\Omega$-stable, while $\Omega(\tuz^{\pm})=\tuz^{\mp}$.

\n ii) $\tuz^{+,0}$, $\tuz^{+,1}$ and $\tuz^{+,c}$ are $\sigma$-stable, $T^{\pm 2}$-stable, $\lambda_m$-stable.

\n \,\,\,\,\,
$\tuz^{+,c}$ is also $T^{\pm 1}$-stable, while $T^{\pm 1}(\tuz^{+,0})=\tuz^{+,1}$.

\n \,\,\,\,\,
$\Omega(\tuz^{+,0})=\tuz^{-,0}$, 
$\Omega(\tuz^{+,1})=\tuz^{-,1}$ and
$\Omega(\tuz^{+,c})=\tuz^{-,c}$.

\n iii) $\tuz^{0,0}$, $\tuz^{0,+}$ and $\tuz^{0,-}$ are $\sigma$-stable and $T^{\pm 1}$-stable.

\n \,\,\,\,\,
$\tuz^{0,0}$ is also $\Omega$-stable and $\lambda_m$-stable; $\Omega(\tuz^{0,\pm})=\tuz^{0,\mp}$; $\tuz^{0,\pm}$ is $\lambda_m$-stable if $m>0$, while $\lambda_m(\tuz^{0,\pm})\subseteq\tuz^{0,\mp}$ if $m<0$.

\n iv)
$\tuz^0$ is $\sigma$-stable, $\Omega$-stable, $T^{\pm 1}$-stable, $\lambda_m$-stable.
\begin{proof}
The only non-trivial assertion is the claim that $\tuz^{0,+}$ is $\lambda_m$-stable when $m>0$, which was proved in lemma \ref{ometiomecap},i).

\n The assertion about $\lambda_m(\tuz^{0,\pm})$ in the general case follows using that $$\Omega(\tuz^{0,\pm})=\tuz^{0,\mp}=\lambda_{-1}(\tuz^{0,\pm}),\ \ \lambda_m\Omega=\Omega\lambda_m\ \ {\rm{and}}\ \  \lambda_{-m}=\lambda_{-1}\lambda_m.$$
Remark that
$$\sigma(\tilde h_{\pm}(u))\!=\!\tilde h_{\pm}(u)^{-1},\,\Omega(\tilde h_{\pm}(u))\!=\!\lambda_{-1}(\tilde h_{\pm}(u))\!=\!\tilde h_{\mp}(u),\,T^{\pm 1}(\tilde h_{\pm}(u))\!=\!\tilde h_{\pm}(u).$$
\end{proof}
\end{proposition}

\vskip .3 truecm

\begin{remark} \label{autstab}
 
\n The stability properties described in proposition \ref{sttuz} imply that:

\n i) $\sigma(\tuz^{0,-}\tuz^{0,0}\tuz^{0,+})=\tuz^{0,+}\tuz^{0,0}\tuz^{0,-}$; in particular
$$\tuz^0=\tuz^{0,-}\tuz^{0,0}\tuz^{0,+}\Leftrightarrow\tuz^0=\tuz^{0,+}\tuz^{0,0}\tuz^{0,-}.$$

\n ii) $T^{\pm 1}(\tuz^{+,1}\tuz^{+,c}\tuz^{+,0})=\tuz^{+,0}\tuz^{+,c}\tuz^{+,1}$ and $\tuz^{+,1}\tuz^{+,c}\tuz^{+,0}$ is $T^{\pm 2}$-stable and $\lambda_m$-stable ($m\in\Z$ odd); in particular:

$$\tuz^+=\tuz^{+,1}\tuz^{+,c}\tuz^{+,0}\Leftrightarrow\tuz^+=\tuz^{+,0}\tuz^{+,c}\tuz^{+,1}.$$

\n iii) $\tuz^0\tuz^+$ is $T^{\pm 1}$-stable and $\lambda_{-1}$-stable, and $\Omega(\tuz^0\tuz^+)=\tuz^-\tuz^0$; in particular
it is enough to prove that $(x_0^+)^{(k)}\tilde h_+(u)\in\tilde h_+(u)\tuz^+[[u]]$ $\forall k\geq 0$ in order to show that 
$$(x_r^+)^{(k)}\tilde h_{\pm}(u)\in\tilde h_{\pm}(u)\tuz^+[[u]], \tilde h_{\pm}(u)(x_r^-)^{(k)}\in\tuz^-[[u]]\tilde h_{\pm}(u)\ \ \forall r\in Z,\ k\in\N,$$
or equivalently that $\tuz^+\tuz^0\subseteq\tuz^0\tuz^+$ and $\tuz^0\tuz^-\subseteq\tuz^-\tuz^0$.

\n iv) $\tuz^-\tuz^0\tuz^+$ is $T^{\pm 1}$-stable and $\lambda_m$-stable ($m\in\Z$ odd); in particular if one shows that 
$(x_0^+)^{(k)}(x_1^-)^{(l)}\in\tuz^-\tuz^0\tuz^+$ it follows that
$\forall r,s\in\Z$ such that $2\not|r+s$ $$(x_r^+)^{(k)}(x_s^-)^{(l)}=T^{-r}\lambda_{r+s}((x_0^+)^{(k)}(x_1^-)^{(l)})\in\tuz^-\tuz^0\tuz^+.$$

\end{remark}

\vskip.3 truecm

\begin{proposition} \label{zkd}

\n The following identities hold in $\tu$:

$$\hat h_+(u)\hat h_-(v)=\hat h_-(v)(1-uv)^{-4c}(1+uv)^{2c}\hat h_+(u)$$
and
$$\tilde h_+(u)\tilde h_-(v)=\tilde h_-(v)(1-uv)^{-4c}(1+uv)^{2c}\tilde h_+(u).$$
In particular $\tuz^0=\tuz^{0,-}\tuz^{0,0}\tuz^{0,+}$ and $\tuz^0$ is an integral form of $\tu^0$.
\begin{proof}

\n Since $[h_r,h_s]=[\varepsilon_rh_r,\varepsilon_sh_s]=\delta_{r+s,0}2r(2+(-1)^{r-1})c$,
the claim is proposition \ref{heise} with $m=4$, $l=-2$.

\end{proof}
\end{proposition}

\begin{lemma} \label{zkp}

\n The following identity holds in $\tu$ for all $r,s\in\Z$:

$$\exp(x_{2r}^+u)\exp(x_{2s+1}^+v)=\exp(x_{2s+1}^+v)\exp(-X_{2r+2s+1}^+uv)\exp(x_{2r}^+u).$$
\begin{proof}
The claim is an immediate consequence of lemma \ref{cle},vii), thanks to the relation $[x_{2r}^+,x_{2s+1}^+]=-X_{2r+2s+1}^+$.
\end{proof}
\end{lemma}

\vskip.3 truecm
\begin{corollary} \label{zppkd}
$\tuz^+=\tuz^{+,1}\tuz^{+,c}\tuz^{+,0}$; then $\tuz^{\pm}$ is an integral form of $\tu^{\pm}$.
\begin{proof}
\n From lemma \ref{zkp} we deduce that:

\n 
i) $(X_{2r+1}^+)^{(k)}\in\tuz^+$ $\forall k\in\N,r\in\Z$; this implies that $$\tuz^{+,c}\subseteq\tuz^+\ \ {\rm{and}}\ \ \tuz^{+,1}\tuz^{+,c}\tuz^{+,0}\subseteq\tuz^+.$$

\n 
ii) $\tuz^{+,0}\tuz^{+,1}\subseteq\tuz^{+,1}\tuz^{+,c}\tuz^{+,0}$, hence $\tuz^{+,1}\tuz^{+,c}\tuz^{+,0}$ is stable by left multiplication by $\tuz^{+,0}$, hence by $\tuz$ (which is generated by $\tuz^{+,0}$ and $\tuz^{+,1}$).

\n Since $1\in\tuz^{+,1}\tuz^{+,c}\tuz^{+,0}$, we deduce $\tuz^+\subseteq\tuz^{+,1}\tuz^{+,c}\tuz^{+,0}$, and the claim follows applying $\Omega$ (see proposition \ref{sttuz},i)).

\end{proof}
\end{corollary}
 
\begin{proposition} \label{xtuzzero}
$\tuz^+\tuz^{0,0}\subseteq\tuz^{0,0}\tuz^+$; more precisely
$$(x_r^+)^{(k)}{h_0\choose l}={h_0-2k\choose l}(x_r^+)^{(k)}\ \ \forall r\in\Z,\ k,l\in\N.$$
\begin{proof}
The claim follows by immediate application of \ref{fru}.\end{proof}
\end{proposition}
\vskip .3 truecm

\begin{proposition} \label{zkopp}
In $\tu$ the following holds:

\n i) $x_0^+\tilde h_+(u)=\tilde h_+(u)(1-uT^{-1})^6(1-u^2T^{-2})^{-3}(1+u^2T^{-2})(x_0^+)$;

\n ii) $(x_0^+)^{(k)}\tilde h_+(u)\in\tilde h_+(u)\tuz^+[[u]]$ $\forall k\in\N$;

\n iii) $\tuz^+\tuz^{0,+}\subseteq\tuz^{0,+}\tuz^+$.
\begin{proof}
i) 
We have that
$[\varepsilon_r h_r, x_0^+]=\varepsilon_r2(2+(-1)^{r-1})x_r^+$
and $$\varepsilon_r2(2+(-1)^{r-1})=\begin{cases}6&{\rm{if}}\ 2\not|r\\2=6-4&{\rm{if}}\ 2|r\ {\rm{and}}\ 4\not|r\\-2=6-4-4&{\rm{if}}\ 4|r,\end{cases}$$
hence
proposition \ref{hh} applies, with $m_1=6$, $m_2=-2$, $m_4=-1$ and implies that 
$$x_0^+\tilde h_+(u)=\tilde h_+(u)(1+uT^{-1})^{-6}(1-u^2T^{-2})^{2}(1-u^4T^{-4})(x_0^+)=$$
$$=\tilde h_+(u)(1-uT^{-1})^6(1-u^2T^{-2})^{-3}(1+u^2T^{-2})(x_0^+).$$

\n ii) Let us underline that $(1-u^2)^{-3}(1+u^2)\in\Z[[u^2]]$, hence from the coefficients of $(1-u)^6$ it can be deduced that
$$(1-u)^6(1-u^2)^{-3}(1+u^2)\in\Z[[u^2]]+2u\Z[[u^2]]$$
and
$$x_0^+\tilde h_+(u)=\tilde h_+(u)\sum_{r\geq 0}a_rx_r^+u^r\ \ {\rm{with}}\ \ a_r\in\Z\ \forall r\geq 0\ {\rm{and}}\ 2|a_r\ \forall r\ {\rm{odd}}.$$
If we define $y_0=\sum_{r\ge0}a_{2r}x_{2r}^+u^{2r}$, $y_1={1\over 2}\sum_{r	\ge0}a_{2r+1}x_{2r+1}^+u^{2r+1}$ we have that, thanks to lemma \ref{cle},viii)
$$\exp(x_0^+v)\tilde h_+(u)=\tilde h_+(u)\exp((y_0+2y_1)v)=$$
$$=\tilde h_+(u)\exp(2y_1v)\exp([y_0,y_1]v^2)\exp(y_0v)\in\tilde h_+(u)\tuz^+[[u,v]]$$
from which the claim follows thanks to remark \ref{whtilde}.

\n iii) From the $T^{\pm 1}$-stability of $\tuz^+$ and the fact that $T^{\pm 1}\big|_{\tuz^{0,+}}=id$ we deduce that for all $r\in\Z,\ k\in\N$
$$(x_r^+)^{(k)}\tilde h_+(u)\in\tilde h_+(u)\tuz^+[[u]].$$
The claim follows recalling that the $(x_r^+)^{(k)}$'s generate $\tuz^+$ and the $\tilde h_k$'s generate $\tuz^{0,+}$.
\end{proof}
\end{proposition}

\begin{corollary} \label{czzp2}
$\tuz^{\pm}\tuz^0=\tuz^0\tuz^{\pm}$.
In particular $\tuz^0\tuz^+$ and $\tuz^-\tuz^0$ are 
subalgebras of $\tuz$.
\begin{proof}
$\tuz^+\tuz^{0,0}\subseteq \tuz^{0,0}\tuz^+$ (see proposition \ref{xtuzzero}) and $\tuz^+\tuz^{0,+}\subseteq \tuz^{0,+}\tuz^+$ (see \ref{zkopp},iii)); moreover
$$\tuz^+\tuz^{0,-}=\lambda_{-1}(\tuz^+\tuz^{0,+})\subseteq\lambda_{-1}(\tuz^{0,+}\tuz^+)=\tuz^{0,-}\tuz^+.$$
Hence $\tuz^+\tuz^0\subseteq\tuz^0\tuz^+$.

\n Applying $\sigma$ we get the reverse inclusion and applying $\Omega$ we obtain the claim for $\tuz^-$.

\end{proof}
\end{corollary}

\n Now that we have described $\tuz^0$, $\tuz^{\pm}$ and the $\Z$-subalgebras generated by $\tuz^0$ and $\tuz^+$ (respectively by $\tuz^0$ and $\tuz^-$), in order to show that $\tuz=
\tuz^-\tuz^0\tuz^+$ it remains to prove that $$\tuz^0\subseteq\tuz\ \ {\rm{and}}\ \ \tuz^+\tuz^-\subseteq
\tuz^-\tuz^0\tuz^+.$$
Before attaching this problem in its generality it is worth evidentiating the existence of some copies of $\hat\goths$ inside $\gothd$, hence of embeddings $\hu\hookrightarrow\tu$, that induce some useful commutation relations in $\tu$.

\vskip .3 truecm

\begin{remark}\label{emgg}
The $\Q$-linear maps $f,F:\hat\goths\to\gothd$ defined by
$$x_r^{\pm}\mapsto x_{2r}^{\pm},\ \ h_r\mapsto h_{2r},\ \ c\mapsto 2c\leqno{f:}$$
$$x_r^{\pm}\mapsto {X_{2r\mp 1}^{\pm}\over 4}
,\ \ h_r\mapsto {h_{2r}\over 2}-\delta_{r,0}{c\over 4},\ \ c\mapsto {c\over 2}\leqno{F:}$$
are Lie-algebra homomorphisms, obviously injective, inducing embeddings $f,F:\hu\hookrightarrow\tu$.

\end{remark}
\begin{corollary}\label{czzp}
$f(\huz^{0,0})\subseteq\tuz^{0,0}\subseteq\tuz$.
\begin{proof}
Since $f(\huz^{\pm})\subseteq\tuz^{\pm,0}\subseteq\tuz$ we have that  $f$ maps $\huz$ (which is generated by $\huz^+$ and $\huz^-$) into $\tuz$; in particular
$f(\huz^{0,0})\subseteq\tuz$. But 
$$f(\huz^{0,0})=f(\Z^{(bin)}[h_0,c])=\Z^{(bin)}[h_0,2c],$$ thus $\Z^{(bin)}[h_0,2c]\subseteq\tuz.$
Since $\tuz$ is $T$-stable and $T(h_0)=h_0-c$ we also have $\Z^{(bin)}[h_0-c]\subseteq\tuz$, so that
$$f(\huz^{0,0})=\Z^{(bin)}[h_0,2c]\subseteq\Z^{(bin)}[h_0,c]=\Z^{(bin)}[h_0,h_0-c]\subseteq\tuz$$
which is the claim because $\tuz^{0,0}=\Z^{(bin)}[h_0,c]$.
\end{proof}
\end{corollary}

\begin{proposition} \label{czzq}
$\tuz^{+,0}\tuz^{-,0}\subseteq\tuz^-\tuz^0\tuz^+$ and $\tuz^{+,1}\tuz^{-,1}\subseteq\tuz^-\tuz^0\tuz^+$.
\begin{proof}
$\tuz^{+,0}\tuz^{-,0}=f(\huz^+\huz^-)\subseteq f(\huz^-\huz^0\huz^+)=\tuz^{-,0}f(\huz^0)\tuz^{+,0}$: we want to prove that 
$f(\huz^0)=f(\huz^{0,-}\huz^{0,0}\huz^{0,+})\subseteq\tuz^0$. 

\n By corollary \ref{czzp} $f(\huz^{0,0})\subseteq\tuz^{0,0}$.

\n On the other hand
$$f(\huz^{0,+})=f(\Z^{(sym)}[h_r|r>0])=\Z^{(sym)}[h_{2r}|r>0]=\lambda_2(\Z[\hat h_k|k>0]),$$
hence $f(\huz^{0,+})\subseteq\Z[\tilde h_k|k>0]=\tuz^{0,+}$ thanks to lemma \ref{ometiomecap} ii).

\n Finally remark that $f\Omega=\Omega f$, thus $f(\huz^{0,-})=f\Omega(\huz^{0,+})\subseteq\Omega\tuz^{0,+}\subseteq\tuz^{0,-}$.

\n It follows that $f(\huz^0)\subseteq\tuz^0$ and $\tuz^{+,0}\tuz^{-,0}\subseteq\tuz^-\tuz^0\tuz^+$.

\n The assertion for $\tuz^{\pm,1}$ follows applying $T$, see proposition \ref{sttuz},i),ii) and iv).

\end{proof}
\end{proposition}

\subsection{$\exp(x_0^+u)\exp(x_1^-v)$ and $\tuz^{0,+}$: here comes the hard work}\label{sottosezione}

\n We shall deal with the commutation between $\tuz^{+,0}$ and $\tuz^{-,1}$ following the strategy already proposed for $\huz$ and recalling remark \ref{autstab},iv): finding an explicit expression involving suitable exponentials for 
$$\exp(x_0^+u)\exp(x_1^-v)\in\tu^{-,1}\tu^{-,c}\tu^{-,0}\tu^{0,+}\tu^{+,1}\tu^{+,c}\tu^{+,0}[[u,v]]$$
and proving that all its coefficients lie in $$\tuz^{-,1}\tuz^{-,c}\tuz^{-,0}\tuz^{0,+}\tuz^{+,1}\tuz^{+,c}\tuz^{+,0}\subseteq\tuz^-\tuz^0\tuz^+.$$

\n Since here there are more factors involved, the computation is more complicated than in the case of $\hat\goths$ and the simplification provided by this approach is even more evident. On the other hand it is not immediately clear from the commutation formula that our element belongs to $\tuz^-\tuz^0\tuz^+$, or better: the factors relative to the (negative, resp.  positive) real root vectors will be evidently elements of $\tuz^-$, resp. $\tuz^+$, while proving that the null part lies indeed in $\tuz^0$ is not evident at all and will require a deeper inspection (see remark \ref{hhdehh}, lemma \ref{contidp} and corollary \ref{hcappucciod}).

\n As we shall see, in order to complete the proof that $\tuz^{0,+}\subseteq\tuz$ (see proposition \ref{samealg}), it is useful to compute also $\exp(x_0^+u)\exp(X_1^-v)$. The two computations ($\exp(x_0^+u)\exp(yv)$ with $y=x_1^-$ or $y=X_1^-$) are essentially the same and will be performed together (see the considerations from remark \ref{gesponenziale} to lemma \ref{wderivcomm}, of which the propositions \ref{xmenogrande} and \ref{x0piux1meno} are straightforward applications); even though $\exp(x_0^+u)\exp(x_1^-v)$ presents more symmetries than $\exp(x_0^+u)\exp(X_1^-v)$ (see remark \ref{gtg},iii)), its interpretation will require more work, since it is not evident the connection with $\tuz^{0,+}$, as just mentioned.

\vskip .3 truecm

\begin{remark}\label{gesponenziale}
Let $G=G(u,v)\in\tu[[u,v]]$ and $y\in L^-$ (see definition \ref{sottoalgebraL}); then $$G(u,v)=\exp(x_0^+u)\exp(yv)$$ if and only if the following two conditions hold:

a) $G(0,v)=\exp(yv)$;

b) ${d\over du}G(u,v)=x_0^+G(u,v)$.

\end{remark}

\vskip .3 truecm
\begin{notation}\label{gabg}
In the following $G^-$, $G^0$, $G^+$ will denote elements of $\tu[[u,v]]$ of the form
$$G^-=\exp({\alpha_-})\exp({\beta_-})\exp({\gamma_-}),$$
$$G^+=\exp(\gamma_+)\exp(\beta_+)\exp(\alpha_+),$$
$$G^0=\exp(\eta)$$
with
$$\alpha_-\in\Q[w^2][[u,v]].x_1^-,\ \beta_-\in\Q[w][[u,v]].X_1^-,\ \gamma_-\in\Q[w^2][[u,v]].x_0^-,$$
$$\alpha_+\in\Q[w^2][[u,v]].x_0^+,\ \beta_+\in\Q[w][[u,v]].X_1^+,\ \gamma_+\in\Q[w^2][[u,v]].x_1^+,$$
$$\eta\in w\Q[w][[u,v]].h_0.$$
$G(u,v)$ will denote the element $G(u,v)=G=G^-G^0G^+$.

\end{notation}
\vskip .3 truecm
\begin{remark}\label{gtg}
Let $G=G^-G^0G^+\in\tu[[u,v]]$ be as in notation \ref{gabg}. Then:

\n i) Of course $${dG\over du}={dG^-\over du}G^0G^++G^-{dG^0\over du}G^++G^-G^0{dG^+\over du}$$
where, considering the commutativity properties, we have that
$${dG^-\over du}=\exp({\alpha_-})\exp({\beta_-}){d(\alpha_-+\beta_-+\gamma_-)\over du}\exp({\gamma_-}),$$
$${dG^+\over du}=\exp(\gamma_+){d(\alpha_++\beta_++\gamma_+)\over du}\exp(\beta_+)\exp(\alpha_+),$$
$${dG^0\over du}={d\eta\over du}G^0
.$$

\n ii) If moreover $G=\exp(x_0^+u)\exp(yv)$ with $y\in L^-$, the property b) of remark \ref{gesponenziale} translates into

$$x_0^+G=\exp({\alpha_-})\exp({\beta_-}){d(\alpha_-+\beta_-+\gamma_-)\over du}\exp({\gamma_-})G^0G^++$$
$$+G^-{d\eta\over du}G^0G^++G^-G^0\exp(\gamma_+){d(\alpha_++\beta_++\gamma_+)\over du}\exp(\beta_+)\exp(\alpha_+).$$

\n iii) If in particular $y=x_1^-$, then $T\lambda_{-1}\Omega(G(u,v))=G(v,u)$; hence $$G^-(u,v)=T\lambda_{-1}\Omega(G^+)(v,u),$$
$$\alpha_-(u,v)=T\lambda_{-1}\Omega(\alpha_+)(v,u),$$
$$\beta_-(u,v)=T\lambda_{-1}\Omega(\beta_+)(v,u),$$
$$\gamma_-(u,v)=T\lambda_{-1}\Omega(\gamma_+)(v,u),$$
$$\eta(u,v)=\eta(v,u).$$
Observe that $T\lambda_{-1}\Omega(X_{2r+1}^+)=-X_{2r+3}^-$ $\forall r\in\Z$.
\end{remark}
\vskip .3 truecm

\n The following lemma is based on lemma \ref{cle}, iv) and on the defining relations of $\tu$ (definition \ref{a22}).
\begin{lemma}\label{derivcomm}
With the notations fixed in \ref{gabg} we have that:
$$x_0^+\exp(\alpha_-)=\leqno{i)}$$
$$=\exp(\alpha_-)\left(x_0^++[x_0^+,\alpha_-]+{1\over 2}[[x_0^+,\alpha_-],\alpha_-]+{1\over 6}[[[x_0^+,\alpha_-],\alpha_-],\alpha_-]\right);$$
$$x_0^+\exp(\alpha_-)\exp(\beta_-)=\exp(\alpha_-)\exp(\beta_-)\cdot\leqno{ii)}$$ $$
\cdot\left(x_0^++[x_0^+,\alpha_-]+{1\over 2}[[x_0^+,\alpha_-],\alpha_-]+{1\over 6}[[[x_0^+,\alpha_-],\alpha_-],\alpha_-]+[x_0^+,\beta_-]\right);$$
$$(x_0^++[x_0^+,\alpha_-])\exp(\gamma_-)=\leqno{iii)}$$
$$=\exp(\gamma_-)\left(x_0^++[x_0^+,\alpha_-]+[x_0^+,\gamma_-]\right)+$$
$$+\left([[x_0^+,\alpha_-],\gamma_-]+{1\over 2}[[x_0^+,\gamma_-],\gamma_-]-{1\over 2}[[[x_0^+,\alpha_-],\gamma_-],\gamma_-]\right)\exp(\gamma_-);$$
iv) $x_0^+\exp(\eta)=\exp(\eta)(y_0+y_1)$ with $$y_0\in\Q[w^2][[u,v]].x_0^+,\ \ y_1\in w\Q[w^2][[u,v]].x_0^+;$$ 
$$(y_0+y_1)\exp(\gamma_+)=\exp(\gamma_+)(y_0+y_1+[y_0,\gamma_+]).\leqno{v)}$$
vi) In conclusion
$$x_0^+G={dG\over du}$$ if and only if the following relations hold:
$${d\alpha_-\over du}=[x_0^+,\beta_-]+[[x_0^+,\alpha_-],\gamma_-]$$
$${d\beta_-\over du}={1\over 6}[[[x_0^+,\alpha_-],\alpha_-],\alpha_-]-{1\over 2}[[[x_0^+,\alpha_-],\gamma_-],\gamma_-]$$
$${d\gamma_-\over du}={1\over 2}[[x_0^+,\alpha_-],\alpha_-]+{1\over 2}[[x_0^+,\gamma_-],\gamma_-]$$
$${d\eta\over du}=[x_0^+,\gamma_-]+[x_0^+,\alpha_-]$$
$${d\alpha_+\over du}=y_0$$
$${d\beta_+\over du}=[y_0,\gamma_+]$$
$${d\gamma_+\over du}=y_1.$$
\begin{proof}
i)-v) are straightforward repeated applications of lemma \ref{cle},iv) remarking that:

\n i) and ii): $[[[x_0^+,\alpha_-],\alpha_-],\alpha_-]\in\tu^{-,c}[[u,v]]$, hence it commutes with both $\alpha_-$ and $\beta_-$ (which are in $\tu^-[[u,v]]$);

\n ii): $\beta_-\in\tu^{-,c}[[u,v]]$, hence it commutes also with $[[x_0^+,\alpha_-],\alpha_-]$ and $[x_0^+,\beta_-]$ (which belong to $\tu^-[[u,v]]$) and with $[x_0^+,\alpha_-]$ (because $[h_{2r+1},\tu^{-,c}]=0$ $\forall r\in\Z$);

\n iii): $[[x_0^+,\gamma_-],\gamma_-]$ and $[[[x_0^+,\alpha_-],\gamma_-],\gamma_-]$ belong respectively to $\tu^{-,0}[[u,v]]$ and $\tu^{-,c}[[u,v]]$, so that they commute with $\gamma_-\in\tu^{-,0}[[u,v]]$; the claim follows from the identities
$$(x_0^++[x_0^+,\alpha_-])\exp(\gamma_-)=\exp(\gamma_-)\cdot
\Big(x_0^++[x_0^+,\alpha_-]+$$
$$+[x_0^+,\gamma_-]+[[x_0^+,\alpha_-],\gamma_-]+{1\over 2}[[x_0^+,\gamma_-],\gamma_-]+{1\over 2}[[[x_0^+,\alpha_-],\gamma_-],\gamma_-]\Big)$$
and 
$$\exp(\gamma_-)[[x_0^+,\alpha_-],\gamma_-]=([[x_0^+,\alpha_-],\gamma_-]-[[[x_0^+,\alpha_-],\gamma_-],\gamma_-])\exp(\gamma_-);$$
iv): lemma \ref{lhlh} implies that $\exp(\eta)^{-1}x_0^+\exp(\eta)\in\Q[w][[u,v]].x_0^+;$

\n v): $\gamma_+\in\tu^{+,1}[[u,v]]$ commutes with both $y_1\in\tu^{+,1}[[u,v]]$ and $[y_0,\gamma_+]\in\tu^{+,c}[[u,v]]$.

\n Point vi) is a consequence of points i)-v) and remark \ref{gtg},i).

\end{proof}
\end{lemma}

\begin{lemma}\label{wderivcomm}
By abuse of notation let $\alpha_{\pm}$, $\beta_{\pm}$, $\gamma_{\pm}$, $\eta$ and $y_0$ (see notation \ref{gabg} and lemma \ref{derivcomm},iv)) denote also the elements of $\Q[w][[u,v]]$ such that
$$\alpha_+=\alpha_+(w^2).x_0^+,\ \ \beta_+=\beta_+(w).X_1^+,\ \ \gamma_+=\gamma_+(w^2).x_1^+,$$
$$\alpha_-=\alpha_-(w^2).x_1^-,\ \ \beta_-=\beta_-(w).X_1^-,\ \ \gamma_-=\gamma_-(w^2).x_0^-,$$
$$\eta=\eta(w).h_0.$$
Then the relations of lemma \ref{derivcomm},vi) become:
$${d\alpha_-(w^2)\over du}=4\beta_-(-w^2)-6\alpha_-(w^2)\gamma_-(w^2),$$
$${d\beta_-(w)\over du}=\alpha_-(-w)(w\alpha_-^2(-w)-3\gamma_-^2(-w)),$$
$${d\gamma_-(w^2)\over du}=-3w^2\alpha_-^2(w^2)-\gamma_-^2(w^2),$$
$${d\eta(w)\over du}=w\alpha_-(w^2)+\gamma_-(w^2),$$
$${d(\alpha_+(w^2)+w\gamma_+(w^2))\over du}=\exp(-4\eta(w)+2\eta(-w)),$$
$${d\beta_+(w)\over du}=-{d\alpha_+(-w)\over du}\gamma_+(-w).$$
\begin{proof}
The claim is obtained using lemma \ref{qwconti}.\end{proof}
\end{lemma}

\begin{proposition}\label{xmenogrande}
$$\exp(x_0^+u)\exp(X_1^-v)=$$
$$=\exp({\alpha_-})\exp({\beta_-})\exp({\gamma_-})\exp(\eta)\exp(\gamma_+)\exp(\beta_+)\exp(\alpha_+)$$
where, with the notations of lemma \ref{wderivcomm},
$$\alpha_-(w)={4uv\over 1-4^2wu^4v^2},\ \ \ \ \alpha_+(w)={u\over 1-4^2wu^4v^2},$$
$$\beta_-(w)={(1+3\cdot 4^2wu^4v^2)v\over (1+4^2wu^4v^2)^2},\ \ \ \ \beta_+(w)={(1-4^2wu^4v^2)u^4v\over (1+4^2wu^4v^2)^2},$$
$$\gamma_-(w)={-4^2wu^3v^2\over 1-4^2wu^4v^2},\ \ \ \ \gamma_+(w)={-4u^3v\over 1-4^2wu^4v^2},$$
$$\eta(w)={1\over 2}\ln(1+4wu^2v).$$
In particular:

\n i) $(x_0^+)^{(k)}(X_1^-)^{(l)}\in\tuz^-\tuz^0\tuz^+$ for all $k,l\in\N$;

\n ii) $\hat h_+(4u)^{{1\over 2}}\in\tuz[[u]]$.
\begin{proof}
We use the notation fixed in \ref{gabg}.

\n It is obvious that $G(0,v)=\exp(X_1^-v)$, so that the condition a) of remark \ref{gesponenziale} is fulfilled, and we need to verify condition b), following lemmas \ref{derivcomm},vi) and \ref{wderivcomm}. 

\n Remark that 
$${d\eta(w)\over du}={4wuv\over 1+4wu^2v}={4wuv(1-4wu^2v)\over 1-4^2w^2u^4v^2}=w\alpha_-(w^2)+\gamma_-(w^2)$$
and 
$$\exp(-4\eta(w)+2\eta(-w))={1-4wu^2v\over(1+4wu^2v)^2},$$
$$\alpha_+(w^2)+w\gamma_+(w^2)={u(1-4wu^2v)\over1-4^2w^2u^4v^2}={u\over 1+4wu^2v},$$
so that 
$${d(\alpha_+(w^2)+w\gamma_+(w^2))\over du}={1+4wu^2v-8wu^2v\over (1+4wu^2v)^2}=\exp(-4\eta(w)+2\eta(-w)).$$
Now let us recall that for all $n,m\in\N$
$${d\over du}{u^n\over (1-au^4)^m}={nu^{n-1}+(4m-n)au^{n+3}\over (1-au^4)^{m+1}},$$
hence, fixing $a=4^2w^2v^2$, we get
$${d\alpha_-(w^2)\over du}={4v(1+3au^4)\over(1-au^4)^2},$$
$${d\beta_-(-w^2)\over du}={-4au^3v(1+3au^4)\over(1-au^4)^3},$$
$${d\gamma_-(w^2)\over du}={-a(3u^2+au^6)\over(1-au^4)^2},$$
$${d\alpha_+(w^2)\over du}={1+3au^4\over(1-au^4)^2},$$
$${d\beta_+(-w^2)\over du}={4vu^3(1+3au^4)\over(1-au^4)^3}.$$

The relations to prove are then equivalent to the following:
$$4v(1+3au^4)=4(1-3au^4)v+6\cdot 4uv\cdot au^3,$$
$$-4au^3v(1+3au^4)=4uv(-w^24^2u^2v^2-3a^2u^6),$$
$$-a(3u^2+au^6)=-3w^2\cdot 4^2u^2v^2-a^2u^6,$$
$$4u^3v(1+3au^4)=(1+3au^4)4u^3v,$$
which are easily verified.

\n Then, since $\alpha_{\pm}$, $\beta_{\pm}$, $\gamma_{\pm}$ have integral coefficients, i) follows from example \ref{dvdpw}, remark \ref{whtilde} and lemma \ref{ometiomecap},iii).

\n ii) follows at once from the above considerations, inverting the exponentials.

\end{proof}
\end{proposition}

\begin{proposition}\label{x0piux1meno}
$$\exp(x_0^+u)\exp(x_1^-v)=$$
$$=\exp({\alpha_-})\exp({\beta_-})\exp({\gamma_-})\exp(\eta)\exp(\gamma_+)\exp(\beta_+)\exp(\alpha_+)$$
where, with the notations of lemma \ref{wderivcomm},

$$\alpha_+(w)={(1+wu^2v^2)u\over 1-6wu^2v^2+w^2u^4v^4},\ \ \ \ \alpha_-(w)={(1+wu^2v^2)v\over 1-6wu^2v^2+w^2u^4v^4},$$
$$\beta_+(w)={(1-4wu^2v^2-w^2u^4v^4)u^3v\over (1+6wu^2v^2+w^2u^4v^4)^2},\ \ \ \ 
\beta_-(w)={(1-4wu^2v^2-w^2u^4v^4)wuv^3\over (1+6wu^2v^2+w^2u^4v^4)^2},$$
$$\gamma_+(w)={(-3+wu^2v^2)u^2v\over 1-6wu^2v^2+w^2u^4v^4},\ \ \ \ 
\gamma_-(w)={(-3+wu^2v^2)wuv^2\over 1-6wu^2v^2+w^2u^4v^4},$$
$$\eta(w)={1\over 2}\ln(1+2wuv-w^2u^2v^2).$$
\begin{proof}
\n We use the notations fixed in \ref{gabg}.

\n \n It is obvious that $G(0,v)=\exp(x_1^-v)$, so that the condition a) of remark \ref{gesponenziale} is fulfilled, and we need to verify condition b), following lemma \ref{wderivcomm}. 

\n First of all remark that $$1-6t^2+t^4=(1+2t-t^2)(1-2t-t^2)$$ and that 
$$1+t^2+(-3+t^2)t=1-3t+t^2+t^3=(1-t)(1-2t-t^2);$$
thus,
replacing $t$ by $wuv$, we get
$$\alpha_+(w^2)+w\gamma_+(w^2)={(1-wuv)u\over 1+2wuv-w^2u^2v^2}$$
and
$$w\alpha_-(w^2)+\gamma_-(w^2)={(1-wuv)wv\over 1+2wuv-w^2u^2v^2}.$$
Hence the relations of lemma \ref{wderivcomm} involving $\eta$ are easily proved:
$${d\eta(w)\over du}={(1-wuv)wv\over 1+2wuv-w^2u^2v^2}=w\alpha_-(w^2)+\gamma_-(w^2)$$
and
$$\exp(-4\eta(w)+2\eta(-w))={1-2wuv-w^2u^2v^2\over(1+2wuv-w^2u^2v^2)^2}$$
while, on the other hand,
$${d\over dt}{t-t^2\over 1+2t-t^2}={1-2t-t^2\over (1+2t-t^2)^2}$$
so that
$${d\over du}(\alpha_+(w^2)+w\gamma_+(w^2))={1-2wuv-w^2u^2v^2\over(1+2wuv-w^2u^2v^2)^2}$$
and
$$\exp(-4\eta(w)+2\eta(-w))={d\over du}(\alpha_+(w^2)+w\gamma_+(w^2)).$$

\n In order to prove the remaining relations remark that for all $n,m\in\N$
$${d\over dt}{t^n\over (1-6t^2+t^4)^m}={nt^{n-1}+6(2m-n)t^{n+1}+(n-4m)t^{n+3}\over (1-6t^2+t^4)^{m+1}},$$
which helps to compute the derivative of $\alpha_{\pm}(w^2)$, $\beta_{\pm}(-w^2)$, $\gamma_-(w^2)$, fixing $t=wuv$ and recalling that ${d\over du}=wv{d\over dt}$:
$${d\alpha_-(w^2)\over du}={wv^2(14t-4t^3-2t^5)\over (1-6t^2+t^4)^2},$$
$${d\beta_-(-w^2)\over du}={w^2v^3(-1-30t^2-12t^4+14t^6-3t^8)\over (1-6t^2+t^4)^3},$$
$${d\gamma_-(w^2)\over du}={w^2v^2(-3-15t^2+3t^4-t^6)\over (1-6t^2+t^4)^2},$$
$${d\alpha_+(w^2)\over du}={1+9t^2-9t^4-t^6\over (1-6t^2+t^4)^2},$$
$${d\beta_+(-w^2)\over du}={w^{-2}v^{-1}(3t^2+26t^4-36t^6+6t^8+t^{10})\over (1-6t^2+t^4)^3}.$$
The relations to prove are then equivalent to the following:
$$14t-4t^3-2t^5=-4(1+4t^2-t^4)t-6(1+t^2)(-3+t^2)t,$$
$$-1-30t^2-12t^4+14t^6-3t^8=(1+t^2)(-(1+t^2)^2-3(-3+t^2)^2t^2), $$
$$-3-15t^2+3t^4-t^6=-3(1+t^2)^2-(-3+t^2)^2t^2,$$
$$3t^2+26t^4-36t^6+6t^8+t^{10}=-(1+9t^2-9t^4-t^6)(-3+t^2)t^2,$$
which are easily verified.

\end{proof}

\end{proposition}

\begin{remark}
Since $(1+2t-t^2)^{-1}\in\Z[[t]]$ proposition \ref{x0piux1meno} implies that $G^{\pm}\in\tuz^{\pm}[[u,v]]$
(see notation \ref{gabg}). Then, in order to prove that $$(x_0^+)^{(k)}(x_1^-)^{(l)}\in\tuz^{-}\tuz^{0}\tuz^{+},$$ we just need to show that $\exp(\eta)\in\tuz^{0}[[u,v]]$. This will imply that $\tuz^{-}\tuz^{0}\tuz^{+}$ is closed under multiplication, hence it is an integral form of $\tu$, obviously containing $\tuz$.

\n In order to prove that $\tuz=\tuz^{-}\tuz^{0}\tuz^{+}$ we need to show in addition that $\tuz^0\subseteq\tuz$.

\n The last part of this paper is devoted to prove that $$\exp\left({1\over 2}\ln(1+2u-u^2).h_0\right)\in\tuz^0[[u]]$$ (see corollary \ref{hcappucciod}) and that $\tuz^0\subseteq\tuz$ (see proposition \ref{samealg}).

\end{remark}

\begin{notation} \label{notedn}

\n In the following $d:\Z_+\to\Q$ denotes the function defined by
$$\sum_{n>0}(-1)^{n-1}{d_n\over n}u^n={1\over 2}\ln(1+2u-u^2)$$
and $\tilde d=\varepsilon d$ (that is $\tilde d_n=\varepsilon_n d_n$ for all $n>0$, where $\varepsilon_n$ has been defined in definition \ref{thuz}).

\n Remark that with this notation we have $\exp(\eta)=\hat h_+^{\{d\}}(uv)$ ($\eta$ as in lemma \ref{wderivcomm} and proposition \ref{x0piux1meno}, $\hat h_+^{\{d\}}(u)$ as in notation \ref{hcappucciof}, where we replace $\hat h^{\{d\}}(u)$ by $\hat h_+^{\{d\}}(u)$ in order to distinguish it from its symmetric $\hat h_-^{\{d\}}(u)=\Omega(\hat h_+^{\{d\}}(u))$).
\end{notation}

\begin{remark}\label{hhdehh}
From $1+2u-u^2=(1+(1+\sqrt 2)u)(1+(1-\sqrt 2)u)$, we get that:

\n i) for all $n\in\Z_+$ $d_n={1\over 2}((1+\sqrt{2})^n+(1-\sqrt{2})^n)$; equivalently $\exists\delta_n\in\Z$ such that
$$\forall n\in\Z_+\ \ \ (1+\sqrt{2})^n=d_n+\delta_n\sqrt{2}.$$

\n ii) $d_n$ is odd for all $n\in\Z_+$; $\delta_n$ is odd if and only if $n$ is odd.

\n iii) $\Z[\hat h_k^{\{d\}}|k>0]\not\subseteq\Z[\hat h_k|k>0]$ (indeed $(\mu*d)(4)=d_4-d_2=17-3=14$, which is not a multiple of 4, see propositions \ref{convoluzioneintera} and \ref{emmepiallaerre}).

\n iv) $\Z[\hat h_k^{\{d\}}|k>0]\subseteq\Z[\tilde h_k|k>0]$ if and only if $\Z[\hat h_k^{\{\tilde d\}}|k>0]\subseteq\Z[\hat h_k|k>0]$.

\end{remark}

\begin{lemma}\label{contidp}
\n Let $p,m,r\in\Z_+$ be such that $p$ is prime and $(m,p)=1$. Then

$${\rm{if}}\ p^r=4\ \ \ p^r=4|d_{4m}+d_{2m},$$
$${\rm{if}}\ p^r\neq 4\ \ \ p^r|d_{p^rm}-d_{p^{r-1}m}.$$

\begin{proof}
\n The claim is obvious for $p^r=2$ since the $d_n$'s are all odd.

\n In general if $n$ is any positive integer it follows from remark \ref{hhdehh} that $$d_{np}+\delta_{np}\sqrt{2}=(d_n+\delta_n\sqrt{2})^p.$$
If $p=2$ this means that
$$d_{2n}=d_n^2+2\delta_n^2,$$
$$\delta_{2n}=2d_n\delta_n,$$
hence
$$2^r||\delta_{2^rm}\ \  {\rm{(recall\ that}}\ \delta_m\ {\rm{is\ odd\ since}}\ m\ {\rm{is\ odd)}}$$
$$d_{2^rm}\equiv
d_{2^{r-1}m}^2\ \ (mod\ 2^{2r-1}),$$
from which it follows that
$$d_{2m}\equiv -1\ \ (mod\ 4),$$
$$d_{2^rm}\equiv 1\ \ (mod\ 2^{r+1})\ \ \ {\rm{if}}\ r>1:$$
indeed, since $d_m$ and $\delta_m$  are odd,
$$d_{2m}\equiv_{(8)}1+2\equiv_{(4)}-1,$$
while if $r\geq 2$ then $2r-1\geq r+1$ and by induction on $r$ we get
$$d_{2^rm}\equiv d^2_{2^{r-1}m}=(\pm 1+2^rk)^2\equiv 1\ \ (mod\ 2^{r+1}).$$
These last relations immediately imply the claim for $p=2$.

\n Now let $p\neq 2$. Then
$$d_{pn}=\sum_{h\geq 0}{p\choose 2h}2^hd_{n}^{p-2h}\delta_n^{2h},$$
$$\delta_{pn}=\sum_{h\geq 0}{p\choose 2h+1}2^hd_{n}^{p-2h-1}\delta_n^{2h+1}.$$
Suppose that $d_n=d+p^{r-1}k$, $\delta_n=\delta+p^{r-1}k'$ with $k=k'=0$ if $r=1$. Then 
$$d_{pn}\equiv\sum_{h\geq 0}{p\choose 2h}2^hd^{p-2h}\delta^{2h}\ \ (mod\ p^r),$$
$$\delta_{pn}\equiv\sum_{h\geq 0}{p\choose 2h+1}2^hd^{p-2h-1}\delta^{2h+1}\ \ (mod\ p^r).$$
The above relations allow us to prove by induction on $r>0$ that if $\zeta_p$ is defined by the properties $\zeta_p\in\{\pm 1\}$, $\zeta_p\equiv_{(p)} 2^{{p-1\over 2}}$ then
$$d_{p^rm}\equiv d_{p^{r-1}m}\ \ (mod\ p^r)\ \ \ {\rm{and}}\ \ \ \delta_{p^rm}\equiv \zeta_p\delta_{p^{r-1}m}\ \ (mod\ p^r):$$
indeed
if $r=1$ 
$$d_{pm}\equiv d_m^p\equiv d_m\ \ (mod\ p),$$
$$\delta_{pm}\equiv 2^{{p-1\over 2}}\delta_m^p\equiv\zeta_p\delta_m\ \ (mod\ p);$$
if $r>1$ then
$$d_{p^rm}\equiv_{(p^r)}\sum_{h\geq 0}{p\choose 2h}2^hd_{p^{r-2}m}^{p-2h}\delta_{p^{r-2}m}^{2h}\equiv_{(p^r)}d_{p^{r-1}m},$$
$$\delta_{p^rm}\equiv_{(p^r)}\zeta_p\sum_{h\geq 0}{p\choose 2h+1}2^hd_{p^{r-2}m}^{p-2h-1}\delta_{p^{r-2}m}^{2h+1}\equiv_{(p^r)}\zeta_p\delta_{p^{r-1}m}.$$
\end{proof}
\end{lemma}
\vskip .3 truecm

\begin{corollary}\label{hcappucciod}
\n $\hat h_n^{\{d\}}\in\Z[\tilde h_k|k>0]$ for all $n>0$.

\n In particular $(x_0^+)^{(k)}(x_1^-)^{(l)}\in\tuz^-\tuz^0\tuz^+$ $\forall k,l\in\N$.

\begin{proof}
\n The claim follows from propositions \ref{convoluzioneintera} and \ref{emmepiallaerre}, remark \ref{hhdehh} and lemma \ref{contidp}, remarking that if $m$ is odd then
$$ d_{4m}+d_{2m}=-(\tilde d_{4m}-\tilde d_{2m})$$
while if $(m,p)=1$ and $p^r\neq 4$ then
$$d_{p^rm}-d_{p^{r-1}m}=\pm(\tilde d_{p^rm}-\tilde d_{p^{r-1}m})
.$$
Thus for all $n>0$ $\hat h_n^{\{\tilde d\}}\in\Z[\hat h_k|k>0]$ and 
$\hat h_n^{\{d\}}\in\Z[\tilde h_k|k>0]$.
\end{proof}
\end{corollary}

\begin{corollary}\label{b22}
$\tuz^+\tuz^-\subseteq\tuz^-\tuz^0\tuz^+$; equivalently $\tuz^-\tuz^0\tuz^+$ is an integral form of $\tu$.
\begin{proof}
The proof is identical to that of proposition \ref{strutmodulo} replacing $\hu$ with $\tu$, having care to remark that in this case, too,
$$(x_r^+)^{(k)}(x_s^-)^{(l)}\in\sum_{m\geq 0}\tu_{\Z,-l+m}^-\tuz^0\tu_{\Z,k-m}^+\ \ \forall r,s\in\Z,\ \forall k,l\in\N:$$
if $r+s$ is even this follows at once comparing proposition \ref{czzq} with the properties of the gradation, while if $r+s$ is odd it is true by proposition \ref{x0piux1meno} and remark \ref{autstab},iv).

\end{proof}
\end{corollary}

\begin{proposition} \label{samealg}
$\tuz^0\subseteq\tuz$ and $\tuz=\tuz^-\tuz^0\tuz^+$.

\begin{proof}
\n Let ${\cal Z}$ be the $\Z$-subalgebra of $\Q[h_r|r>0]$ generated by the coefficients of
$\hat h_+^{\{d\}}(u)$ and of $\hat h_+(4u)^{1/2}$.
Remark that, by propositions \ref{xmenogrande} and \ref{x0piux1meno}, ${\cal Z}\subseteq\tuz$.

\n We have already proved that ${\cal Z}\subseteq\Z[\tilde h_k|k>0]$ (see lemma \ref{ometiomecap},iii) and corollary \ref{hcappucciod}). Let us prove, by induction on $j$, that $\tilde h_j\in{\cal Z}$ for all $j>0$.

\n If $j=1$ the claim depends on the equality $\tilde h_1=h_1=\hat h^{\{d\}}_1$ (since $\varepsilon_1=d_1=1$).

\n Let $j>1$ and suppose that $\tilde h_1,...,\tilde h_{j-1}\in{\cal Z}$.

\n We notice that if $a:\Z_+\to\Z$ is such that $\hat h_j^{\{a\}}\in{\cal{Z}}$ then $a_j\tilde h_j\in{\cal{Z}}$: indeed it is always true that
$$\tilde h_j-{\varepsilon_jh_j\over j}\in\Q[h_1,...,h_{j-1}]$$
and
$$\hat h_j^{\{a\}}-{a_jh_j\over j}\in\Q[h_1,...,h_{j-1}]$$
from which we get that
$$\hat h_j^{\{a\}}-\varepsilon_ja_j\tilde h_j\in\Q[h_1,...,h_{j-1}];$$
but the condition $\hat h_j^{\{a\}}\in{\cal{Z}}\subseteq\Z[\tilde h_k|k>0]$ and the inductive hypothesis 
$\Z[\tilde h_1,...,\tilde h_{j-1}]\subseteq{\cal{Z}}$
imply that 
$$\hat h_j^{\{a\}}-\varepsilon_ja_j\tilde h_j\in\Q[h_1,...,h_{j-1}]
\cap\Z[\tilde h_k|k>0]=\Z[\tilde h_1,...,\tilde h_{j-1}]
\subseteq{\cal{Z}}$$
hence $a_j\tilde h_j\in{\cal{Z}}$.

\n This in particular holds for $a=d$ and for $\hat h^{\{a\}}(u)=\hat h_+(4u)^{{1\over 2}}$, hence
$$d_j\tilde h_j\in{\cal Z}\ \ {\rm{and}}\ \ 2^{2j-1}\tilde h_j\in{\cal Z}.$$
But $(d_j,2^{2j-1})=1$ because $d_j$ is odd, hence $\tilde h_j\in{\cal Z}$.

\n Then $\tuz^{0,+}=\Z[\tilde h_k|k>0]={\cal Z}\subseteq\tuz$ and, applying $\Omega$, $\tuz^{0,-}\subseteq\tuz$. The claim follows recalling corollary \ref{czzp}.
\end{proof}
\end{proposition}

\n We can now collect all the results obtained till now in the main theorem of this work.

\vskip .3 truecm
\begin{theorem}\label{trmA22}
The $\Z$-subalgebra $\tuz$ of $\tu$ generated by $$\{(x_r^+)^{(k)},(x_r^-)^{(k)}|r\in\Z,k\in\N\}$$ is an integral form of $\tu$. 

\n More precisely
$$\tuz\cong
\tuz^{-,1}\otimes\tuz^{-,c}\otimes\tuz^{-,0}\otimes\tuz^{0,-}\otimes\tuz^{0,0}\otimes\tuz^{0,+}\otimes\tuz^{+,1}\otimes\tuz^{+,c}\otimes\tuz^{+,0}$$

and a $\Z$-basis of $\tuz$ is given by the product
$$B^{-,1}B^{-,c}B^{-,0}B^{0,-}B^{0,0}B^{0,+}B^{+,1}B^{+,c}B^{+,0}$$
where $B^{\pm,0}$, $B^{\pm,1}$, $B^{\pm,c}$, $B^{0,\pm}$ and $B^{0,0}$ are the $\Z$-bases respectively of $\tuz^{\pm,0}$, $\tuz^{\pm,1}$, $\tuz^{\pm,c}$, $\tuz^{0,\pm}$ and $\tuz^{0,0}$ given as follows:
$$B^{\pm,0}=\Big\{{(\bf{x}}^{\pm,0})^{({\bf{k}})}=\prod_{r\in\Z}(x_{2r}^{\pm})^{(k_r)}|{\bf{k}}:\Z\to\N\,\, {\rm{is\, finitely\, supported}}\Big\}$$
$$B^{\pm,1}=\Big\{{(\bf{x}}^{\pm,1})^{({\bf{k}})}=\prod_{r\in\Z}(x_{2r+1}^{\pm})^{(k_r)}|{\bf{k}}:\Z\to\N\,\, {\rm{is\, finitely\, supported}}\Big\}$$
$$B^{\pm,c}=\Big\{{(\bf{X}}^{\pm})^{({\bf{k}})}=\prod_{r\in\Z}(X_{2r+1}^{\pm})^{(k_r)}|{\bf{k}}:\Z\to\N\,\, {\rm{is\, finitely\, supported}}\Big\}$$
$$B^{0,\pm}=\Big\{{\tilde{{\bf{h}}}_{\pm}^{\bf{k}}}=\prod_{l\in\N}{\tilde h_{\pm l}^{k_l}}|{\bf{k}}:\N\to\N\,\,{\rm{is\, finitely\, supported}}\Big\}$$
$$B^{0,0}=\Big\{{h_0\choose k}{c\choose\tilde k}|k,\tilde k\in\N\Big\}.$$
\end{theorem}

\vskip .5 truecm

\clearpage

\appendix\label{appendix}

\addcontentsline{toc}{section}{Appendices}

\renewcommand{\thesubsection}{\Alph{subsection}}
\subsection{Straightening formulas of $A^{(2)}_2$}\label{appendA}

\newtheorem{mydefinition}{Definiton}
\numberwithin{mydefinition}{subsection}
\newtheorem{mytheorem}[mydefinition]{Theorem}
\newtheorem{myremark}[mydefinition]{Remark}
\newtheorem{mynotation}[mydefinition]{Notation}
\newtheorem{mylemma}[mydefinition]{Lemma}
\newtheorem{mycorollary}[mydefinition]{Corollary}

\numberwithin{equation}{subsection}

For the sake of completeness we collect here the commutation formulas of $A^{(2)}_2$, inserting also the formulas that we didn't need for the proof of theorem
\ref{trmA22}.

\n Notation \ref{pppm} and remark \ref{epppm} will help writing some of the following straightening relations and to understand the origin of some apparently misterious terms.

\begin{mynotation}\label{pppm}

Given $p(t)\in\Q[[t]]$ let us define $p_+(t), p_-(t)\in\Q[[t^2]]$ and $p_0(t)\in\Q[[t]]$ by $$p(t)=p_+(t)+tp_-(t),\ p_0(t^2)={1\over 2}p_+(t)p_-(t).$$
Remark that the maps $p(t)\mapsto p_+(t)$ and $p(t)\mapsto p_-(t)$ are homomorphisms of $\Q[[t^2]]$-modules while $q(t)\in\Q[[t^2]], \tilde q(t^2)=q(t)\Rightarrow
(qp)_0(t)=\tilde q(t)^2p_0(t)$.
\end{mynotation}
\begin{myremark}\label{epppm}
Given $p(t)\in\Q[[t]]$ we have that
$$\exp(p(uw).x_0^+)=$$
$$=\exp(p_+(uw).x_0^+)\exp(up_0(-u^2w).X_1^+)\exp(up_-(uw).x_1^+)=$$
$$=\exp(up_-(uw).x_1^+)\exp(-up_0(-u^2w).X_1^+)\exp(p_+(uw).x_0^+).$$

\end{myremark}
\vskip .3 truecm
\n We shall now list a complete set of {\bf{straightening formulas}} in $\tuz$.

\vskip .3 truecm
\n I) Zero commutations regarding $\tuz^{0,0}$:
$${c\choose k}\ {\rm{is\ central\ in\ }}\tuz;$$
$${h_0\choose k}\ {\rm{is\ central\ in\ }}\tuz^0:\left[{h_0\choose k},\tilde h_l\right]=0\ \forall k\geq 0,\ l\neq 0.$$

\n II) Relations in $\tuz^{0,+}$ (from which those in $\tuz^{0,-}$ follow as well):
$$\tuz^{0,+}\ {\rm{is\ commutative}}:\ [\tilde h_k,\tilde h_l]=0\ \forall k,l>0;$$ 
$$\tilde \lambda_m(\tilde h_+(-u^m))=\prod_{j=1}^{m}\tilde h_+(-\omega^ju)\ \forall m\in \Z_+$$
${\rm{where}}\ \omega\ {\rm{is\ a\ primitive\ }}m^{{\rm{th}}}\ {\rm{root\ of\ }}1$, that is $$\tilde \lambda_m(\tilde h_k)=(-1)^{(m-1)k}\sum_{(k_1,...,k_m):\atop k_1+...+k_m=mk} \omega^{\sum_{j=1}^m j k_j}\tilde h_{k_1} \dots \tilde h_{k_m} ;$$ 
if $m$ is odd
$$\lambda_m(\tilde h_k)=\tilde \lambda_m(\tilde h_k)\ \forall k\geq 0;$$
if $m$ is even
$$\lambda_{m}(\hat h_+(u))=\tilde\lambda_{m}(\tilde h_+((-1)^{m\over 2}u)^{-1});$$

$$\hat h_+(u)=\tilde h_+(u)\tilde\lambda_4(\tilde h_+(-u^4)^{-{1\over 2}})=\tilde h_+(u)^{{1\over 2}}\tilde h_+(-u)^{-{1\over 2}}\tilde h_+(iu)^{-{1\over 2}}\tilde h_+(-iu)^{-{1\over 2}},$$

$$\hat h_+^{\{d\}}(u)=\hat h_+((1+\sqrt{2})u)^{{1\over 2}}\hat h_+((1-\sqrt{2})u)^{{1\over 2}}=\prod_{m>0}\tilde\lambda_m(\tilde h_+(u^m))^{k_m}$$
where the $k_m$'s are integers defined by the identity
$$1+2u-u^2=(1-2u-u^2)(1+6u^2+u^4)\prod_{m>0}(1+u^m)^{4k_m}.$$

\n The corresponding relations in $\tuz^{0,-}$ are obtained applying $\Omega$, that is just replacing $\tilde h_k$, $\tilde h_+(u)$ and $\hat h_+(u)$ with $\tilde h_{-k}$, $\tilde h_-(u)$ and $\hat h_-(u)$.

\n III) Other straightening relations in $\tuz^0$:
$$\tilde h_+(u)\tilde h_-(v)=\tilde h_-(v)(1-uv)^{-4c}(1+uv)^{2c}\tilde h_+(u).$$

\n IV) Commuting elements and straightening relations in $\tuz^+$ (and in $\tuz^-$):
$$(X_{2r+1}^+)^{(k)}\ {\rm{is\ central\ in\ }}\tuz^+:$$ $$[(X_{2r+1}^+)^{(k)},(x_s^+)^{(l)}]=0=[(X_{2r+1}^+)^{(k)},(X_{2s+1}^+)^{(l)}]\ \forall r,s\in\Z,\ k,l\in\N;$$
$${\rm{if}}\ r+s\ {\rm{is\ even\ }}[(x_r^+)^{(k)},(x_s^+)^{(l)}]=0\ \forall k,l\in\N;$$
$${\rm{if}}\ r+s\ {\rm{is\ odd\ }}
\exp(x_r^+ u)\exp(x_{s}^+ v)=\exp(x_s^+ v)\exp((-1)^sX_{r+s}^+ uv)\exp(x_{r}^+ u).$$

\n All the relations in $\tuz^-$ are obtained from those in $\tuz^+$ applying the antiautomorphism $\Omega$; in particular
if $r+s$ is odd
$$\exp(x_r^- u)\exp(x_{s}^- v)=\exp(x_s^- v)\exp((-1)^rX_{r+s}^- uv)\exp(x_{r}^- u).$$

\n V) Straightening relations for $\tuz^+\tuz^{0,0}$ (and for $\tuz^{0,0} \tuz^-$): $\forall r\in\Z,\ k,l\in\N$

$$(x_r^+)^{(k)}{h_0\choose l}={h_0-2k\choose l}(x_r^+)^{(k)},$$
$$(X_{2r+1}^+)^{(k)}{h_0\choose l}={h_0-4k\choose l}(X_{2r+1}^+)^{(k)},$$
and
$${h_0\choose l}(x_r^-)^{(k)}=(x_r^-)^{(k)}{h_0-2k\choose l},$$
$${h_0\choose l}(X_{2r+1}^-)^{(k)}=(X_{2r+1}^-)^{(k)}{h_0-4k\choose l}.$$

\n VI) Straightening relations for $\tuz^+\tuz^{0,+}$ (and for $\tuz^+\tuz^{0,-}$, $\tuz^{0,\pm} \tuz^-$): 

$$(X_{2r+1}^+)^{(k)}\tilde h_+(u)=\tilde h_+(u)\left((1-u^2T^{-1})^2X_{2r+1}^+\right)^{(k)}$$
and $$(x_r^+)^{(k)}\tilde h_+(u)=\tilde h_+(u)\left({(1-uT^{-1})^6(1+u^2T^{-2})\over(1-u^2T^{-2})^3}x_r^+\right)^{(k)};$$
the expression for $\left({(1-uT^{-1})^6(1+u^2T^{-2})\over(1-u^2T^{-2})^3}x_r^+\right)^{(k)}$ can be straightened more explicitly: setting $p(t)=(1-t)^6$ we have $$
p_+(t)=1+15t^2+15t^4+t^8,$$ $$
p_-(t)=-6-20t^2-6t^4,$$ $$
p_0(t)=-(1+15t+15t^2+t^4)(3+10t+3t^2),$$
so that 
$$\exp(x_r^+v)\tilde h_+(u)=\tilde h_+(u)\exp\left({(1-uT^{-1})^6(1+u^2T^{-2})\over(1-u^2T^{-2})^3}x_r^+v\right)=$$
$$=\tilde h_+(u)\exp\left({p_-(uT^{-1})(1+u^2T^{-2})\over(1-u^2T^{-2})^3}x_{r+1}^+uv\right)\cdot$$
$$\cdot\exp\left({(-1)^{r-1}p_0(-u^2T^{-1})(1-u^2T^{-1})^2\over(1+u^2T^{-1})^6}X_{2r+1}^+uv^2\right)\cdot$$
$$\cdot\exp\left({p_+(uT^{-1})(1+u^2T^{-2})\over(1-u^2T^{-2})^3}x_r^+v\right).$$
Applying the homomorphism $\lambda_{-1}$ (that is $x_s^+\mapsto x_{-s}^+$, $X_s^+\mapsto X_{-s}^+$, $\tilde h_+\mapsto\tilde h_-$, $T^{-1}\mapsto T$) one immediately gets the expression for $(X_{2r+1}^+)^{(k)}\tilde h_-(u)$ and for $\exp(x_{r}^+v)\tilde h_-(u)$.

\n Applying the antiautomorphism $\Omega$ ($x_s^+\mapsto x_{-s}^-$, $X_s^+\mapsto X_{-s}^-$, $\tilde h_+\leftrightarrow\tilde h_-$) one gets analogously the expression for
$\tilde h_{\pm}(u)(X_{2r+1}^-)^{(k)}$ and for $\tilde h_{\pm}(u)\exp(x_{r}^-v)$.

\n VII) Straightening relations for $\tuz^+\tuz^-$:

\n VII,a) $\goths$-like relations: $\forall r\in\Z$
$${\rm{exp}}(x_r^+u){\rm{exp}}(x_{-r}^-v)={\rm{exp}}\Big({x_{-r}^-v\over 1+uv}\Big)(1+uv)^{h_0+rc}{\rm{exp}}\Big({x_r^+u\over 1+uv}\Big),$$
$${\rm{exp}}(X_{2r+1}^+u){\rm{exp}}(X_{-2r-1}^-v)={\rm{exp}}\Big({X_{-2r-1}^-v\over 1+4^2uv}\Big)(1+4^2uv)^{{h_0\over 2}+{(2r+1)c\over 4}}{\rm{exp}}\Big({X_{2r+1}^+u\over 1+4^2uv}\Big).$$

\n VII,b) $\hat\goths$-like relations: 

\n if $r+s\neq 0$ is even
$$\exp(x_{r}^+u)\exp(x_s^-v)=$$
$$=\exp\left({1\over 1+uvT^{r+s}}x_s^-v\right)
\lambda_{r+s}(\hat h_+(uv))
\exp\left({1\over 1+uvT^{-r-s}}x_r^+v\right),$$
while $\forall r+s\neq 0$
$$\exp(X_{2r+1}^+u)\exp(X_{2s-1}^-v)=$$
$$=\exp\left({1\over 1+4T^{s+r}uv}X_{2s-1}^-v\right)\cdot$$
$$\cdot\lambda_{2(r+s)}(\hat h_+(4^2uv)^{{1\over 2}})\cdot$$
$$\cdot\exp\left({1\over1+4uvT^{-s-r}}X_{2r+1}^+u\right).$$

\n VII,c) Straightening relations for $\tuz^{+,0}\tuz^{-,c}$ (and $\tuz^{+,1}\tuz^{-,c}$, $\tuz^{+,c}\tuz^{-,{0\atop 1}}$):
$$\exp(x_0^+u)\exp(X_{1}^-v)=$$
$$=\exp\left({4\over 1-4^2w^2u^4v^2}x_1^-uv\right)\exp\left({-4^2w^2\over 1-4^2w^2u^4v^2}x_0^-u^3v^2\right)\cdot$$
$$\cdot\exp\left({1+3\cdot 4^2wu^4v^2\over (1+4^2wu^4v^2)^2}X_1^-v\right)
\hat h_+(4u^2v)^{{1\over 2}}\exp\left({1-4^2w^{-1}u^4v^2\over (1+4^2w^{-1}u^4v^2)^2}X_1^+u^4v\right)\cdot$$
$$\cdot\exp\left({-4\over 1-4^2w^{-2}u^4v^2}x_1^+u^3v\right)\exp\left({1\over 1-4^2w^{-2}u^4v^2}x_0^+u\right),$$
which can be written in a more compact way observing that 
$${1\over 1-4^2t^2}=\left({1\over 1+4t}\right)_+,\ {-4\over 1-4^2t^2}=\left({1\over 1+4t}\right)_-,\ \left({1\over 1+4t}\right)_0={-2\over (1-4^2t)^2}:$$

$$\exp(x_0^+u)\exp(X_{1}^-v)=$$
$$=\exp\left({4\over 1+4wu^2v}x_1^-uv\right)\exp\left({1\over 1+4^2wu^4v^2}X_1^-v\right)\cdot$$
$$\cdot\hat h_+(4u^2v)^{{1\over 2}}\exp\left({1\over 1+4w^{-1}u^2v}x_0^+u\right)\exp\left(-{1\over 1+4^2w^{-1}u^4v^2}X_1^+u^4v\right);$$
that is more symmetric but less explicit in terms of the given basis of $\tuz$.

\n Applying the homomorphism $T^{-r}\lambda_{2r+2s+1}$ (that is $x_l^{\pm}\mapsto x_{l(2r+2s+1)\pm r}^{\pm}$, $X_1^{\pm}\mapsto (-1)^rX_{2r+2s+1\pm 2r}^{\pm}$, $\hat h_k\mapsto \lambda_{2r+2s+1}(\hat h_k)$, $w\big|_{L^{\pm}}\mapsto T^{\mp(2r+2s+1)}$) one deduces the expression for $\exp(x_r^+u)\exp(X_{2s+1}^-v)$.

\n Applying $\Omega$ one analogously gets the expression for 
$\exp(X_{2r+1}^+u)\exp(x_s^-v)$.

\n VII,d) The remaining relations:

$$\exp(x_0^+u)\exp(x_1^-v)=$$
$$=\exp \left ({{1+w^2u^2v^2}\over{1-6w^2u^2v^2+w^4u^4v^4}}x_1^-v\right) \exp \left({{-3+w^2u^2v^2}\over{1-6w^2u^2v^2+w^4u^4v^4}}x_2^-uv^2\right) \cdot$$
$$\cdot \exp \left (-{{1-4wu^2v^2-w^2u^4v^4}\over{(1+6wu^2v^2+w^2u^4v^4)^2}}X_3^-uv^3\right ) \hat h_+^{\{d\}}(uv)\cdot$$
$$\cdot\exp \left({{1-4wu^2v^2-w^2u^4v^4}\over{(1+6wu^2v^2+w^2u^4v^4)^2}}X_1^+u^3v\right)\cdot$$
$$\cdot \exp \left({{-3+w^2u^2v^2}\over{1-6w^2u^2v^2+w^4u^4v^4}}x_1^+u^2v\right)\exp \left({{1+ w^2u^2v^2}\over{(1-6w^2u^2v^2+w^4u^4v^4)^2}}x_0^+u\right)$$
or, as well,
$$\exp(x_0^+u)\exp(x_1^-v)=$$
$$=\exp\left({1-wuv\over 1+2wuv -w^{2}u^2v^2} x_{1}^{-}v\right) \exp\left({1\over 2(1+6wu^2v^2 +w^2u^4v^4)} X_{3}^{-}uv^3\right) \cdot$$
$$\cdot\hat h_+^{\{d\}}(uv)\cdot $$ 
$$\cdot   \exp\left({1-wuv\over 1+2wuv -w^2u^2v^2}x_0^{+}u\right)\left({-1\over 2(1+6wu^2v^2 +w^2u^4v^4)} X_{1}^{+}u^3v\right).$$

\n The general straightening formula for $\exp(x_r^+u)\exp(x_s^-v)$ when $r+s$ is odd is obtained from the case $r=0$, $s=1$ applying $T^{-r}\lambda_{r+s}$, remarking that $w\big|_{L^{\pm}}\mapsto T^{\mp(r+s)}$.

\subsection{Garland's description of $\uz^{im,+}$}\label{appendB}
In this appendix we focus on the imaginary positive part $\uz^{im,+}$ of $\uz=\uz(\gothg)$ (see section \ref{intr}) when $\gothg$ is an affine Kac-Moody algebra of rank 1 (that is $\gothg=\hat\goths$ or $\gothg=\gothd$: these cases are enough to understand also the cases of higher rank): we aim at a better understanding of Garland's (and Mitzman's) basis of $\uz^{im,+}$
and of its connection with the basis consisting of the monomials in the $\hat h_k$'s, basis which arises naturally from the description of $\uz^{im,+}$ as $\Z^{(sym)}[h_r|r>0]=\Z[\hat h_k|k>0]$.

\n First of all let us fix some notations and recall Garland's description of $\uz^{im,+}$.

\begin{definition}\label{bun}
With the notations of example \ref{rvsf} and proposition \ref{tmom} let us define the following elements and subsets in $\Q[h_r|r>0]$:

i) $b_{{\bf{k}}}=\prod_{m>0}\lambda_m(\hat h_{k_m})$ where ${\bf{k}}:\Z_+\to\N$ is finitely supported;
\begin{equation} \label{GarlandBase}{\rm{ii)}}\,\,
B_{\lambda}=\left\{b_{{\bf{k}}}
|{\bf{k}}:\Z_+\to\N\,\,{\rm{is\, finitely\, supported}}\right\};\,\,\,\,\,\,\,\,\,\,\,\,\,\,\,\,\,\,\,\,\,\,\,\,\,\,\,\,\,\,\,\,\,\,\,\,\,\,\,\,\,\,\,\,\,\,\,\,\,\,\,\,\,\,\,\,\,\,\,\,\,
\end{equation}

iii) $\Z_{\lambda}[h_r|r>0]=\sum_{{\bf{k}}}\Z b_{{\bf{k}}}$ is the $\Z$-submodule of $\Q[h_r|r>0]$ generated by $B_{\lambda}$.

\end{definition}

\n Then, with our notation, Garland's description of $\uz^{im,+}$ can be stated as follows:

\begin{theorem}
$\uz^{im,+}$ is a free $\Z$-module with basis $B_{\lambda}$. 

\n Equivalently:

\n i) $\uz^{im,+}=\Z_{\lambda}[h_r|r>0]$;

\n ii) $B_{\lambda}$ is linearly independent.

\end{theorem}

\begin{remark}
Once proved that $\uz^{im,+}$ is the $\Z$-subalgebra of $\u$ generated by $\{\lambda_m(\hat h_k)|m>0,k\geq 0\}$ (hence by $B_{\lambda}$ or equivalently by $\Z_{\lambda}[h_r|r>0]$), proceeding in two different directions leads to the two descriptions of $\uz^{im,+}$ that we want to compare:

\n $\tic$) $\Z_{\lambda}[h_r|r>0]$ is a $\Z$-subalgebra of $\Q[h_r|r>0]$ (that is $\Z_{\lambda}[h_r|r>0]$ is closed under multiplication): this implies that $$\uz^{im,+}=\Z_{\lambda}[h_r|r>0];$$
it also implies that $\Z[\hat h_k|k>0]\subseteq\Z_{\lambda}[h_r|r>0]$;

\n $\tac$) $\Z[\hat h_k|k>0]$ is $\lambda_m$-stable for all $m>0$ (see proposition \ref{tmom}): this implies that $$\uz^{im,+}=\Z[\hat h_k|k>0];$$
it also implies that $\Z_{\lambda}[h_r|r>0]\subseteq\Z[\hat h_k|k>0]$.

Hence $\tic$) and $\tac$) imply that $\uz^{im,+}=\Z_{\lambda}[h_r|r>0]=\Z[\hat h_k|k>0]$.

\n $\tic$) has been proved in \cite{HG} by induction on a suitably defined degree. The first step of the induction is
the second assertion of \cite{HG}-lemma 5.11(b), proved in \cite{HG}-section 9: for all $k,l\in\N$ $\hat h_k\hat h_l-{k+l\choose k}\hat h_{k+l}$ is a linear combination {\it {with integral coefficients}} of elements of $B_{\lambda}$ {\it{of degree lower}} than the degree of $\hat h_{k+l}$.\\ In the proof the author uses that $B_{\lambda}$ is a $\Q$-basis of $\Q[h_r|r>0]$ and concentrates on the integrality of the coefficients: he studies the action of $\frak h$ on $\hat{{\frak sl_3}}^{\otimes N}$ where $\frak h$ is the commutative Lie-algebra  with basis $\{h_r|r>0\}$ and $N\in\N$ is large enough ($N$ is the maximum among the degrees of the elements of $B_{\lambda}$ appearing in $\hat h_k\hat h_l$ with non-integral coefficient, assuming that such an element exists): $\frak h$ is a subalgebra of $\hat\goths$ and there is an embedding of $\hat\goths$ in $\hat{{\frak sl_3}}$ for every vertex of the Dynkin diagram of ${\frak sl}_3$, so that fixing a vertex of the Dynkin diagram of ${\frak sl}_3$ induces an embedding ${\frak h}\subseteq\hat\goths\hookrightarrow\hat{{\frak sl_3}}$, hence an action of ${\frak h}$ on $\hat{{\frak sl_3}}$. 
But the integral form of $\hat{{\frak sl_3}}$ defined as the $\Z$-span of a Chevalley basis is $\uz(\hat{\frak sl_3})$-stable;
since the stability under $\uz(\hat{\frak sl_3})$ is preserved by tensor products (\cite{HG}-section 6), the author can finally deduce the desired integrality property of $\hat h_k\hat h_l$ from the study of the $\frak h$-action on $\hat{{\frak sl_3}}^{\otimes N}$.

\vskip .3 truecm
\n Garland's argument has been sometimes misunderstood: it is the case for instance of \cite{JM} where the authors affirm (in lemma 1.5) that \cite{HG}-lemma 5.11(b) implies that $\uz^{im,+}=\Z[\hat h_k|k>0]$, while, as discussed above, it just implies the inclusion $\Z[\hat h_k|k>0]\subseteq\uz^{im,+}=\Z_{\lambda}[h_r|r>0]$.

\n On the other hand Garland's argument strongly involves many results of the (integral) representation theory of the Kac-Moody algebras, while $\tic$) is a property of the algebra $\Q[h_r|r>0]$ and of its integral forms that can be stated in a way completely independent of the Kac-Moody algebra setting:
$$\Z^{(sym)}[h_r|r>0]\subseteq\Z_{\lambda}[h_r|r>0].$$
\n The above considerations motivate the present appendix, whose aim is to propose a self-contained proof of $\tic$), independent of the Kac-Moody algebra context: on one hand we think that a direct proof can help evidentiating the essential structure of the integral form of $\Q[h_r|r>0]$ arising from our study; on the other hand the idea of isolating the single pieces  
and glueing them together after studying them separately is much in the spirit of this work, so that it is natural for us to explain also Garland's basis of $\uz^{im,+}$ through this approach; and finally we hope that presenting a different proof can also help to clarify the steps which appear more difficult in Garland's proof.
\vskip .3 truecm
\n In the following we go back to the description of $\Z[\hat h_k|k>0]$ as the algebra of the symmetric functions and we show that $B_{\lambda}$ is a basis of $\Z[\hat h_k|k>0]$ by comparing it with a well known $\Z$-basis of this algebra.

\end{remark}

\begin{remark}\label{rcdcp}
Recall that $\Z[\hat h_k|k>0]$ is the algebra  of the symmetric functions and that $\forall n\in\N$ the projection $\pi_n:\Z[\hat h_k|k>0]\to\Z[x_1,...,x_n]^{{\cal{S}}_n}$ induces an isomorphism $\Z[\hat h_1,...,\hat h_n]\cong\Z[x_1,...,x_n]^{{\cal{S}}_n}$ through which $\hat h_k$ corresponds to the $k^{{\rm{th}}}$ elementary symmetric polynomial $e_k^{[n]}$, while $\pi_n(\hat h_k)=0$ if $k>n$ and $h_r$ corresponds to the sum of the $r^{{\rm{th}}}$-powers $\sum_{i=1}^nx_i^r$ $\forall r>0$ (see example \ref{rvsf}).

\n Then it is well known and obvious that:

\n i) $\forall{\bf{k}}:\Z_+\to\N$ finitely supported $\exists!(\sigma x)_{{\bf{k}}}\in\Z[\hat h_k|k>0]$ such that 
$$\pi_n((\sigma x)_{\bf{k}})=\sum_{{a_1,...,a_n\atop\#\{i|a_i=m\}=k_m\ \forall m>0}}\prod_{i=1}^nx_i^{a_i}\in\Z[x_1,...,x_n]^{{\cal{S}}_n}\ \ \forall n\in\N;$$

\n ii) $\{(\sigma x)_{\bf{k}}|{\bf{k}}:\Z_+\to\N$ finitely supported$\}$ is a $\Z$-basis of $\Z[\hat h_k|k>0]$.

\n (It is the basis that  in \cite{IM} is denoted by $\{m_{\lambda}|\lambda=(\lambda_1\geq\lambda_2\geq...\geq 0)\}$:  $m_{\lambda}=(\sigma x)_{{\bf{k}}}$ where $\forall m>0$ $k_m=\#\{i|\lambda_i=m\}$).

\end{remark}

\begin{notation}
As in remark \ref{rcdcp}, for all ${\bf{k}}:\Z_+\to\N$ finitely supported let us denote by $(\sigma x)_{\bf{k}}$ the limit of the elements $$\sum_{{a_1,...,a_n\atop\#\{i|a_i=m\}=k_m\ \forall m>0}}\prod_{i=1}^nx_i^{a_i}\ \ (n\in\N).$$
By abuse of notation, when $n\geq\sum_{m>0}k_m$ we shall write 
$$(\sigma x)_{\bf{k}}=\sum_{{a_1,...,a_n\atop\#\{i|a_i=m\}=k_m\ \forall m>0}}\prod_{i=1}^nx_i^{a_i},$$
which is justified because, under the hypothesis that $n\geq\sum_{m>0}k_m$, ${\bf{k}}$ is determined by the set $\{(a_1,...,a_n)|\#\{i=1,...,n|a_i=m\}=k_m$ $\forall m>0\}$.

\end{notation}
\begin{definition} \label{basi}
$\forall n\in\N$ define $B_{\lambda}^{(n)}$, $B_x^{(n)}$, $\Z_{\lambda}^{(n)}$, $\Z_x^{(n)}\subseteq\Q[h_r|r>0]=\Q[\hat h_k|k>0]$ as follows:
$$B_{\lambda}^{(n)}=\left\{b_{\bf{k}}=\prod_{m>0}\lambda_m(\hat h_{k_m})\in B_{\lambda}|\sum_{m>0}k_m\leq n\right\},$$
$$B_x^{(n)}=\left\{(\sigma x)_{\bf{k}}|\sum_{m>0}k_m\leq n\right\},$$
$\Z_{\lambda}^{(n)}$ is the $\Z$-module generated by $B_{\lambda}^{(n)}$,
$\Z_x^{(n)}$ is the $\Z$-module generated by $B_x^{(n)}$.

\end{definition}

\begin{remark}\label{vdbx}
By the very definition of $B_x^{(n)}$ we have that:

\n i) $B_x^{(n)}$ is a basis of $\Z_x^{(n)}\subseteq\Z[\hat h_k|k>0]=\sum_{n'\in\N}\Z_x^{(n')}$, see remark \ref{rcdcp}, ii);

\n ii) $h\in\Z_x^{(n)}$ means that for all $N\geq n$ each monomial in the $x_i$'s appearing in $\pi_N(h)$ with nonzero coefficient involves no more that $n$ indeterminates $x_i$; hence  in particular $$h\in\Z_x^{(n)}, h'\in\Z_x^{(n')}\Rightarrow hh'\in\Z_x^{(n+n')}.$$
\end{remark}
\begin{lemma} \label{sxkp}
Let $n,n',n''\in\N$ and ${\bf{k}}',{\bf{k}}'':\Z_+\to\N$ be such that $n'+n''=n$, $\sum_{m>0}k'_m=n'$, $\sum_{m>0}k''_m=n''$. Then:

\n i) $(\sigma x)_{{\bf{k}}'}\cdot(\sigma x)_{{\bf{k}}''}\in\Z(\sigma x)_{{\bf{k}}'+{\bf{k}}''}\oplus\Z_x^{(n-1)}$;

\n ii) if $k'_m
k''_m=0$ $\forall m>0$ then $(\sigma x)_{{\bf{k}}'}
(\sigma x)_{{\bf{k}}''}-(\sigma x)_{{\bf{k}}'+{\bf{k}}''}\in\Z_x^{(n-1)}$.
\begin{proof}
That $(\sigma x)_{{\bf{k}}'}\cdot(\sigma x)_{{\bf{k}}''}$ lies in $\Z_x^{(n)}$ follows from remark \ref{vdbx},ii), so we just need to: 

\n i) prove that if $\prod_{i=1}^nx_i^{a_i}$ with $a_i\neq 0$ $\forall i=1,...,n$ is the product of two monomials $M'$ and $M''$ appearing with nonzero coefficient respectively in $(\sigma x)_{{\bf{k}}'}$ and in $(\sigma x)_{{\bf{k}}''}$ then $\#\{i|a_i=m\}=k'_m+k''_m$ for all $m>0$;

\n ii) compute the coefficient of $(\sigma x)_{{\bf{k}}'+{\bf{k}}''}$ in the expression of $(\sigma x)_{{\bf{k}}'}\cdot(\sigma x)_{{\bf{k}}''}$ as a linear combination of the $(\sigma x)_{{\bf{k}}}$'s when $\forall m>0$ $k'_m$ and $k''_m$ are not simultaneously non zero, and find that it is 1.

\n i) is obvious because the condition $a_i\neq 0$ $\forall i=1,...,n$ implies that the indeterminates involved in $M'$ and those involved in $M''$ are disjoint sets.

\n For ii) it is enough to show that, under the further condition on $k'_m$ and $k''_m$, the monomial $\prod_{i=1}^nx_i^{a_i}$ chosen in i) uniquely determines $M'$ and $M''$ such that $\prod_{i=1}^nx_i^{a_i}=M'M''$: indeed $$M'=\prod_{i:k'_{a_i}\neq 0}x_i^{a_i}\ \ {\rm{and}}\ M''=\prod_{i:k''_{a_i}\neq 0}x_i^{a_i}.$$

\end{proof}
\end{lemma}

\begin{lemma}\label{sldftv}
Let ${\bf{k}}:\Z_+\to\N$, $n\in\N$ be such that $\sum_{m>0}k_m=n$. Then:

\n i) if $\exists m>0$ such that $k_{m'}= 0$ for all $m'\neq m$ (equivalently $k_m=n$)
we have
$$(\sigma x)_{{\bf{k}}}=\lambda_m(\hat h_{n})=b_{{\bf{k}}}\in\Z_x^{(n)}\cap\Z_{\lambda}^{(n)};$$
\n ii) in general $b_{{\bf{k}}}-(\sigma x)_{{\bf{k}}}\in\Z_x^{(n-1)}$.
\begin{proof}
i) $\forall N\geq n$ we have $$(\sigma x)_{{\bf{k}}}=\sum_{1\leq i_1<...<i_n\leq N}x_{i_1}^m\cdot...\cdot x_{i_n}^m=\lambda_m\left(\sum_{1\leq i_1<...<i_n\leq N}x_{i_1}\cdot...\cdot x_{i_n}\right)=\lambda_m(e_n^{[N]})$$ so that $(\sigma x)_{{\bf{k}}}=\lambda_m(\hat h_{n})$.

\n ii) $b_{{\bf{k}}}=\prod_{m>0}\lambda_m(\hat h_{k_m})=\prod_{m>0}(\sigma x)_{{\bf{k}}^{[m]}}$ where $k^{[m]}_{m'}=\delta_{m,m'}k_m$ $\forall m,m'>0$; thanks to lemma \ref{sxkp},ii) we have that
$\prod_{m>0}(\sigma x)_{{\bf{k}}^{[m]}}-(\sigma x)_{\sum_m{\bf{k}}^{[m]}}\in\Z_x^{(n-1)}$; but $\sum_{m>0}{\bf{k}}^{[m]}={\bf{k}}$ and the claim follows.

\end{proof}
\end{lemma}
\begin{theorem}\label{gmrvs}
$B_{\lambda}$ is a $\Z$-basis of $\Z[\hat h_k|k>0]$ (thus $\Z[\hat h_k|k>0]=\Z_{\lambda}[h_r|r>0]$).
\begin{proof}
We prove by induction on $n$ that $B_{\lambda}^{(n)}$ is a $\Z$-basis of $\Z_x^{(n)}=\Z_{\lambda}^{(n)}$ $\forall n\in\N$, the case $n=0$ being obvious.

\n Let $n>0$: by the inductive hypothesis $B_{\lambda}^{(n-1)}$ and $B_x^{(n-1)}$ are both $\Z$-bases of $\Z_x^{(n-1)}=\Z_{\lambda}^{(n-1)}$; by definition 
$B_x^{(n)}\setminus B_x^{(n-1)}$ represents a $\Z$-basis of $\Z_x^{(n)}/\Z_x^{(n-1)}$ while $B_{\lambda}^{(n)}\setminus B_{\lambda}^{(n-1)}$ represents a set of generators of the $\Z$-module $\Z_{\lambda}^{(n)}/\Z_{\lambda}^{(n-1)}$.

\n Now lemma \ref{sldftv},ii) implies that if $\sum_{m>0}k_m=n$ then $b_{{\bf{k}}}$ and $(\sigma x)_{{\bf{k}}}$ represent the same element in $\Q[\hat h_k|k>0]/\Z_x^{(n-1)}=\Q[\hat h_k|k>0]/\Z_{\lambda}^{(n-1)}$.

\n Hence $B_{\lambda}^{(n)}\setminus B_{\lambda}^{(n-1)}$ represents a $\Z$-basis of $\Z_x^{(n)}/\Z_x^{(n-1)}=\Z_x^{(n)}/\Z_{\lambda}^{(n-1)}$, that is $B_{\lambda}^{(n)}$ is a $\Z$-basis of $\Z_x^{(n)}$; but $B_{\lambda}^{(n)}$ generates $\Z_{\lambda}^{(n)}$ and the claim follows.

\end{proof}
\end{theorem}

\subsection{Comparison with the Mitzman integral form}\label{appendC}
In the present appendix we compare the integral form $\tuz=$ $\duz(\gothd)$ of $\tu$ described in section \ref{ifa22} with the integral form $\uz(\gothd)$ of the same algebra $\tu$ introduced and studied by Mitzman in \cite{DM}, that we denote here by $\um$ and that is easily defined as the $\Z$-subalgebra of $\tu$ generated by the divided powers of the Kac-Moody generators $e_i,f_i$ ($i=0,1$): see also remark \ref{rfink}.

\n More precisely:

\begin{mydefinition}\label{ka22}

\n $\tu$ is the enveloping algebra of the Kac-Moody algebra whose generalized Cartan matrix is $A_2^{(2)}=(a_{i,j})_{i,j\in\{0,1\}}=\left(\begin{matrix}2&-1\\-4&2\end{matrix}\right)$ (see \cite{VK}): it has generators $\{e_i,f_i,h_i|i=0,1\}$ and relations
$$[h_i,h_j]=0,\ \ [h_i,e_j]=a_{i,j}e_j,\ \ [h_i,f_j]=-a_{i,j}f_j,\ \ [e_i,f_j]=\delta_{i,j}h_i\ \ (i,j\in \{0,1\})$$
$$({\rm{ad}}e_i)^{1-a_{i,j}}(e_j)=0=({\rm{ad}}f_i)^{1-a_{i,j}}(f_j)\ \ (i\neq j\in \{0,1\}).$$

\end{mydefinition}

\begin{mydefinition} \label{mif}
The Mitzman integral form $\um$ of $\tu$ is the $\Z$-subalgebra of $\tu$ generated by $\{e_i^{(k)},f_i^{(k)}|i=0,1,\ k\in\N\}$.

\end{mydefinition}

\begin{myremark}
The Kac-Moody presentation of $\tu$ (definition \ref{ka22}) and its presentation given in definition \ref{a22} are identified through the following isomorphism:
$$e_1\mapsto x_0^+,\ \ f_1\mapsto x_0^-,\ \ h_1\mapsto h_0,\ \ e_0\mapsto{1\over 4}X_1^-,\ \ f_0\mapsto{1\over 4}X_{-1}^+,\ \ h_0\mapsto {1\over 4}c-{1\over 2} h_0.$$

\end{myremark}

\begin{mynotation} \label{mitNota}
In order to avoid in the following any confusion and heavy notations, we set:
$$y_{2r+1}^{\pm}={1\over 4}X_{2r+1}^{\pm},\ \ \k_r={1\over 2}h_r,\ \ \tilde c={1\over 4}c$$
where the $X_{2r+1}^{\pm}$'s, the $h_r$'s and $c$ are those introduced in definition \ref{a22}
(thus $e_0=y_1^-$, $f_0=y_{-1}^+$, while the Kac-Moody $h_0$ and $h_1$ appearing in definition \ref{ka22} are respectively $\tilde c-\k_0$ and $2\k_0$; moreover
$\um$ is the $\Z$-subalgebra of $\tu$ generated by $\{(x_0^{\pm})^{(k)},(y_{\mp 1}^{\pm})^{(k)}|k\in\N\}$).
\end{mynotation}

\begin{myremark}\label{stau}
$\um$ is $\Omega$-stable, $\exp(\pm{\rm{ad}}e_i)$-stable and $\exp(\pm{\rm{ad}}f_i)$-stable.
In particular $\um$ is stable under the action of 
$$\tau_0=\exp({\rm{ad}}e_0)\exp(-{\rm{ad}}f_0)\exp({\rm{ad}}e_0)=\exp({\rm{ad}}y_1^-)\exp(-{\rm{ad}}y_{-1}^+)\exp({\rm{ad}}y_1^-)$$ and
$$\tau_1=\exp({\rm{ad}}e_1)\exp(-{\rm{ad}}f_1)\exp({\rm{ad}}e_1)=\exp({\rm{ad}}x_0^+)\exp(-{\rm{ad}}x_0^-)\exp({\rm{ad}}x_0^+)$$
(cfr. \cite{JH}).

\begin{proof}
The claim for $\Omega$ follows at once from the definitions; the remaining claim are an immediate consequence of the identity
$({\rm{ad}}a)^{(n)}(b)=\sum_{r+s=n}(-1)^sa^{(r)}ba^{(s)}$.
\end{proof}
\end{myremark}

\begin{myremark}\label{emgfg}
Recalling the embedding $F:\hu\to\tu$ defined in remark \ref{emgg}, theorem \ref{trm} implies that the $\Z$-subalgebra of $\tu$ generated by the divided powers of the $y_{2r+1}^{\pm}$'s is the tensor product of the $\Z$-subalgebras 
$\Z^{(div)}[y_{2r+1}^{\pm}|r\in\Z]$, $\Z^{(sym)}[\k_{\pm r}|r>0]$, $\Z^{(bin)}[\k_0-\tilde c, 2\tilde c]$.

\end{myremark}

\n Mitzman completely described the integral form generated by the divided powers of the Kac-Moody generators in all the twisted cases; in case $A_2^{(2)}$ his result can be stated as follows, using our notations (see examples \ref{dvdpw}, \ref{binex} and \ref{rvsf}, definition \ref{bun} and notation \ref{mitNota}):

\begin{mytheorem}\label{mitz}
$\um\cong\um^-\otimes_{\Z}\um^0\otimes_{\Z}\um^+$ where
$$\um^{\pm}\cong\Z^{(div)}[x_{2r}^{\pm}|r\in\Z]\otimes_{\Z}\Z^{(div)}[y_{2r+1}^{\pm}|r\in\Z]\otimes_{\Z}\Z^{(div)}[x_{2r+1}^{\pm}|r\in\Z]\cong$$
$$\cong\Z^{(div)}[x_{2r+1}^{\pm}|r\in\Z]\otimes_{\Z}\Z^{(div)}[y_{2r+1}^{\pm}|r\in\Z]\otimes_{\Z}\Z^{(div)}[x_{2r}^{\pm}|r\in\Z],$$
$$\um^0\cong\Z_{\lambda}[\k_{-r}|r>0]\otimes_{\Z}\Z^{(bin)}[2\k_0,\tilde c-\k_0]\otimes_{\Z}\Z_{\lambda}[\k_r|r>0].$$
The isomorphisms are all induced by the product in $\tu$.

Remark that $\Z^{(bin)}[2\k_0,\tilde c-\k_0]=\Z^{(bin)}[\k_0-\tilde c,2\tilde c]$ (see example \ref{binex}) and $\Z_{\lambda}[\k_r|r>0]=\Z^{(sym)}[\k_r|r>0]$
(see theorem \ref{gmrvs}).
\end{mytheorem}

\vskip.3 truecm
\begin{myremark} \label{tmnv}
\n As in the case of $\hat\goths$ (see remark \ref{tmfv}) we can evidentiate the relation between the
elements $\hat \k_k$'s with $k>0$ and the elements $p_{n,1}$'s ($n>0$) defined in \cite {F} following Garland's $\Lambda_{k}$'s.

\n Setting 
$$\sum_{n\geq 0}p_nu^n=P(u)=\hat \k(-u)^{-1}$$ we have on one hand
$\Z[\hat \k_k|k>0]=\Z[p_{n}|n>0]$ and on the other hand $$p_0=1,\ \ p_n={1\over n}\sum_{r=1}^n\k_rp_{n-r}\ \forall n>0,$$
hence $p_n=p_{n,1}$ $\forall n\geq 0$ (see \cite{F}) and $\Z[\hat \k_k|k>0]=\Z[p_{n,1}|n>0]$.

\end{myremark}

\begin{mycorollary}\label{yins}
$\tuz\subsetneq\um$.

\n More precisely:
$$\Z^{(div)}[X_{2r+1}^{\pm}|r\in\Z]\subsetneq\Z^{(div)}[y_{2r+1}^{\pm}|r\in\Z],$$ 
so that $\tuz^+\subsetneq\um^+$ and $\tuz^-\subsetneq\um^-$;
$$\Z^{(bin)}[h_0,c]=\Z^{(bin)}[2\k_0,4\tilde c]\subsetneq\Z^{(bin)}[2\k_0,\tilde c-\k_0]$$
and (see definition \ref{thuz})
$$\Z^{(sym)}[\varepsilon_rh_r|r>0]\subsetneq\Z^{(sym)}[\k_r|r>0]$$
(and similarly for the negative part of $\um^0$), so that $\tuz^0\subsetneq\um^0$.
\begin{proof}
For $\Z^{(div)}$ and $\Z^{(bin)}$ the claim is obvious. 
For $\Z^{(sym)}$ the inequality follows at once from the fact that $\k_1={h_1\over 2}$ does not belong to $\Z^{(sym)}[\varepsilon_rh_r|r>0]$ while
the inclusion follows from propositions \ref{convoluzioneintera} and \ref{emmepiallaerre} remarking that for all $r>0$ $\varepsilon_rh_r=2\varepsilon_r\k_r$.

\n Then the assertion for $\tuz$ and $\um$ follows from theorems \ref{trmA22} and \ref{mitz}.
\end{proof}
\end{mycorollary}

\begin{myremark} \label{mizRem}
Theorem \ref{mitz} can be deduced from the commutation formulas discussed in this paper and collected in appendix \ref{appendA}, thanks to the triangular decompositions (see remark \ref{tefp}) 
and to the following observations:

\n 
i) $\um^0$ is a $\Z$-subalgebra of $\tu$: 

\n indeed, since the map $h_r\mapsto \k_r$, $c\mapsto\tilde c$ defines an automorphism of $\tu^0$, proposition \ref{zkd} implies that
$$\hat \k_+(u)\hat \k_-(v)=\hat \k_-(v)(1-uv)^{-4\tilde c}(1+uv)^{2\tilde c}\hat \k_+(u).$$

\n ii) $\um^+$  and $\um^-$ are $\Z$-subalgebras of $\tu$: 

\n indeed the $[(x_{2r}^+)^{(k)},(x_{2s+1}^+)^{(l)}]$'s (the only non trivial commutators in $\um^+$) lie in $\tuz^+\subseteq\um^+$; on the other hand $\um^-=\Omega(\um^+)$. 

\n iii) $\exp\left(\sum_{r>0}a_rx_r^+\right)\in\um^+$ if $a_r\in\Z$ for all $r>0$:

\n see lemma \ref{cle},viii), formula  (\ref{sdivmv}) and the relation $[x_{2r}^+,x_{2s+1}^+]=-4y_{2r+2s+1}^+$.

\n iv) $\um^0\um^+$ and $\um^-\um^0$ are $\Z$-subalgebras of $\tu$:

\n that $(y_{2r+1}^+)^{(k)}\um^0\subseteq\um^0\um^+$ follows from remark \ref{emgfg}; moreover by propositions \ref{bdm} and \ref{hh} we get
$$(x_r^+)^{(k)}{\k_0-\tilde c\choose l}={\k_0-\tilde c-k\choose l}(x_r^+)^{(l)},$$
$$(x_r^+)^{(k)}\hat \k_+(u)=\hat \k_+(u)\left({1-uT^{-1}\over(1+uT^{-1})^2}x_r^+\right)^{(k)},$$
$$\lambda_{-1}(x_r^+)=x_{-r}^+,\ \ 
\lambda_{-1}(\hat \k_+(u))=\hat \k_-(u).$$
On the other hand $\um^-\um^0=\Omega(\um^0\um^+)$.

\n v) $\um^-\um^0\um^+$ is a $\Z$-subalgebra of $\tu$:
$$(x_r^+)^{(k)}(x_s^-)^{(l)}\in\tuz=\tuz^-\tuz^0\tuz^+\subseteq\um^-\um^0\um^+$$
(see theorem \ref{trmA22} and corollary \ref{yins}),
$$(y_{2r+1}^+)^{(k)}(y_{2s+1}^-)^{(l)}\in\um^-\um^0\um^+$$
(see remark \ref{emgfg}), and
$$\exp(x_0^+u)\exp(y_1^-v)=$$
$$=\exp({\alpha_-})\exp({\beta_-})\exp({\gamma_-})\hat \k_+(u^2v)
\exp(\gamma_+)\exp(\beta_+)\exp(\alpha_+)$$
where
$$\alpha_-={uv\over 1-w^2u^4v^2}.x_1^-,\ \ \ \ \alpha_+={u\over 1-w^2u^4v^2}.x_0^+,$$
$$\beta_-={(1+3\cdot wu^4v^2)v\over (1+wu^4v^2)^2}.y_1^-,\ \ \ \ \beta_+={(1-wu^4v^2)u^4v\over (1+wu^4v^2)^2}.y_1^+,$$
$$\gamma_-={-w^2u^3v^2\over 1-w^2u^4v^2}.x_0^-,\ \ \ \ \gamma_+={-u^3v\over 1-w^2u^4v^2}.x_1^+$$
\n (see proposition \ref{xmenogrande} recalling definition \ref{qwmodulo} and remark \ref{whtilde}), so that 
$(x_0^+)^{(k)}(y_1^-)^{(l)}$ lies in $\um^-\um^0\um^+$ for all $k,l\geq 0$; from this it follows that 
$(x_r^+)^{(k)}(y_{2s+1}^-)^{(l)}$ and $(y_{2s+1}^+)^{(l)}(x_r^-)^{(k)}$ lie in $\um^-\um^0\um^+$ for all $r,s\in\Z$, $k,l\geq 0$ because
$\um^-\um^0\um^+$ is stable under $T^{\pm 1}$, $\lambda _m$ ($m\in\Z$ odd) and $\Omega$, and
$$x_r^+=T^{-r}\lambda_{2r+2s+1}(x_0^+),\ \ y_{2s+1}^-=(-1)^rT^{-r}\lambda_{2r+2s+1}(y_1^-),$$
$$y_{2s+1}^+=\Omega(y_{-2s-1}^-),\ \ x_r^-=\Omega(x_{-r}^+);$$

\n vi) $\um\subseteq\um^-\um^0\um^+
$:

\n it follows from v) since $(x_0^{\pm})^{(k)}\in\Z^{(div)}[x_{2r}^{\pm}|r\!\in\!\Z]$ and $(y_{\mp 1}^{\pm})^{(k)}\in\Z^{(div)}[y_{2r+1}^{\pm}|r\!\in\!\Z]$.

\n vii) $\um^{\pm}\subseteq\um$: 

\n this follows from remark \ref{stau}, observing that
$$\tau_0(x_r^+)=(-1)^{r-1}x_{r+1}^-,\ \ \tau_1(x_r^-)=x_r^+,\ \ \tau_1(y_{2r+1}^-)=y_{2r+1}^+,\ \ \tau_0(y_{2r+1}^+)=-y_{2r+3}^-.$$

\n viii) $\um^0\subseteq\um$:

\n it follows from vii), v) and the stability under $\Omega$.

\n ix) $\um^-\um^0\um^+\subseteq\um$:

\n this is just vii) and viii) together.

Then $\um=\um^-\um^0\um^+$, which is the claim.

\end{myremark}

\begin{myremark}\label{rfink}
As one can see from remark \ref{mizRem},vii), $$\{x_r^{\pm},y_{2r+1}^{\pm},\k_s,2\k_0,\tilde c-\k_0|r,s\in\Z, s\neq 0\}$$
is, up to signs, a Chevalley basis of $\gothd$ (see \cite{DM}).

\n It is actually through these basis elements that Mitzman introduces, following \cite{HG}, the integral form of $\tu$, as the $\Z$-subalgebra of $\tu$ generated by
$$\{(x_r^{\pm})^{(k)},(y_{2r+1}^{\pm})^{(k)}|r\in\Z,k\in\N\};$$
but this $\Z$-subalgebra is precisely the algebra $\um$ introduced in definition \ref{mif}: indeed
it turns out to be generated over $\Z$ just by $\{e_i^{(k)},f_i^{(k)}|i=0,1,\ k\geq 0\}$, that is by $\{(x_0^{\pm})^{(k)},(y_{\mp 1}^{\pm})^{(k)}|k\geq 0\}$, thanks to remarks \ref{stau} and \ref{mizRem},vii).

\end{myremark}

\subsection{List of Symbols}\label{appendD}

\underline{Lie Algebras and Commutative Algebras}:
	\begin{abbrv}
	\item[$S^{(div)}$]  \hfill Example \ref{dvdpw} \\
	\item[$S^{(bin)}$]  \hfill Example \ref{binex} \\
	\item[$S^{(sym)}$]  \hfill Example \ref{rvsf} \\
	\item[$\frak {sl_2}$]  \hfill  Definition \ref{sl2} \\
	\item[$\hat{\goths}$]  \hfill  Definition \ref{hs2} \\
	\item[$\gothd$] \hfill Definition \ref{a22} \\
	\end{abbrv}
	
\underline{Enveloping Algebras}:
	\begin{abbrv}
	\item[$  \uz ^{re, \pm}, \;\uz^{im,\pm}, \; \uz^{\frak h},\; \duz, \; \duz^{imm, \pm}$]      \hfill Section \ref{intr} \\
	\item[$\u, \; \uz$]  \hfill Definiton \ref{sl2} \\
	\item[$\u^+,\; \u^-,\; \u^0$] \hfill Theorem \ref{trdc} \\
	\item[$\hu,\; \hu^{+},\;  \hu^{-}, \; \hu^{0}, \; \hu^{0,\pm},\;  \hu^{0,0}$ ]     \hfill Definition \ref{hs2} \\
	\item[$\huz,\; \huz^\pm, \;  \huz^{0,\pm}, \; \huz^{0,0}$ ]     \hfill Definition \ref{hhuz} \\
	\item[$\tu,\; \tu^{\pm},\; \tu^{0},\; \tu^{\pm,0},\; \tu^{\pm,1},\; \tu^{\pm,c}, \;  \tu^{0,\pm},\; \tu^{0,0}$] \hfill Definition \ref{a22} \\
	\item[$\tuz,\; \tuz^{\pm},\; \tuz^{0},\; \tuz^{\pm,0},\; \tuz^{\pm,1},\; \tuz^{\pm,c} \;  \tuz^{0,\pm},\; \tuz^{0,0}$] \hfill Definition \ref{thuz}  \\

	\item[$\um$]   \hfill  Definiton \ref{mif}\\
	\item[$\um^-, \;\um^0, \;\um^+ \;$]   \hfill  Theorem \ref{mitz}\\
	\end{abbrv}
	
\underline{Bases}:
	\begin{abbrv}
	\item[$B^{re,\pm}, \; B^{im,\pm}, \;B^{\frak h}$]      \hfill Section \ref{intr} 
	\item[$B^\pm,\; B^{0, \pm},\; B^{0,0} $ ]     \hfill Theorem \ref{trm} \\
	\item[$B^{\pm,0},\; B^{\pm,1},\; B^{\pm,c} $ ]     \hfill Theorem \ref{trmA22} \\
	\item[$B_\lambda \; B_\lambda^{[n]} \; B_x \;B_x^{[n]}$] \hfill Definitions \ref{bun} and \ref{basi} \\
	\end{abbrv}

\underline{Elements and their generating series}:
	\begin{abbrv}
	\item[$\Lambda_r(\xi(k))$]      \hfill Section \ref{intr} 
	\item[$a^{(k)}, \; \exp(au)$]  \hfill Notation \ref{ntdvd} \\
	\item[$\binom{a}{k}, \; (1+u)^a$]  \hfill Notation \ref{ntbin} \\
	\item[$\hat p(u),\; \hat p_r$]  \hfill Example \ref{rvsf} \\
	\item[$\hat h_r^{\{a\}},\; \hat h_+^{\{a\}}(u)$]  \hfill Notation \ref{hcappucciof} \\
	\item[$x^\pm_r, \; h_r, \;c$]     \hfill Definition \ref{hs2}  and Definition \ref{a22}\\
	\item[$X^\pm_{2r+1}$]     \hfill Definition \ref{a22}\\
	\item[$x^\pm(u),\; h_\pm(u), \; \hat h_\pm(u), \; \hat h_r$]     \hfill Notation \ref{hgens} \\
	\item[$\tilde h_{\pm}(u), \tilde h_{\pm r}$] \hfill Definition \ref{thuz} \\
	\item[$e_i, f_i, h_i$] \hfill Remark \ref{ka22}\\
	\item[$y_{2r+1}^\pm , \; \k_r, \; \tilde c$]  \hfill  Notation \ref{mitNota}\\
	\item[$\k_\pm(u)$] \hfill Remark \ref{mizRem}\\
	\end{abbrv}
	
\underline{Anti/auto/homomorphisms}:
	\begin{abbrv}
	\item[$\lambda_m, \; \lambda_m ^{[n]}$]  \hfill Proposition \ref{tmom}  \\	
	\item[$ev$]  \hfill Equation \ref{evaluation}  \\	
	\item[$\sigma, \; \Omega,\; T, \; \lambda_m$] \hfill Definition \ref{hto} and Definition \ref{tto} \\
	\item[$\tilde \lambda_m$] \hfill Lemma \ref{ometiomecap} \\
	\end{abbrv}

\underline{Other symbols}:
	\begin{abbrv}
	\item[$ \uno,  \; \uno^{(m)},  \; \uno_r, \; \uno^{(m)}_r$]  \hfill Notation \ref{hcappucciof} \\
	\item[$L_a, \; R_a$]  \hfill Notation \ref{lard} \\
	\item[$\varepsilon_r$] \hfill Definition \ref{thuz} \\
	\item[$L, \; L^\pm,\; L^0, \; L^{\pm, 0}, \; L^{\pm, 1}, \; L^{\pm, c}$] \hfill Definition \ref{sottoalgebraL} \\
	\item[$w.$] \hfill Definition \ref{qwmodulo} \\
	\item[$d, \; \tilde d, \; d_n, \; \tilde d_n$] \hfill Notation \ref{notedn} \\
	\item[$\delta_n$] \hfill Remark \ref{hhdehh} \\
	\end{abbrv}

\end{document}